\documentclass[12pt,reqno]{amsart}
\usepackage[margin=1in]{geometry}
\usepackage{amsmath,amssymb,amsthm,graphicx,amsxtra, setspace}
\usepackage[utf8]{inputenc}
\usepackage{mathrsfs}
\usepackage{hyperref}
\usepackage{upgreek}
\usepackage{mathtools}
\usepackage[mathcal]{euscript}
\usepackage{xcolor}
\allowdisplaybreaks

\usepackage[pagewise]{lineno}

\newtheorem{theorem}{Theorem}[section]
\newtheorem{lemma}[theorem]{Lemma}

\newtheorem{assumption}[theorem]{Assumption}
\newtheorem{corollary}[theorem]{Corollary}
\newtheorem{definition}[theorem]{Definition}
\newtheorem{example}[theorem]{Example}
\newtheorem{remark}[theorem]{Remark}

\let\originalleft\left
\let\originalright\right
\renewcommand{\left}{\mathopen{}\mathclose\bgroup\originalleft}
\renewcommand{\right}{\aftergroup\egroup\originalright}

\renewcommand{\d}{\/\mathrm{d}\/}

\def\w{\textbf{W}^{\varepsilon}_{{\theta}^{\varepsilon}}}

\def\e{\varepsilon}

\def\L{\mathbb{L}}
\def\A{\mathrm{A}}

\def\C{\mathrm{C}}
\def\f{\boldsymbol{f}}

\def\B{\mathrm{B}}
\def\D{\mathrm{D}}
\def\y{\boldsymbol{y}}

\def\x{\boldsymbol{x}}

\def\g{\boldsymbol{g}}
\def\p{\boldsymbol{p}}
\def\h{\boldsymbol{h}}
\def\z{\boldsymbol{z} }
\def\v{\boldsymbol{v}}
\def\w{\boldsymbol{w}}
\def\W{\mathrm{W}}

\def\N{\mathbb{N}}

\def\V{\mathbb{V}}
\def\wi{\widetilde}

\def\u{\mathrm{U}}

\def\u{\boldsymbol{u}}
\def\H{\mathbb{H}}

\newcommand{\R}{\mathbb{R}}

\renewcommand{\d}{\/\mathrm{d}\/}

\newcommand{\Addresses}{{
		\footnote{
			\noindent \textsuperscript{1,2}Department of Mathematics, Indian Institute of Technology Roorkee-IIT Roorkee,
			Haridwar Highway, Roorkee, Uttarakhand 247667, INDIA.\par\nopagebreak
			\noindent  \textit{e-mail:} \texttt{Manil T. Mohan: maniltmohan@ma.iitr.ac.in, maniltmohan@gmail.com.}
			
			\textit{e-mail:} \texttt{Kush Kinra: kkinra@ma.iitr.ac.in.}
			
			\noindent \textsuperscript{*}Corresponding author.
			
			\textit{Key words:} Wong-Zakai approximation, random pullback attractor, upper semicontinuity, stochastic convective Brinkman-Forchheimer equations, bounded domains, whole domain.
			
			Mathematics Subject Classification (2020): Primary 35B41, 35R60, 35Q35; Secondary 37L55, 37N10.

}}}

\begin{document}
	
	\title[Attractors for random CBF equations with colored noise]{Long term behavior of 2D and 3D non-autonomous random convective Brinkman-Forchheimer equations driven by colored noise
		\Addresses}
	
	\author[K. Kinra and M. T. Mohan]
	{Kush Kinra\textsuperscript{1} and Manil T. Mohan\textsuperscript{2*}}

	\maketitle
	
	\begin{abstract}
		The long time behavior of Wong-Zakai approximations of 2D as well as 3D non-autonomous stochastic convective Brinkman-Forchheimer (CBF) equations with non-linear diffusion terms on bounded and unbounded ($\mathbb{R}^d$ for $d=2,3$) domains is discussed in this work. To establish the existence of random pullback attractors, the concept of asymptotic compactness (AC) is used. In bounded domains, AC is proved via compact Sobolev embeddings. In unbounded domains, due to the lack of compact embeddings, the ideas of energy equations and uniform tail estimates are exploited to prove AC. In the literature, CBF equations are also known as \emph{Navier-Stokes equations (NSE) with damping}, and it is interesting to see that the modification in NSE by linear and nonlinear damping provides  better results than that available for NSE. The presence of linear damping term helps to establish the  results in the whole domain $\R^d$. The nonlinear damping term supports to obtain better results in  3D and also for a large class of nonlinear diffusion terms. Moreover, we prove the existence of a unique random pullback attractor for stochastic CBF equations with additive white noise. Finally, for additive as well as multiplicative noise case, we establish the convergence of solutions and upper semicontinuity of random pullback attractors for Wong-Zakai approximations of stochastic CBF equations towards the random pullback attractors for stochastic CBF equations when correlation time of colored noise converges to zero.  
	\end{abstract}

	\section{Introduction} \label{sec1}\setcounter{equation}{0}
	In this article, we are interested in the long time behavior of   stochastic convective Brinkman-Forchheimer (CBF) equations with non-autonomous forcing term (deterministic term) and nonlinear diffusion term. This stochastic mathematical model describes  the motion of incompressible fluid flows in a saturated porous medium (cf. \cite{PAM}).   The CBF equations apply to flows when the velocities are sufficiently high and porosities are not too small, that is, when the Darcy law for a porous medium does not hold. Due to the lack of Darcy's law, it is sometimes referred as \emph{non-Darcy model} also (cf. \cite{PAM}). Given $\mathfrak{s}\in\R$, we consider the non-autonomous stochastic CBF equations in $\mathcal{O}\subseteq\mathbb{R}^d\ (d=2,3)$ as
	\begin{equation}\label{1}
		\left\{
		\begin{aligned}
			\frac{\partial \u}{\partial t}-\mu \Delta\u+(\u\cdot\nabla)\u+\alpha\u+\beta|\u|^{r-1}\u+\nabla p&=\boldsymbol{f}+S(t,x,\u)\circ\frac{\d \W}{\d t}, \\&\qquad \qquad \ \ \ \text{ in } \mathcal{O}\times(\mathfrak{s},\infty), \\ \nabla\cdot\u&=0, \qquad \ \  \text{ in } \ \mathcal{O}\times(\mathfrak{s},\infty), \\
			\u(x,\mathfrak{s})&=\u_{\mathfrak{s}}(x),\  \ \ x\in \mathcal{O} \text{ and }\mathfrak{s}\in\R,
		\end{aligned}
		\right.
	\end{equation}
	where $\u(x,t) :
	\mathcal{O}\times(\mathfrak{s},\infty)\to \mathbb{R}^d$ denotes the velocity field, $p(x,t):
	\mathcal{O}\times(\mathfrak{s},\infty)\to\R$ represents the pressure field, $\f(x,t):
	\mathcal{O}\times(\mathfrak{s},\infty)\to \mathbb{R}^d$ is an external body force, $S(t,x,\u)$ is a nonlinear diffusion term (see sections \ref{sec3} and \ref{sec4}), the symbol $\circ$ means that the stochastic integral should be understood in the sense of Stratonovich integral and $\W=\W(t,\omega)$ is an one-dimensional two-sided Wiener process defined on a probability space $(\Omega,\mathscr{F},\mathbb{P})$. Here  $\Omega$ is given  by
	\begin{align*}
		\Omega=\{\omega\in\mathrm{C}(\R;\R):\omega(0)=0\}, 
	\end{align*}
	$\mathscr{F}$ is the Borel sigma-algebra induced by the compact-open topology of $\Omega$, and $\mathbb{P}$ is the two-sided Gaussian measure on $(\Omega,\mathscr{F})$. Consequently, $\W$ has a form $\W(t,\omega)=\omega(t)$. Moreover, $\u(\cdot,\cdot)$ satisfies $$\u=\mathbf{0}\ \text{ on  }\  \partial\mathcal{O}\times(\mathfrak{s},\infty)\ \text{ or }\ \u(x,\mathfrak{s})\to 0\  \ \text{ as }\  \  |x|\to\infty,$$ when $\mathcal{O}\subset\R^d$ is bounded with smooth boundary or  $\mathcal{O}=\R^d$. The positive constants $\mu,\alpha, \beta>0$ are the Brinkman (effective viscosity), Darcy (permeability of porous medium) and Forchheimer (proportional to the porosity of the material) coefficients, respectively. The absorption exponent $r\in[1,\infty)$, $r=3$ is called the critical exponent and the model \eqref{1} with $r>3$ is referred as the CBF equations with fast growing nonlinearites (\cite{KT2}). The system \eqref{1} is also known as damped Navier-Stokes equations (NSE) (cf. \cite{HZ}) because it is classical NSE for $\alpha=\beta=0$. Moreover, it is proved in \cite{HR} (see Proposition 1.1, \cite{HR}) that CBF equations \eqref{1} have the same scaling as NSE only when $r=3$ and $\alpha=0$ but no scale invariance property for other values of $\alpha$ and $r$, therefore it is sometimes called NSE modified by an absorption term (\cite{SNA}) or  tamed NSE (\cite{MRXZ}).

	Let us now discuss some solvability results available in the literature  for the stochastic system \eqref{1} and similar models. There are a good number of works available on the solvability results for the system \eqref{1} and related models (see \cite{MRXZ1,WLMR,WL,HBAM,MTM1,LHGH1}, etc). In particular, 3D NSE with a Brinkman-Forchheimer type term subject to an anisotropic viscosity driven by multiplicative noise and stochastic tamed 3D NSE on whole domain is considered in \cite{HBAM} and \cite{MRXZ1}, respectively. An improvement of the work \cite{MRXZ1} for a slightly simplified system is done in \cite{ZBGD}. On bounded domains, the existence of martingale solutions and strong solutions (in the probabilistic sense) for stochastic CBF equations are established in \cite{LHGH1} and \cite{MTM1}, respectively. The existence of a unique pathwise strong solution to 3D stochastic NSE is a well known open problem, and the same is open for 3D stochastic CBF equations with  $r\in[1,3)$. Therefore, we do not consider $d=3$ with $1\leq r<3$ for any $\mu, \beta>0$ and $d=r=3$ for $2\beta\mu<1$ (see \cite{MTM1}). In the sequel, if we write $r\in[1,3),$ then it indicates that the two dimensional case is considered.
	
	In \cite{BCF,CDF}, authors introduced the concept of random attractors for stochastic partial differential equations (SPDEs) and applied it to some fluid flow models including 2D stochastic NSE. Later, this concept is used to establish the existence of random attractors for several fluid flow models (see \cite{FY,KM3,You} etc and references therein). In unbounded domains, due to the absence of compact Sobolev embeddings, random dynamical systems (RDS)  are no more compact. In the deterministic case, the difficulty discussed above was resolved by different methods, cf.  \cite{Abergel, Ghidaglia,MTM3,Rosa}, etc for the autonomous case and \cite{CLR1, CLR2,KM6}, etc for non-autonomous case. For SPDEs, the methods available in the deterministic case have also been generalized by several authors (see for example, \cite{ BLL,BLW,BCLLLR,Wang}, etc).  The concept of an asymptotically compact cocycle was introduced in \cite{CLR1} and the authors proved the existence of attractors for  non-autonomous 2D NSE. Later, several authors used this method to prove the existence of random attractors on unbounded domains, see for example \cite{BLL,BLW,KM,LXS,Wang} etc. Two of the methods used to prove asymptotic compactness without Sobolev embeddings are as follows: the method of energy equations (see \cite{GGW,KM,PeriodicWang} etc) and the method of uniform tail estimates (see \cite{SandN_Wang,WLW} etc).

	For the existence of a unique random attractor (dynamics of almost all sample paths), it is required to define RDS  or random cocycle through the solution operator of the equation \eqref{1}. To the best our knowledge, it is still an open problem to define RDS for \eqref{1}, when $S(t,x,\u)$ is a general nonlinear diffusion term. In this direction, when $S(t,x,\u)$ is a nonlinear diffusion term, the concept of weak pullback mean random attractors was introduced in \cite{Wang1} and this theory was applied to several physical models (cf. \cite{GU,KM4,Wang2} etc). On the other hand, using random transformations, SPDEs can generate RDS for either $S(t,x,\u)=\u$ or $S(t,x,\u)$ is independent of $\u$. Indeed, the existence of a unique random attractor for \eqref{1} is available for such  $S(t,x,\u)$ only (see \cite{KM1,KM3} for bounded domains and \cite{KM,KM6} for unbounded domains). 
	
	Wong and Zakai introduced the concept of approximating stochastic differential equations (SDEs) by deterministic differential equations in \cite{WZ,WZ1}, where they have approximated one dimensional Brownian motion. Later, this work was extended to SDEs in higher dimensions, see for example \cite{IW,DK.IM,Sussmann} etc and references therein. Due to the lack of differentiability of sample paths of Wiener process $\W$, the authors in \cite{UO,WU} first introduced the Ornstein-Uhlenbeck process (or colored noise) to approximate $\W$. The Wong-Zakai approximations (with the help of colored noise)  can be used to analyze the pathwise random dynamics of \eqref{1} with nonlinear diffusion term. The approach through approximations obeys the well-known framework given by Wong and Zakai in \cite{WZ,WZ1}. Despite of analyzing the corresponding random versions of \eqref{1} acquired  by random transformations, we will rather analyze the system obtained by the Wong-Zakai approximations of \eqref{1} with nonlinear diffusion term. With reference to the nonlinear diffusion term, the idea of Wong-Zakai approximations for random pullback attractors is used for a wide class of physically admissible models, see \cite{GLW,GW,LW} etc for bounded domains, \cite{GGW,WLS,WLW} etc for unbounded domains and \cite{WSLW,YSW} etc for discrete systems.
	
	On  $(\Omega,\mathscr{F},\mathbb{P})$, consider Wiener shift operator $\{\vartheta_{t}\}_{t\in\R}$ defined by 
	\begin{equation}\label{vartheta}
		\vartheta_{t}\omega(\cdot)=\omega(\cdot+t) -\omega(t), \ \ \   t\in\R\ \text{ and }\ \omega\in\Omega.
	\end{equation}
	It is well known from \cite{Arnold} that the Gaussian measure $\mathbb{P}$ is ergodic and invariant for $\vartheta_{t}$. Hence, $(\Omega,\mathscr{F},\mathbb{P},\{\vartheta_{t}\}_{t\in\R})$ is a metric dynamical system (see \cite{Arnold}). Let us define the colored noise $\mathcal{Z}_{\delta}:\Omega\to\R$ with the given correlation time $\delta\in\R\backslash\{0\}$ such that 
	\begin{align*}
		\mathcal{Z}_{\delta}(\omega)=\frac{\omega(\delta)}{\delta} \ \text{\ or equivalantly  }\	\mathcal{Z}_{\delta}(\vartheta_{t}\omega)=\frac{1}{\delta}\left(\omega(t+\delta)-\omega(t)\right).
	\end{align*} 
	The properties of Wiener process  imply that $\mathcal{Z}_{\delta}(\vartheta_{t}\omega)$ is stationary process which follows normal distribution. It is observed that white noise can be approximated by $\mathcal{Z}_{\delta}(\vartheta_{t}\omega)$ in the sense given in Lemma \ref{WnCn} below (see Lemma 2.1, \cite{GGW} also). In \cite{LW1,SLW}, the authors used this approximation to study the chaotic behavior of RDS. 
	
	Let us now consider the following random partial differential equations as Wong-Zakai approximations of \eqref{1} for $\mathfrak{s}\in\R$ and $\delta\in\R\backslash\{0\}$:
	\begin{equation}\label{2}
		\left\{
		\begin{aligned}
			\frac{\partial \u}{\partial t}-\mu \Delta\u+(\u\cdot\nabla)\u+\alpha\u+\beta|\u|^{r-1}\u+\nabla p&=\boldsymbol{f}+S(t,x,\u)\mathcal{Z}_{\delta}(\vartheta_{t}\omega), \ \text{ in } \mathcal{O}\times(\mathfrak{s},\infty), \\ \nabla\cdot\u&=0, \ \ \ \ \ \ \ \ \ \ \ \ \ \ \ \ \ \ \ \ \ \ \ \ \ \ \ \text{ in } \ \mathcal{O}\times(\mathfrak{s},\infty), \\
			\u(x,\mathfrak{s})&=\u_{\mathfrak{s}}(x),\ \ \ \ \ \ \ \ \ \ \ \ \ \ \ \ \ \ x\in \mathcal{O} \text{ and }\mathfrak{s}\in\R.
		\end{aligned}
		\right.
	\end{equation}
	Also, $\u(x,t)$ satisfies $\u=\mathbf{0}$ on   $\partial\mathcal{O}\times(\mathfrak{s},\infty)$ or  $\u(x,\mathfrak{s})\to 0$ as  $|x|\to\infty,$ when $\mathcal{O}\subset\R^d$ is bounded with smooth boundary  or $\mathcal{O}=\R^d$. In bounded domains, as we are using compact Sobolev embeddings, we assume that the nonlinear diffusion term $S(t,x,\u)$ satisfies Assumption \ref{NDT1}. Due to the absence of compact Sobolev embeddings in unbounded domains, we use the method of energy equations and method of uniform tail estimates, which require a stronger assumption than Assumption \ref{NDT1} on $S(t,x,\u)$ (see Assumptions \ref{NDT2} and \ref{NDT3} for the method of energy equations and method of uniform tail estimates, respectively). It can be verified that Assumptions \ref{NDT2} and \ref{NDT3} are not same. In particular, both Assumptions cover different classes of nonlinear functions. For functional setting, one can refer to section \ref{sec2}. 
	\vskip 2mm
	\noindent
	\textbf{Aims and scopes of the work:}   The major aims and novelties of this work are:
	\begin{itemize}
		\item [(i)] Existence of a unique random pullback attractor for the system \eqref{2} under Assumption \ref{NDT1} on bounded domains.
		\item [(ii)] Existence of a unique random pullback attractor for  the system \eqref{2} under Assumptions \ref{NDT2} and \ref{NDT3} on the whole domain.
		\item [(iii)] Existence of a unique random pullback attractor for the system \eqref{1} with $S(t,x,\u)=e^{\sigma t}\textbf{g}(x)$, for given $\sigma\geq0$ and $\textbf{g}\in\D(\A)$ (additive noise).
		\item [(iv)] For $S(t,x,\u)=e^{\sigma t}\textbf{g}(x)$, for given $\sigma\geq0$ and $\textbf{g}\in\D(\A)$ as well as $S(t,x,\u)=\u$, convergence of solutions and upper semicontinuity of random pullback attractors for \eqref{2} towards the random pullback attractor for \eqref{1} as $\delta\to0$ (additive and multiplicative noise).
	\end{itemize}
	Moreover, 	we prove the existence of a unique random pullback attractor of the system \eqref{2} in two dimensions  with  $\alpha=\beta=0$ (2D NSE) under Assumption \ref{NDT4} on Poincar\'e domains.
	\vskip 2mm
	\noindent
	\textbf{Advantages of damping term:} CBF equations are also known as NSE with damping (cf. \cite{HZ}). The damping arises from the resistance to the motion of the flow or by  friction effects. Due to the presence of the damping term $\alpha\u+\beta|\u|^{r-1}\u$, we can establish better results than which are available for NSE. In particular, in this work, linear damping $\alpha\u$ helps us to obtain results on the whole domain (for NSE, the global as well as random attractor on the whole domain is still an open problem) and the nonlinear damping $\beta|\u|^{r-1}\u$ (for $r>1$) helps us to cover a large class of nonlinearity on $S(t,x,\u)$ (see subsection \ref{subsec4.3}). The authors in \cite{GGW} established the existence of a unique random attractor for 2D NSE with $S(t,x,\u)=e^{\sigma t}\textbf{g}(x)$, for given $\sigma>0$ and $\textbf{g}\in\D(\A)$,  based on the following assumption on $\textbf{g}(\cdot)$:
	\begin{assumption}[See section 3, \cite{GGW}]\label{GA}
		Assume that the function \emph{$\textbf{g}(\cdot)$} satisfies the following condition: there exists a strictly positive constant $\aleph$ such that 
		\begin{align*}
			|b(\u,\emph{\textbf{g}},\u)|\leq \aleph\|\u\|^2_{\H}, \ \  \text{ for all } \ \u\in\H.
		\end{align*}
	\end{assumption}
	We also point out that the requirement of Assumption \ref{GA} is not necessary for $r>1$ ($\textbf{g}\in\D(\A)$ is enough).	For $r=1$, one gets linear damping, and the existence of  a random attractor for \eqref{SCBF_Add} and upper semicontinuity of random attractors for \eqref{WZ_SCBF_Add} can be proved  on the whole domain via same arguments  as it has been done for 2D stochastic NSE on Poincar\'e domains in \cite{GGW} under Assumption \ref{GA}. For $r>1$, due to the presence of  nonlinear damping term $\beta|\u|^{r-1}\u$, we provide a different treatment to avoid Assumption \ref{GA} on $\textbf{g}$ (see Lemma \ref{LemmaUe_add}). 
	
	It is well-known that passing limit in nonlinear terms is not an easy task in any analysis. In this work, while proving the pullback asymptotic compactness of solutions of  the system \eqref{SCBF_Add} (Lemma \ref{Asymptotic_UB_Add}) or establishing the uniform compactness of random attractors of the system \eqref{WZ_SCBF_Add} (Lemma \ref{precompact}), it is necessary to show the convergence of the nonlinear terms appearing in the corresponding energy equations.   In our model, we have two nonlinear terms given by $(\u\cdot\nabla)\u$ and $\beta|\u|^{r-1}\u$. The required convergence is obtained by breaking down the integral containing  nonlinear terms into two parts such that one part is defined in a bounded domain $\mathcal{O}_k$ with radius $k$, and another part is defined in the complement of $\mathcal{O}_k$. We then show that the nonlinear term is convergent in $\mathcal{O}_k$ and its tail on the complement of $\mathcal{O}_k$ is uniformly small when $k$ is sufficiently large, from which the required convergence follows (see Corollaries \ref{convergence_b} and \ref{convergence_c} for terms $(\u\cdot\nabla)\u$ and  $\beta|\u|^{r-1}\u$, respectively).
	
	The remaining sections are arranged as follows. In next section, we provide the necessary function spaces required to establish the results of this work, and we define linear operator, bilinear operator and nonlinear operator  to obtain an  abstract formulation of the system \eqref{2}. Furthermore, we recall some well known inequalities, properties of operators, and some properties of white and colored noises. Finally, we discuss the solvability of the system \eqref{2} in the same section. In section \ref{sec3}, we prove the existence of a unique pullback random attractor for the system \eqref{2} under Assumption \ref{NDT1} on bounded domains in which we prove the asymptotic compactness using compact Sobolev embeddings (Theorem \ref{WZ_RA_B}). In section \ref{sec4}, we prove the existence of a unique pullback random attractor for the system \eqref{2} under Assumptions \ref{NDT2} and \ref{NDT3} on the whole domain in which we  prove the asymptotic compactness using the idea of energy equations and uniform tail estimates, for Assumptions \ref{NDT2} and \ref{NDT3}, respectively (Theorems \ref{WZ_RA_UB} and \ref{WZ_RA_UB_GS}). Section \ref{sec6} is devoted for the stochastic CBF equations perturbed by additive white noise. Firstly, we prove the existence of a unique pullback random attractor for stochastic CBF equations driven by additive white noise (Theorem \ref{RA_add}). Next, we demonstrate the convergence of solutions and upper semicontinuity of random pullback attractors for Wong-Zakai approximations of stochastic CBF equations towards the solution and random pullback attractor of stochastic CBF equations, respectively, as the correlation time $\delta$ converges to zero, using the fundamental theory introduced in \cite{non-autoUpperWang}  (Theorem \ref{Main_T_add}). Since the existence of a unique random pullback attractor for stochastic CBF equations driven by linear multiplicative noise is established in \cite{KM6}, we prove only the convergence of solutions and upper semicontinuity of random pullback attractors for its Wong-Zakai approximations towards its solution and random pullback attractor, respectively, as the correlation time $\delta\to0$ (Theorem \ref{Main_T_Multi}) in  section \ref{sec7}. In Appendix \ref{sec5}, we prove the existence of a unique pullback random attractor for the system \eqref{WZ_NSE} under Assumption \ref{NDT4} on Poincar\'e domains (bounded or unbounded) (Theorem \ref{WZ_RA_UB_GS_NSE}).

	\section{Mathematical Formulation}\label{sec2}\setcounter{equation}{0}
	In this section, first we provide the necessary function spaces needed to obtain the results of this work. Next, we define some operators to set up an abstract formulation. Finally, we recall some properties of white noise as well as colored noise. We fix  $\mathcal{O}$ as either a bounded subset of $\R^d$ with $\mathrm{C}^2$-boundary or whole domain $\R^d$.
	\subsection{Function spaces} 
	We define the space $$\mathcal{V}:=\{\u\in\C_0^{\infty}(\mathcal{O};\mathbb{R}^d):\nabla\cdot\u=0\},$$ where $\C_0^{\infty}(\mathcal{O};\mathbb{R}^d)$ denotes the space of all infinite times differentiable functions ($\mathbb{R}^d$-valued) with compact support in $\mathbb{R}^d$. Let $\H$, $\V$ and $\wi\L^p$ denote the completion of $\mathcal{V}$ in $\mathrm{L}^2(\mathcal{O};\mathbb{R}^d)$, $\mathrm{H}^1(\mathcal{O};\mathbb{R}^d)$ and $\mathrm{L}^p(\mathcal{O};\mathbb{R}^d)$, $p\in(2,\infty)$,  norms, respectively. The space $\H$ is endowed with the norm $\|\u\|_{\H}^2:=\int_{\mathcal{O}}|\u(x)|^2\d x,$ the norm on the space $\widetilde{\L}^{p}$ is defined by $\|\u\|_{\wi \L^p}^p:=\int_{\mathcal{O}}|\u(x)|^p\d x,$ for $p\in(2,\infty)$ and the norm on the space $\V$ is given by $\|\u\|^2_{\V}=\int_{\mathcal{O}}|\u(x)|^2\d x+\int_{\mathcal{O}}|\nabla\u(x)|^2\d x.$ The inner product in the Hilbert space $\H$ is denoted by $( \cdot, \cdot)$. The duality pairing between the spaces $\V$ and $\V'$, and $\widetilde{\L}^p$ and its dual $\widetilde{\L}^{\frac{p}{p-1}}$ is represented by $\langle\cdot,\cdot\rangle.$ It should be noted that $\H$ can be identified with its own dual $\H'$. We endow the space $\V\cap\widetilde{\L}^{p}$ with the norm $\|\u\|_{\V}+\|\u\|_{\widetilde{\L}^{p}},$ for $\u\in\V\cap\widetilde{\L}^p$ and its dual $\V'+\widetilde{\L}^{p'}$ with the norm $$\inf\left\{\max\left(\|\v_1\|_{\V'},\|\v_1\|_{\widetilde{\L}^{p'}}\right):\v=\v_1+\v_2, \ \v_1\in\V', \ \v_2\in\widetilde{\L}^{p'}\right\}.$$ Moreover, we have the following continuous  embedding also:
	$$\V\cap\widetilde{\L}^{p}\hookrightarrow\V\hookrightarrow\H\equiv\H'\hookrightarrow\V'\hookrightarrow\V'+\widetilde\L^{\frac{p}{p-1}}.$$ One can define equivalent norms on $\V\cap\widetilde\L^{p}$ and $\V'+\widetilde\L^{\frac{p}{p-1}}$ as 
	\begin{align*}
		\|\u\|_{\V\cap\widetilde\L^{p}}=\left(\|\u\|_{\V}^2+\|\u\|_{\widetilde\L^{p}}^2\right)^{\frac{1}{2}}\ \text{ and } \ \|\u\|_{\V'+\widetilde\L^{\frac{p}{p-1}}}=\inf_{\u=\v+\w}\left(\|\v\|_{\V'}^2+\|\w\|_{\widetilde\L^{\frac{p}{p-1}}}^2\right)^{\frac{1}{2}}.
	\end{align*}
	
	\subsection{Projection operator}
	Let $\mathcal{P}_p: \L^p(\mathcal{O}) \to\widetilde{\L}^p$ be the Helmholtz-Hodge (or Leray) projection  (cf.  \cite{DFHM,MTSS}, etc). For $p=2$, $\mathcal{P}:=\mathcal{P}_2$ becomes an orthogonal projection  and for $2<p<\infty$, it is a bounded linear operator.  
	\vskip 2mm
	\noindent
	\textbf{Case I:} \textit{When $\mathcal{O}$ is a bounded subset of $\R^d$ with $\mathrm{C}^2$-boundary.} Since $\mathcal{O}$ has $\mathrm{C}^2$-boundary, $\mathcal{P}$ maps $\H^1(\mathcal{O})$ into itself (see Remark 1.6, \cite{Temam}).
	\vskip 2mm
	\noindent
	\textbf{Case II:} \textit{When $\mathcal{O}=\R^d$.} The projection operator $\mathcal{P}:\L^2(\R^d) \to\H$ can be expressed in terms of the Riesz transform (cf. \cite{MTSS}). Moreover, $\mathcal{P}$ and $\Delta$ commutes, that is, $\mathcal{P}\Delta=\Delta\mathcal{P}$. 
	\subsection{Linear operator}
	We define the Stokes operator
	\begin{equation*}
		\A\u:=-\mathcal{P}\Delta\u,\;\u\in\D(\A):=\V\cap\H^{2}(\mathcal{O}).
	\end{equation*}
	\subsection{Bilinear operator}
	Let us define the \emph{trilinear form} $b(\cdot,\cdot,\cdot):\V\times\V\times\V\to\R$ by $$b(\u,\v,\w)=\int_{\mathbb{R}^d}(\u(x)\cdot\nabla)\v(x)\cdot\w(x)\d x=\sum_{i,j=1}^d\int_{\mathbb{R}^d}\u_i(x)\frac{\partial \v_j(x)}{\partial x_i}\w_j(x)\d x.$$ If $\u, \v$ are such that the linear map $b(\u, \v, \cdot) $ is continuous on $\V$, the corresponding element of $\V'$ is denoted by $\B(\u, \v)$. We also denote $\B(\u) = \B(\u, \u)=\mathcal{P}[(\u\cdot\nabla)\u]$.
	An integration by parts gives 
	\begin{equation}\label{b0}
		\left\{
		\begin{aligned}
			b(\u,\v,\v) &= 0,\ \text{ for all }\ \u,\v \in\V,\\
			b(\u,\v,\w) &=  -b(\u,\w,\v),\ \text{ for all }\ \u,\v,\w\in \V.
		\end{aligned}
		\right.\end{equation}
	\begin{remark}
		The following estimate on $b(\cdot,\cdot,\cdot)$ is helpful in the sequel (see Chapter 2, section 2.3, \cite{Temam1}). For all   $\u, \v, \w\in \V$,
		\begin{align}\label{b1}
			|b(\u,\v,\w)|&\leq C
			\begin{cases}
				\|\u\|^{\frac{1}{2}}_{\H}\|\nabla\u\|^{\frac{1}{2}}_{\H}\|\nabla\v\|_{\H}\|\w\|^{\frac{1}{2}}_{\H}\|\nabla\w\|^{\frac{1}{2}}_{\H},\ \ \ 	\text{ for } d=2,\\
				\|\u\|^{\frac{1}{4}}_{\H}\|\nabla\u\|^{\frac{3}{4}}_{\H}\|\nabla\v\|_{\H}\|\w\|^{\frac{1}{4}}_{\H}\|\nabla\w\|^{\frac{3}{4}}_{\H}, \ \ \ \text{ for } d=3.
			\end{cases}
		\end{align}
	\end{remark}
	\begin{remark}
		Note that $\langle\B(\v,\u-\v),\u-\v\rangle=0$ and it implies that
		\begin{align}\label{441}
			\langle \B(\u)-\B(\v),\u-\v\rangle &=\langle\B(\u-\v,\u),\u-\v\rangle\nonumber\\&=-\langle\B(\u-\v,\u-\v),\u\rangle=-\langle\B(\u-\v,\u-\v),\v\rangle.
		\end{align}
	\end{remark}
	\subsection{Nonlinear operator}
	Consider the nonlinear operator $\mathcal{C}(\u):=\mathcal{P}(|\u|^{r-1}\u)$. It is immediate that $\langle\mathcal{C}(\u),\u\rangle =\|\u\|_{\widetilde{\L}^{r+1}}^{r+1}$ and the map $\mathcal{C}(\cdot):\V\cap\widetilde{\L}^{r+1}\to\V'+\widetilde{\L}^{\frac{r+1}{r}}$. For all $\u\in\wi\L^{r+1}$, the map is Gateaux differentiable with Gateaux derivative 
	\begin{align}\label{29}
		\mathcal{C}'(\u)\v&=\left\{\begin{array}{cl}\mathcal{P}(\v),&\text{ for }r=1,\\ \left\{\begin{array}{cc}\mathcal{P}(|\u|^{r-1}\v)+(r-1)\mathcal{P}\left(\frac{\u}{|\u|^{3-r}}(\u\cdot\v)\right),&\text{ if }\u\neq \mathbf{0},\\\mathbf{0},&\text{ if }\u=\mathbf{0},\end{array}\right.&\text{ for } 1<r<3,\\ \mathcal{P}(|\u|^{r-1}\v)+(r-1)\mathcal{P}(\u|\u|^{r-3}(\u\cdot\v)), &\text{ for }r\geq 3,\end{array}\right.
	\end{align}
	for all $\v\in\V\cap\widetilde{\L}^{r+1}$. Moreover, for any $r\in [1, \infty)$ and $\u_1, \u_2 \in \V\cap\widetilde{\L}^{r+1}$, we have (see subsection 2.4, \cite{MTM1})
	\begin{align}\label{MO_c}
		\langle\mathcal{C}(\u_1)-\mathcal{C}(\u_2),\u_1-\u_2\rangle\geq\frac{1}{2}\||\u_1|^{\frac{r-1}{2}}(\u_1-\u_2)\|_{\H}^2+\frac{1}{2}\||\u_2|^{\frac{r-1}{2}}(\u_1-\u_2)\|_{\H}^2 \geq 0
	\end{align}
	and 
	\begin{align}\label{a215}
		\|\u-\v\|_{\wi\L^{r+1}}^{r+1}\leq 2^{r-2}\||\u|^{\frac{r-1}{2}}(\u-\v)\|_{\H}^2+2^{r-2}\||\v|^{\frac{r-1}{2}}(\u-\v)\|_{\H}^2,
	\end{align}
	for $r\geq 1$ (replace $2^{r-2}$ with $1,$ for $1\leq r\leq 2$).
	\subsection{Inequalities}
	The following inequalities are frequently used  in the paper. 
	\begin{lemma}\label{Holder}
		Assume that $\frac{1}{p}+\frac{1}{p'}=1$ with $1\leq p,p'\leq\infty$, $\u_1\in\L^{p}(\mathcal{O})$ and $\u_2\in\L^{p'}(\mathcal{O})$. Then we get 
		\begin{align*}
			\|\u_1\u_2\|_{\L^1(\mathcal{O})}\leq\|\u_1\|_{\L^{p}(\mathcal{O})}\|\u_2\|_{\L^{p'}(\mathcal{O})}. 
		\end{align*}
	\end{lemma}
	\begin{lemma}\label{Interpolation}
		Assume $1\leq s_1\leq s\leq s_2\leq \infty$, $\theta\in(0,1)$ such that $\frac{1}{s}=\frac{a}{s_1}+\frac{1-a}{s_2}$ and $\u\in\L^{s_1}(\mathcal{O})\cap\L^{s_2}(\mathcal{O})$, then we have 
		\begin{align*}
			\|\u\|_{\L^s(\mathcal{O})}\leq\|\u\|_{\L^{s_1}(\mathcal{O})}^{a}\|\u\|_{\L^{s_2}(\mathcal{O})}^{1-a}. 
		\end{align*}
	\end{lemma}
	\begin{lemma}\label{Young}
		For all $a,b,\varepsilon>0$ and for all $1<p,p'<\infty$ with $\frac{1}{p}+\frac{1}{p'}=1$, we obtain 
		\begin{align*}
			ab\leq\frac{\varepsilon}{p}a^p+\frac{1}{p'\varepsilon^{p'/p}}b^{p'}.
		\end{align*}
	\end{lemma}
	\subsection{White noise and colored noise}
	In this subsection, we recall some properties of white noise and colored noise.
	\begin{lemma}[Lemma 2.1, \cite{GGW}]\label{WnCn}
		Assume that the correlation time $\delta\in(0,1]$. There exists a $\{\vartheta_t\}_{t\in\R}$-invariant subset $\widetilde{\Omega}\subseteq\Omega$ of full measure, such that for $\omega\in\widetilde{\Omega}$,
		\begin{itemize}
			\item [(i)] \begin{align}\label{N1}
				\lim\limits_{t\to\pm\infty}\frac{\omega(t)}{t}=0;
			\end{align}
			\item[(ii)] the mapping
			\begin{align}\label{N2}
				(t,\omega)\mapsto\mathcal{Z}_{\delta}(\vartheta_t\omega)=-\frac{1}{\delta^2}\int_{-\infty}^{0}e^{\frac{\xi}{\delta}}\vartheta_t\omega(\xi)\d\xi
			\end{align}
			is a stationary solution (also called a \emph{colored noise} or an \emph{Ornstein-Uhlenbeck process}) of one-dimensional stochastic differential equation $$\d\mathcal{Z}_{\delta}+\frac{1}{\delta}\mathcal{Z}_{\delta}\d t=\frac{1}{\delta}\d\W$$ with continuous trajectories satisfying
			\begin{align}
				\lim\limits_{t\to\pm\infty}\frac{\left|\mathcal{Z}_{\delta}(\vartheta_t\omega)\right|}{t}&=0, \ \ \ \text{ for every } 0<\delta\leq1,\label{N3}\\
				\lim\limits_{t\to\pm\infty}\frac{1}{t}\int_0^t\mathcal{Z}_{\delta}(\vartheta_{\xi}\omega)\d\xi=\mathbb{E}[\mathcal{Z}_{\delta}]&=0, \ \ \ \text{ uniformly for } 0<\delta\leq1;\label{N4}
			\end{align}
		\end{itemize}
		and
		\begin{itemize}
			\item[(iii)] for arbitrary $T>0,$
			\begin{align}\label{N5*}
				\lim_{\delta\to0}\int_{0}^{t}\mathcal{Z}_{\delta}(\vartheta_{\xi}\omega)\d\xi=\omega(t) \ \text{ uniformly for } \ t\in[\mathfrak{s},\mathfrak{s}+T].
			\end{align}
		\end{itemize}
	\end{lemma}
	\begin{remark}
		For convenience, we use $\Omega$ itself in place of $\widetilde{\Omega}$ throughout the work.
	\end{remark}
	\subsection{Abstract formulation}
	In this subsection, we describe an abstract formulation and solution of \eqref{2}. Taking orthogonal projection $\mathcal{P}$ to the first equation of \eqref{2}, we obtain
	\begin{equation}\label{WZ_SCBF}
		\left\{
		\begin{aligned}
			\frac{\d \u(t)}{\d t}+\mu \A\u(t)+\B(\u(t))+\alpha\u(t)+\beta\mathcal{C}(\u(t))&=\boldsymbol{f}(t) + S(t,x,\u(t))\mathcal{Z}_{\delta}(\vartheta_t\omega), \ \ t>\mathfrak{s}, \\
			\u(x,\mathfrak{s})&=\u_{\mathfrak{s}}(x), \ \ \ \ \ \ \ \ \ x\in \mathbb{R}^n \text{ and }\mathfrak{s}\in\R.
		\end{aligned}
		\right.
	\end{equation}
Strictly speaking, one has to write $\mathcal{P}\f$ for $\f$ and $\mathcal{P}S(\cdot,\cdot,\cdot)$  for $S(\cdot,\cdot,\cdot)$ in \eqref{WZ_SCBF}.
	\begin{definition}\label{def3.1}
		Let us assume that $\mathfrak{s}\in\R,$ $ \omega\in\Omega,$ $ \u_{\mathfrak{s}}\in \H$, $\f\in \mathrm{L}^2_{\emph{loc}}(\R,\V')$, $T>0$ be any fixed time. Then, the function $\u(\cdot)$ is called a solution (in the weak sense) of the system \eqref{WZ_SCBF} on time interval $[\mathfrak{s},\mathfrak{s}+T]$, if $$\u\in \mathrm{L}^{\infty}(\mathfrak{s},\mathfrak{s}+T;\H)\cap\mathrm{L}^2(\mathfrak{s}, \mathfrak{s}+T;\V)\cap\mathrm{L}^{r+1}(\mathfrak{s},\mathfrak{s}+T;\widetilde{\L}^{r+1}),$$ with $ \partial_t\u\in \mathrm{L}^{2}(\mathfrak{s},\mathfrak{s}+T;\V')+\mathrm{L}^{\frac{r+1}{r}}(\mathfrak{s},\mathfrak{s}+T;\widetilde{\L}^{\frac{r+1}{r}})$ satisfying: 
		\begin{enumerate}
			\item [(i)]  	for any $\psi\in \V\cap\widetilde{\L}^{r+1},$ 
			\begin{align*}
				(\u(t), \psi) &= (\u_{\mathfrak{s}}, \psi)  -  \mu \int\limits_{\mathfrak{s}}^{\mathfrak{s}+T}(\nabla\u(\xi),\nabla\psi)\d\xi-\int\limits_{\mathfrak{s}}^{\mathfrak{s}+T}b(\u(\xi),\u(\xi),\psi)\d\xi-\alpha\int\limits_{\mathfrak{s}}^{\mathfrak{s}+T}(\u(\xi),\psi)\d\xi\nonumber\\&\quad-\beta \int\limits_{\mathfrak{s}}^{\mathfrak{s}+T}\left\langle|\u(\xi)|^{r-1}\u(\xi) +\f(\xi) , \psi\right\rangle\d\xi + \int\limits_{\mathfrak{s}}^{\mathfrak{s}+T}(S(\xi,x,\u(s))\mathcal{Z}_{\delta}(\vartheta_{\xi}\omega), \psi)\d\xi,
			\end{align*}
			for a.e. $t\in[\mathfrak{s},\mathfrak{s}+T]$,
			\item [(ii)] the initial data is satisfied in the following sense: 
			$$\lim\limits_{t\downarrow \mathfrak{s}}\int_{\mathcal{O}}\u(t,x)\psi(x)\d x=\int_{\mathcal{O}}\u_{\mathfrak{s}}(x)\psi(x)\d x, $$ for all $\psi\in\H$,
			\item [(iii)] the energy equality:
			\begin{align*}
				&\|\u(t)\|_{\H}^2+2\mu\int_{\mathfrak{s}}^t\|\u(s)\|_{\V}^2\d s+2\alpha\int_{\mathfrak{s}}^t\|\u(s)\|_{\H}^2\d s+2\beta\int_{\mathfrak{s}}^t\|\u(s)\|_{\wi\L^{r+1}}^{r+1}\d s\nonumber\\&=\|\u_{\mathfrak{s}}\|_{\H}^2+2\int_0^t(\f(s),\u(s))\d s+2\int_0^t(S(s,x,\u(s))\mathcal{Z}_{\delta}(\vartheta_s\omega),\u(s))\d s,
			\end{align*}
		for all $t\in[\mathfrak{s},T]$. 
		\end{enumerate}
	\end{definition}
	
	By a standard Galerkin method (see \cite{MTM,MTM1}), one can obtain that if any one of the assumptions  \ref{NDT1}, \ref{NDT2} and \ref{NDT3} are fulfilled, then for all $t>\mathfrak{s}, \ \mathfrak{s}\in\R,$ and for every $\u_{\mathfrak{s}}\in\H$ and $\omega\in\Omega$, \eqref{WZ_SCBF} has a unique solution in the sense of Definition \ref{def3.1}. Moreover, $\u(t,\mathfrak{s},\omega,\u_{\mathfrak{s}})$ is continuous with respect to initial data $\u_{\mathfrak{s}}(x)$ (see Lemmas \ref{Continuity}, \ref{ContinuityUB1} and \ref{ContinuityUB2}) and $(\mathscr{F},\mathscr{B}(\H))$-measurable in $\omega\in\Omega.$ 
	
	Now, we define a cocycle $\Phi:\R^+\times\R\times\Omega\times\H\to\H$ for the system \eqref{WZ_SCBF} such that for given $t\in\R^+, \mathfrak{s}\in\R, \omega\in\Omega$ and $\u_{\mathfrak{s}}\in\H$, let
	\begin{align}\label{Phi1}
		\Phi(t,\mathfrak{s},\omega,\u_{\mathfrak{s}}) =\u(t+\mathfrak{s},\mathfrak{s},\vartheta_{-\mathfrak{s}}\omega,\u_{\mathfrak{s}}).
	\end{align}
	Then $\Phi$ is a continuous cocycle on $\H$ over $(\Omega,\mathscr{F},\mathbb{P},\{\vartheta_{t}\}_{t\in\R})$, where $\{\vartheta_t\}_{t\in\R}$ is given by \eqref{vartheta}, that is,
	\begin{align}\label{Phi2}
		\Phi(t+\mathfrak{s},s,\omega,\u_{\mathfrak{s}})=\Phi(t,\mathfrak{s}+s,\vartheta_{\mathfrak{s}}\omega,\Phi(\mathfrak{s},s,\omega,\u_{\mathfrak{s}})).
	\end{align}
	
	Assume that $D=\{D(\mathfrak{s},\omega):\mathfrak{s}\in\R,\omega\in\Omega\}$ is a family of non-empty subsets of $\H$ satisfying, for every $c>0, \mathfrak{s}\in\R$ and $\omega\in\Omega$, 
	\begin{align}\label{D_1}
		\lim_{t\to\infty}e^{-ct}\|D(\mathfrak{s}-t,\vartheta_{-t}\omega)\|_{\H}=0,
	\end{align}
	where $\|D\|_{\H}=\sup\limits_{\x\in D}\|\x\|_{\H}.$ Let $\mathfrak{D}$ be the collection of all tempered families of bounded non-empty subsets of $\H$, that is,
	\begin{align}\label{D_11}
		\mathfrak{D}=\big\{D=\{D(\mathfrak{s},\omega):\mathfrak{s}\in\R\text{ and }\omega\in\Omega\}:D \text{ satisfying } \eqref{D_1}\big\}.
	\end{align}

	\section{Random pullback attractors for Wong-Zakai approximations: bounded domains} \label{sec3}\setcounter{equation}{0}
	In this section, we prove the existence of unique random $\mathfrak{D}$-pullback attractor for the system \eqref{WZ_SCBF} on bounded domains with nonlinear diffusion term. Throughout the section, we assume that $\mathcal{O}$ is bounded subset of $\R^d$.
	\subsection{Nonlinear diffusion term}
	The nonlinear diffusion term $S(t,x,\u)$ appearing in \eqref{WZ_SCBF} satisfies the following assumption:
	\begin{assumption}\label{NDT1}
		We assume that the nonlinear diffusion term $$S(t,x,\u)=e^{\sigma t}\left[\kappa\u+\mathcal{S}(\u)+\h(x)\right],$$ where $\sigma\geq0, \kappa\geq0$ and $\h\in\H$. Also, $\mathcal{S}:\V\to\H$ is a continuous function satisfying 
		\begin{align}
			\|\mathcal{S}(\u)-\mathcal{S}(\v)\|_{\H}&\leq s_1\|\u-\v\|_{\V}, \ \ \ \ \ \ \ \ \ \ \text{ for all } \u,\v\in\V,\label{S2}\\
			|\left(\mathcal{S}(\u)-\mathcal{S}(\v),\w\right)|&\leq s_2\|\u-\v\|_{\H}\|\w\|_{\V}, \ \ \ \text{ for all } \u,\v,\w \in\V,\label{S3}\\
			|\left(\mathcal{S}(\u),\u\right)|&\leq s_3 + s_4\|\u\|_{\V}^{1+s_5}, \ \ \ \ \ \ \text{ for all }\u\in\V,\label{S4}
		\end{align}
		where $s_1$ and $s_2$ both are non-negative constants, $s_3\geq0,  s_4\geq0$ and $s_5\in[0,1)$.
	\end{assumption}
	\begin{remark}
		One can take $\h\in\widetilde{\L}^{\frac{r+1}{r}}$ also.
	\end{remark}
	\begin{example}
		Let us discuss an example for such a nonlinear diffusion term which satisfies the above conditions. Let $\mathcal{S}:\V\to\H$ be a nonlinear operator defined by $\mathcal{S}(\u)=\sin\u+\B(\g_1,\u)$ for all $\u\in\V$, where $\g_1\in\D(\A)$ is a fixed element. It is easy to verify that $\mathcal{S}$ satisfies \eqref{S2}-\eqref{S4} (see \cite{GLW}).
	\end{example}
	\subsection{Deterministic nonautonomous forcing term}
	In the sequel, the following assumptions are needed on the external forcing term $\f$. 
	\begin{assumption}\label{DNFT1}
		There exists a number $\gamma\in[0,\alpha)$ such that:
		\begin{itemize}
			\item [(i)] 
			\begin{align}\label{forcing1}
				\int_{-\infty}^{\mathfrak{s}} e^{\gamma\xi}\|\f(\cdot,\xi)\|^2_{\V'}\d \xi<\infty, \ \ \text{ for all }\  \mathfrak{s}\in\R.
			\end{align}
			\item [(ii)] for every $c>0$
			\begin{align}\label{forcing2}
				\lim_{\tau\to-\infty}e^{c\tau}\int_{-\infty}^{0} e^{\gamma\xi}\|\f(\cdot,\xi+\tau)\|^2_{\V'}\d \xi=0.
			\end{align}
		\end{itemize}
		It should be noted that \eqref{forcing2} implies \eqref{forcing1} for $\f\in \mathrm{L}^2_{\emph{loc}}(\mathbb{R};\V')$.
	\end{assumption}
	
	\begin{lemma}\label{Continuity}
		For $d=2$ with $r\geq1$, $d=3$ with $r>3$ and $d=r=3$ with $2\beta\mu>1$, assume that $\f\in \mathrm{L}^2_{\emph{loc}}(\mathbb{R};\V')$ and Assumption \ref{NDT1} is fulfilled. Then, the solution of \eqref{WZ_SCBF} is continuous in initial data $\u_{\mathfrak{s}}(x).$
	\end{lemma}
	\begin{proof}
		Let $\u_{1}(\cdot)$ and $\u_{2}(\cdot)$ be two solutions of \eqref{WZ_SCBF}. Then $\mathfrak{X}(\cdot)=\u_{1}(\cdot)-\u_{2}(\cdot)$ with $\mathfrak{X}(\mathfrak{s})=\u_{1,\mathfrak{s}}(x)-\u_{2,\mathfrak{s}}(x)$ satisfies
		\begin{align}\label{Conti1}
			\frac{\d\mathfrak{X}(t)}{\d t}&=-\mu \A\mathfrak{X}(t)-\alpha\mathfrak{X}(t)-\left\{\B\big(\u_{1}(t)\big)-\B\big(\u_{2}(t)\big)\right\} -\beta\left\{\mathcal{C}\big(\u_1(t)\big)-\mathcal{C}\big(\u_2(t)\big)\right\}\nonumber\\&\quad+e^{\sigma t}\left[\kappa\mathfrak{X}(t)+\mathcal{S}(\u_1(t))-\mathcal{S}(\u_2(t))\right]\mathcal{Z}_{\delta}(\vartheta_t\omega),
		\end{align}
		in $\V'+\widetilde{\L}^{\frac{r+1}{r}}$ for a.e. $t\in[0,T]$. Taking the inner product with $\mathfrak{X}(\cdot)$ to the first equation in \eqref{Conti1}, we obtain
		\begin{align}\label{Conti2}
			\frac{1}{2}\frac{\d}{\d t} \|\mathfrak{X}(t)\|^2_{\H} &=-\mu \|\nabla\mathfrak{X}(t)\|^2_{\H} - \alpha\|\mathfrak{X}(t)\|^2_{\H} -\left\langle\B\big(\u_1(t)\big)-\B\big(\u_2(t)\big), \mathfrak{X}(t)\right\rangle \nonumber\\&\quad-\beta\left\langle\mathcal{C}\big(\u_1(t)\big)-\mathcal{C}\big(\u_2(t)\big),\mathfrak{X}(t)\right\rangle + \kappa e^{\sigma t}\mathcal{Z}_{\delta}(\vartheta_{t}\omega)\|\mathfrak{X}(t)\|^2_{\H}\nonumber\\&\quad+e^{\sigma t}\mathcal{Z}_{\delta}(\vartheta_{t}\omega)\left(\mathcal{S}(\u_1(t))-\mathcal{S}(\u_2(t)),\mathfrak{X}(t)\right) ,
		\end{align}
		for a.e. $t\in[\mathfrak{s},\mathfrak{s}+T] \text{ with } T>0$. By \eqref{S3}, we get
		\begin{align}\label{ContiS}
				\kappa & e^{\sigma t}\mathcal{Z}_{\delta}(\vartheta_{t}\omega)\|\mathfrak{X}\|^2_{\H}+e^{\sigma t}\mathcal{Z}_{\delta}(\vartheta_{t}\omega)\left(\mathcal{S}(\u_1)-\mathcal{S}(\u_2),\mathfrak{X}\right)\nonumber\\&\leq \kappa e^{\sigma t}\left|\mathcal{Z}_{\delta}(\vartheta_{t}\omega)\right|\|\mathfrak{X}\|^2_{\H}+s_2 e^{\sigma t}\left|\mathcal{Z}_{\delta}(\vartheta_{t}\omega)\right|\|\mathfrak{X}\|_{\H} \|\mathfrak{X}\|_{\V}\nonumber\\&\leq \frac{\alpha}{4}\|\mathfrak{X}\|^2_{\H}+ \frac{\min\{\mu,\alpha\}}{4}\|\mathfrak{X}\|^2_{\V} +Ce^{2\sigma t}\left|\mathcal{Z}_{\delta}(\vartheta_{t}\omega)\right|^2\|\mathfrak{X}\|^2_{\H}.
		\end{align}
		\vskip 2mm
		\noindent
		\textbf{Case I:} \textit{When $d=2$ and $r\geq1$.} Using \eqref{b1}, \eqref{441} and Lemma \ref{Young}, we obtain
		\begin{align}\label{Conti3}
			\left| \left\langle\B\big(\u_1\big)-\B\big(\u_2\big), \mathfrak{X}\right\rangle\right|&=\left|\left\langle\B\big(\mathfrak{X},\mathfrak{X} \big), \u_2\right\rangle\right|\leq\frac{\mu}{4}\|\nabla\mathfrak{X}\|^2_{\H}+C\|\u_2\|^4_{\widetilde{\L}^4}\|\mathfrak{X}\|^2_{\H}\nonumber\\&\leq\frac{\mu}{4}\|\nabla\mathfrak{X}\|^2_{\H}+C\|\u_2\|^2_{\H}\|\nabla\u_2\|^2_{\H}\|\mathfrak{X}\|^2_{\H},
		\end{align}
		and from \eqref{MO_c}, we have 
		\begin{align}\label{Conti4}
			-\beta\left\langle\mathcal{C}\big(\u_1\big)-\mathcal{C}\big(\u_2\big),\mathfrak{X}\right\rangle\leq 0.
		\end{align}
		Making use of \eqref{ContiS}-\eqref{Conti4} in \eqref{Conti2}, we get
		\begin{align*}
			&	\frac{\d}{\d t} \|\mathfrak{X}(t)\|^2_{\H}  \leq C\bigg\{e^{2\sigma t}\left|\mathcal{Z}_{\delta}(\vartheta_{t}\omega)\right|^2+ \|\u_2(t)\|^2_{\H}\|\nabla\u_2(t)\|^2_{\H}\bigg\}\|\mathfrak{X}(t)\|^2_{\H},
		\end{align*}
		for a.e. $t\in[\mathfrak{s},\mathfrak{s}+T]$ and an application of Gronwall's inequality completes the proof.
		\vskip 2mm
		\noindent
		\textbf{Case II:} \textit{When $d= 3$ and $r>3$.} The nonlinear term $|(\B(\mathfrak{X},\mathfrak{X}),\u_2)|$ is estimated using Lemmas \ref{Holder} and \ref{Young} as 
		\begin{align}\label{3d-ab12}
			\left|(\B(\mathfrak{X},\mathfrak{X}),\u_2)\right|&\leq \||\u_2||\mathfrak{X}|\|_{\H}\|\nabla\mathfrak{X}\|_{\H}\leq\frac{\mu}{4}\|\nabla\mathfrak{X}\|_{\H}^2+\frac{1}{\mu}\||\u_2||\mathfrak{X}|\|_{\H}^2\nonumber\\&\leq\frac{\mu}{4}\|\nabla\mathfrak{X}\|_{\H}^2+\frac{\beta}{2}\||\mathfrak{X}||\u_2|^{\frac{r-1}{2}}\|^2_{\H}+\eta_1\|\mathfrak{X}\|^2_{\H},
		\end{align}
		where $\eta_1= \frac{r-3}{2\mu(r-1)}\left[\frac{4}{\beta\mu (r-1)}\right]^{\frac{2}{r-3}}$  (see \cite{MTM1} for details). Using \eqref{441} and \eqref{MO_c}, we have
		\begin{align}\label{Conti7}
			\left|\left\langle\B\big(\u_1\big)-\B\big(\u_2\big), \mathfrak{X}\right\rangle\right|=\left|\left\langle\B\big(\mathfrak{X},\mathfrak{X} \big), \u_2\right\rangle\right|
		\end{align}
		and
		\begin{align}\label{Conti8}
			-\beta \left\langle\mathcal{C}\big(\u_1\big)-\mathcal{C}\big(\u_2\big),\mathfrak{X}\right\rangle\leq - \frac{\beta}{2}\||\mathfrak{X}||\u_2|^{\frac{r-1}{2}}\|^2_{\H},
		\end{align}
		respectively. Combining \eqref{Conti2}-\eqref{ContiS} and \eqref{3d-ab12}-\eqref{Conti8}, we get
		\begin{align*}
			&	\frac{\d}{\d t} \|\mathfrak{X}(t)\|^2_{\H}  \leq\biggl\{e^{2\sigma t}\left|\mathcal{Z}_{\delta}(\vartheta_{t}\omega)\right|^2+2\eta_1\biggr\}\|\mathfrak{X}(t)\|^2_{\H},
		\end{align*} for a.e. $t\in[\mathfrak{s},\mathfrak{s}+T]$, which completes the proof using Gronwall's inequality.
		\vskip 2mm
		\noindent
		\textbf{Case III:} \textit{When $d=3$ and $r=3$ with $2\beta\mu>1$.} Since $2\beta\mu>1$, there exists $0<\theta<1$ such that
		\begin{align}\label{Theta}
			2\beta\mu\geq\frac{1}{\theta}.
		\end{align} Using \eqref{b1} and \eqref{441}, we have
		\begin{align}\label{Conti5}
			\left|\left\langle\B\big(\u_1\big)-\B\big(\u_2\big), \mathfrak{X}\right\rangle\right|&=\left|\left\langle\B\big(\mathfrak{X},\mathfrak{X} \big), \u_2\right\rangle\right|\leq\theta\mu\|\nabla\mathfrak{X}\|^2_{\H}+\frac{1}{4\theta\mu}\||\u_2||\mathfrak{X}|\|^2_{\H},
		\end{align}
		where $\theta$ is the same as in \eqref{Theta}. Again by \eqref{S3}, we obtain
		\begin{align}\label{ContiS1}
				\kappa & e^{\sigma t}\mathcal{Z}_{\delta}(\vartheta_{t}\omega)\|\mathfrak{X}\|^2_{\H}+e^{\sigma t}\mathcal{Z}_{\delta}(\vartheta_{t}\omega)\left(\mathcal{S}(\u_1)-\mathcal{S}(\u_2),\mathfrak{X}\right)\nonumber\\&\leq \frac{\alpha}{2}\|\mathfrak{X}\|^2_{\H}+ (1-\theta)\mu\|\mathfrak{X}\|^2_{\V} +Ce^{2\sigma t}\left|\mathcal{Z}_{\delta}(\vartheta_{t}\omega)\right|^2\|\mathfrak{X}\|^2_{\H}.
		\end{align}
	Using \eqref{Conti8}-\eqref{Conti5} in \eqref{Conti2}, we get
		\begin{align*}
			&	\frac{\d}{\d t} \|\mathfrak{X}(t)\|^2_{\H}  \leq 2(1-\theta)\mu\|\mathfrak{X}(t)\|^2_{\H} +Ce^{2\sigma t}\left|\mathcal{Z}_{\delta}(\vartheta_{t}\omega)\right|^2\|\mathfrak{X}(t)\|^2_{\H} ,
		\end{align*} for a.e. $t\in[\mathfrak{s},\mathfrak{s}+T]$.
		After applying Gronwall's inequality, one can conclude the proof.
	\end{proof}
	Next result supports us to show the existence of random $\mathfrak{D}$-pullback absorbing set for continuous cocycle $\Phi$.
	\begin{lemma}\label{LemmaUe}
		For $d=2$ with $r\geq1$, $d=3$ with $r>3$ and $d=r=3$ with $2\beta\mu\geq1$, assume that $\f\in \mathrm{L}^2_{\emph{loc}}(\mathbb{R};\V')$ satisfies \eqref{forcing1} and Assumption \ref{NDT1} is fulfilled. Then for every $0<\delta\leq1$, $\mathfrak{s}\in\R,$ $ \omega\in \Omega$ and $D=\{D(\mathfrak{s},\omega):\mathfrak{s}\in\R, \omega\in\Omega\}\in\mathfrak{D},$ there exists $\mathcal{T}=\mathcal{T}(\delta, \mathfrak{s}, \omega, D)>0$ such that for all $t\geq \mathcal{T}$ and $\tau\geq \mathfrak{s}-t$, the solution $\u$ of the system \eqref{WZ_SCBF} with $\omega$ replaced by $\vartheta_{-\mathfrak{s}}\omega$ satisfies 
		\begin{align}\label{ue}
			&	\|\u(\tau,\mathfrak{s}-t,\vartheta_{-\mathfrak{s}}\omega,\u_{\mathfrak{s}-t})\|^2_{\H} \nonumber\\&\leq\frac{4}{\min\{\mu,\alpha\}} \int_{-\infty}^{\tau-\mathfrak{s}} e^{\int_{\tau-\mathfrak{s}}^{\xi}\left(\alpha-2\kappa e^{\sigma (\zeta+\mathfrak{s})}\mathcal{Z}_{\delta}(\vartheta_{\zeta}\omega)\right)\d\zeta} \|\f(\cdot,\xi+\mathfrak{s})\|^2_{\V'}\d \xi\nonumber\\&\quad+\int_{-\infty}^{\tau-\mathfrak{s}} e^{\int_{\tau-\mathfrak{s}}^{\xi}\left(\alpha-2\kappa e^{\sigma (\zeta+\mathfrak{s})}\mathcal{Z}_{\delta}(\vartheta_{\zeta}\omega)\right)\d\zeta}\bigg\{2s_3e^{\sigma (\xi+\mathfrak{s})}\left|\mathcal{Z}_{\delta}(\vartheta_{\xi}\omega)\right|\nonumber\\&\qquad+2s_6 e^{2\sigma (\xi+\mathfrak{s})}\left|\mathcal{Z}_{\delta}(\vartheta_{\xi}\omega)\right|^2+s_7\left[e^{\sigma(\xi+\mathfrak{s})}\left|\mathcal{Z}_{\delta}(\vartheta_{\xi}\omega)\right|\right]^{\frac{2}{1-s_5}}\bigg\}\d\xi,
		\end{align}
		where $\u_{\mathfrak{s}-t}\in D(\mathfrak{s}-t,\vartheta_{-t}\omega)$, $s_6=\frac{2}{\alpha}\|\h\|^2_{\H}$, $s_7=s_4(1-s_5)\left[\frac{4s_4(1+s_5)}{\min\{\mu,\alpha\}}\right]^{\frac{1+s_5}{1-s_5}}$ and, $s_3, s_4$ and $s_5$ are the constants appearing in \eqref{S4}.
	\end{lemma}
	\begin{proof}
		From the first equation of the system \eqref{WZ_SCBF}, we obtain
		\begin{align}\label{ue0}
			&	\frac{\d}{\d t} \|\u\|^2_{\H} +2\mu\|\nabla\u\|^2_{\H} + 2\alpha\|\u\|^2_{\H} + 2\beta\|\u\|^{r+1}_{\wi \L^{r+1}}\nonumber\\&= 2\left\langle\f,\u\right\rangle +2e^{\sigma t}\mathcal{Z}_{\delta}(\vartheta_{t}\omega)\left(\kappa\u+\mathcal{S}(\u)+\h,\u\right)\nonumber\\&\leq 2\|\f\|_{\V'}\|\u\|_{\V}+2\kappa e^{\sigma t}\mathcal{Z}_{\delta}(\vartheta_{t}\omega)\|\u\|^2_{\H}+2e^{\sigma t}\left|\mathcal{Z}_{\delta}(\vartheta_{t}\omega)\right|\left(s_3+s_4\|\u\|^{1+s_5}_{\V}+\|\h\|_{\H}\|\u\|_{\H}\right)\nonumber\\&\leq\frac{\alpha}{2}\|\u\|^2_{\H}+ \frac{\min\{\mu,\alpha\}}{2}\|\u\|^2_{\V}+\frac{4\|\f\|^2_{\V'}}{\min\{\mu,\alpha\}}+2\kappa e^{\sigma t}\mathcal{Z}_{\delta}(\vartheta_{t}\omega)\|\u\|^2_{\H}+2s_3e^{\sigma t}\left|\mathcal{Z}_{\delta}(\vartheta_{t}\omega)\right|\nonumber\\&\quad+s_6 e^{2\sigma t}\left|\mathcal{Z}_{\delta}(\vartheta_{t}\omega)\right|^2+s_7\left[e^{\sigma t}\left|\mathcal{Z}_{\delta}(\vartheta_{t}\omega)\right|\right]^{\frac{2}{1-s_5}},
		\end{align}
		where we have used \eqref{b0}, \eqref{S4}, Lemmas \ref{Holder} and \ref{Young}, and the constants $s_6=\frac{2}{\alpha}\|\h\|^2_{\H}$ and $s_7=s_4(1-s_5)\left[\frac{4s_4(1+s_5)}{\min\{\mu,\alpha\}}\right]^{\frac{1+s_5}{1-s_5}}$. We rewrite  \eqref{ue0} as
		\begin{align}\label{ue1}
			&\frac{\d}{\d t} \|\u\|^2_{\H}+ \left(\alpha-2\kappa e^{\sigma t}\mathcal{Z}_{\delta}(\vartheta_{t}\omega)\right)\|\u\|^2_{\H}\nonumber\\& \leq \frac{4\|\f\|^2_{\V'}}{\min\{\mu,\alpha\}}+2s_3e^{\sigma t}\left|\mathcal{Z}_{\delta}(\vartheta_{t}\omega)\right|+s_6 e^{2\sigma t}\left|\mathcal{Z}_{\delta}(\vartheta_{t}\omega)\right|^2+s_7\left[e^{\sigma t}\left|\mathcal{Z}_{\delta}(\vartheta_{t}\omega)\right|\right]^{\frac{2}{1-s_5}}.
		\end{align}
		Making use of variation of constant formula in \eqref{ue1} and replacing $\omega$ by $\vartheta_{-\mathfrak{s}}\omega$, we get
		\begin{align}\label{ue2}
			&	\|\u(\tau,\mathfrak{s}-t,\vartheta_{-\mathfrak{s}}\omega,\u_{\mathfrak{s}-t})\|^2_{\H} \nonumber\\&\leq e^{\int_{\tau}^{\mathfrak{s}-t}\left(\alpha-2\kappa e^{\sigma \zeta}\mathcal{Z}_{\delta}(\vartheta_{\zeta-\mathfrak{s}}\omega)\right)\d\zeta}\|\u_{\mathfrak{s}-t}\|^2_{\H}+ \frac{4}{\min\{\mu,\alpha\}}\int_{\mathfrak{s}-t}^{\tau} e^{\int_{\tau}^{\xi}\left(\alpha-2\kappa e^{\sigma \zeta}\mathcal{Z}_{\delta}(\vartheta_{\zeta-\mathfrak{s}}\omega)\right)\d\zeta} \|\f(\cdot,\xi)\|^2_{\V'}\d \xi\nonumber\\&\quad+\int_{\mathfrak{s}-t}^{\tau} e^{\int_{\tau}^{\xi}\left(\alpha-2\kappa e^{\sigma \zeta}\mathcal{Z}_{\delta}(\vartheta_{\zeta-\mathfrak{s}}\omega)\right)\d\zeta}\bigg\{2s_3e^{\sigma \xi}\left|\mathcal{Z}_{\delta}(\vartheta_{\xi-\mathfrak{s}}\omega)\right|+s_6 e^{2\sigma \xi}\left|\mathcal{Z}_{\delta}(\vartheta_{\xi-\mathfrak{s}}\omega)\right|^2\nonumber\\&\qquad+s_7\left[e^{\sigma\xi}\left|\mathcal{Z}_{\delta}(\vartheta_{\xi-\mathfrak{s}}\omega)\right|\right]^{\frac{2}{1-s_5}}\bigg\}\d\xi\nonumber\\&\leq e^{\int_{\tau-\mathfrak{s}}^{-t}\left(\alpha-2\kappa e^{\sigma (\zeta+\mathfrak{s})}\mathcal{Z}_{\delta}(\vartheta_{\zeta}\omega)\right)\d\zeta}\|\u_{\mathfrak{s}-t}\|^2_{\H}\nonumber\\&\quad+ \frac{4}{\min\{\mu,\alpha\}}\int_{-t}^{\tau-\mathfrak{s}} e^{\int_{\tau-\mathfrak{s}}^{\xi}\left(\alpha-2\kappa e^{\sigma (\zeta+\mathfrak{s})}\mathcal{Z}_{\delta}(\vartheta_{\zeta}\omega)\right)\d\zeta} \|\f(\cdot,\xi+\mathfrak{s})\|^2_{\V'}\d \xi\nonumber\\&\quad+\int_{-t}^{\tau-\mathfrak{s}} e^{\int_{\tau-\mathfrak{s}}^{\xi}\left(\alpha-2\kappa e^{\sigma (\zeta+\mathfrak{s})}\mathcal{Z}_{\delta}(\vartheta_{\zeta}\omega)\right)\d\zeta}\bigg\{2s_3e^{\sigma (\xi+\mathfrak{s})}\left|\mathcal{Z}_{\delta}(\vartheta_{\xi}\omega)\right|+s_6 e^{2\sigma (\xi+\mathfrak{s})}\left|\mathcal{Z}_{\delta}(\vartheta_{\xi}\omega)\right|^2\nonumber\\&\qquad+s_7\left[e^{\sigma(\xi+\mathfrak{s})}\left|\mathcal{Z}_{\delta}(\vartheta_{\xi}\omega)\right|\right]^{\frac{2}{1-s_5}}\bigg\}\d\xi.
		\end{align}
		Now, we need to estimate each term on the right hand side of \eqref{ue2}. Depending on $\sigma$, there are two possible cases. 
		\vskip1mm
		\noindent
		\textbf{Case I:} \textit{When $\sigma=0$.} By \eqref{N4}, we have
		\begin{align}\label{ue3}
			\lim_{\xi\to-\infty}\frac{1}{\xi}\int_{0}^{\xi}\left(\alpha-2\kappa \mathcal{Z}_{\delta}(\vartheta_{\zeta}\omega)\right)\d\zeta=\alpha-2\kappa\mathbb{E}[\mathcal{Z}_{\delta}]=\alpha.
		\end{align}
		Since $\gamma<\alpha$ (see Assumption \ref{DNFT1}), from \eqref{ue3}, we infer that there exists $\xi_0=\xi_0(\delta,\omega)<0$ such that for all $\xi\leq\xi_0$,
		\begin{align}\label{ue4}
			\int_{0}^{\xi}\left(\alpha-2\kappa \mathcal{Z}_{\delta}(\vartheta_{\zeta}\omega)\right)\d\zeta<\gamma\xi.
		\end{align}
		Taking \eqref{forcing1} and \eqref{ue4} into account, we obtain
		\begin{align}\label{ue5}
			\int_{-\infty}^{\xi_0}e^{\int_{0}^{\xi}\left(\alpha-2\kappa \mathcal{Z}_{\delta}(\vartheta_{\zeta}\omega)\right)\d\zeta}\|\f(\cdot,\xi+\mathfrak{s})\|^2_{\V'}\d\xi<\int_{-\infty}^{\xi_0}e^{\gamma\xi}\|\f(\cdot,\xi+\mathfrak{s})\|^2_{\V'}\d\xi<\infty,
		\end{align}
		and similarly from \eqref{N3} and \eqref{ue4}, we get
		\begin{align}\label{ue6}
			\int_{-\infty}^{\xi_0}e^{\int_{0}^{\xi}\left(\alpha-2\kappa \mathcal{Z}_{\delta}(\vartheta_{\zeta}\omega)\right)\d\zeta}\bigg\{2s_3\left|\mathcal{Z}_{\delta}(\vartheta_{\xi}\omega)\right|+s_6\left|\mathcal{Z}_{\delta}(\vartheta_{\xi}\omega)\right|^2+s_7\left|\mathcal{Z}_{\delta}(\vartheta_{\xi}\omega)\right|^{\frac{2}{1-s_5}}\bigg\}\d\xi<\infty.
		\end{align}
		\vskip1mm
		\noindent
		\textbf{Case II:} \textit{When $\sigma>0$.} From \eqref{N3}, we get
		\begin{align*}
			\lim_{\zeta\to-\infty}\left(\alpha-2\kappa e^{\sigma(\zeta+\mathfrak{s})}\mathcal{Z}_{\delta}(\vartheta_{\zeta}\omega)\right)=\alpha>\gamma,
		\end{align*}
		therefore, there exists $\zeta_0=\zeta_0(\mathfrak{s},\delta,\omega)<0$ such that for all $\zeta\leq\zeta_0$,
		\begin{align}\label{ue7}
			\alpha-2\kappa e^{\sigma(\zeta+\mathfrak{s})}\mathcal{Z}_{\delta}(\vartheta_{\zeta}\omega)>\gamma.
		\end{align}
		By \eqref{forcing1} and \eqref{ue7}, we have
		\begin{align}\label{ue8}
			\int_{-\infty}^{\zeta_0}e^{\int_{\zeta_0}^{\xi}\left(\alpha-2\kappa e^{\sigma(\zeta+\mathfrak{s})} \mathcal{Z}_{\delta}(\vartheta_{\zeta}\omega)\right)\d\zeta}\|\f(\cdot,\xi+\mathfrak{s})\|^2_{\V'}\d\xi< e^{-\gamma\zeta_0}\int_{-\infty}^{\zeta_0}e^{\gamma\xi}\|\f(\cdot,\xi+\mathfrak{s})\|^2_{\V'}\d\xi<\infty,
		\end{align}
		and similarly by \eqref{N3} and \eqref{ue7}, it is easy to obtain
		\begin{align}\label{ue9}
			&	\int_{-\infty}^{\xi_0}e^{\int_{0}^{\xi}\left(\alpha-2\kappa e^{\sigma(\zeta+\mathfrak{s})} \mathcal{Z}_{\delta}(\vartheta_{\zeta}\omega)\right)\d\zeta}\bigg\{2s_3e^{\sigma(\xi+\mathfrak{s})}\left|\mathcal{Z}_{\delta}(\vartheta_{\xi}\omega)\right|+s_6e^{2\sigma(\xi+\mathfrak{s})}\left|\mathcal{Z}_{\delta}(\vartheta_{\xi}\omega)\right|^2\nonumber\\&\qquad+s_7\left[e^{\sigma(\xi+\mathfrak{s})}\left|\mathcal{Z}_{\delta}(\vartheta_{\xi}\omega)\right|\right]^{\frac{2}{1-s_5}}\bigg\}\d\xi<\infty.
		\end{align}
		It follows from \eqref{ue5}-\eqref{ue6} and \eqref{ue8}-\eqref{ue9} that for any $\sigma\geq0$,
		\begin{align}\label{ue10}
			&\int_{-\infty}^{\tau-\mathfrak{s}} e^{\int_{\tau-\mathfrak{s}}^{\xi}\left(\alpha-2\kappa e^{\sigma (\zeta+\mathfrak{s})}\mathcal{Z}_{\delta}(\vartheta_{\zeta}\omega)\right)\d\zeta} \|\f(\cdot,\xi+\mathfrak{s})\|^2_{\V'}\d \xi+\int_{-\infty}^{\tau-\mathfrak{s}} e^{\int_{\tau-\mathfrak{s}}^{\xi}\left(\alpha-2\kappa e^{\sigma (\zeta+\mathfrak{s})}\mathcal{Z}_{\delta}(\vartheta_{\zeta}\omega)\right)\d\zeta}\nonumber\\&\quad\times\bigg\{2s_3e^{\sigma (\xi+\mathfrak{s})}\left|\mathcal{Z}_{\delta}(\vartheta_{\xi}\omega)\right|+s_6 e^{2\sigma (\xi+\mathfrak{s})}\left|\mathcal{Z}_{\delta}(\vartheta_{\xi}\omega)\right|^2+s_7\left[e^{\sigma(\xi+\mathfrak{s})}\left|\mathcal{Z}_{\delta}(\vartheta_{\xi}\omega)\right|\right]^{\frac{2}{1-s_5}}\bigg\}\d\xi<\infty.
		\end{align}
		Finally, since $\u_{\mathfrak{s}-t}\in D(\mathfrak{s}-t,\vartheta_{-t}\omega)$ and $D\in\mathfrak{D}$, making use of \eqref{ue4} and \eqref{ue7}, we have for $\sigma\geq0$,
		\begin{align*}
			e^{\int_{0}^{-t}\left(\alpha-2\kappa e^{\sigma (\zeta+\mathfrak{s})}\mathcal{Z}_{\delta}(\vartheta_{\zeta}\omega)\right)\d\zeta}\|\u_{\mathfrak{s}-t}\|^2_{\H}\leq e^{\int_{0}^{-t}\left(\alpha-2\kappa e^{\sigma (\zeta+\mathfrak{s})}\mathcal{Z}_{\delta}(\vartheta_{\zeta}\omega)\right)\d\zeta}\|D(\mathfrak{s}-t,\vartheta_{-t}\omega)\|^2_{\H}\to0,
		\end{align*}
		as $t\to\infty$, there exists $\mathcal{T}=\mathcal{T}(\delta,\mathfrak{s},\omega,D)>0$ such that for all $t\geq\mathcal{T}$,
		\begin{align}\label{ue11}
			&e^{\int_{\tau-\mathfrak{s}}^{-t}\left(\alpha-2\kappa e^{\sigma (\zeta+\mathfrak{s})}\mathcal{Z}_{\delta}(\vartheta_{\zeta}\omega)\right)\d\zeta}\|\u_{\mathfrak{s}-t}\|^2_{\H}\nonumber\\&\leq s_6\int_{-\infty}^{\tau-\mathfrak{s}} e^{\int_{\tau-\mathfrak{s}}^{\xi}\left(\alpha-2\kappa e^{\sigma (\zeta+\mathfrak{s})}\mathcal{Z}_{\delta}(\vartheta_{\zeta}\omega)\right)\d\zeta} e^{2\sigma (\xi+\mathfrak{s})}\left|\mathcal{Z}_{\delta}(\vartheta_{\xi}\omega)\right|^2\d\xi.
		\end{align}
		From \eqref{ue2} along with \eqref{ue10}-\eqref{ue11}, one can conclude the proof.
	\end{proof}
	
	Further, we prove $\mathfrak{D}$-pullback asymptotic compactness using compactness of Sobolev embeddings on bounded domains.
	
	\begin{lemma}\label{PCB}
		For $d=2$ with $r\geq1$, $d=3$ with $r>3$ and $d=r=3$ with $2\beta\mu\geq1$, assume that $\f\in \mathrm{L}^2_{\emph{loc}}(\mathbb{R};\V')$ satisfies \eqref{forcing1} and Assumption \ref{NDT1} is fulfilled. Then for every $0<\delta\leq1, \omega\in\Omega, \mathfrak{s}\in\R$ and $t>\mathfrak{s}$, the solution $\u(t,\mathfrak{s},\omega,\cdot):\H\to\H$ is compact, that is, for every bounded set $B$ in $\H$, the image $\u(t,\mathfrak{s},\omega,B)$ is precompact in $\H$.
	\end{lemma}
	\begin{proof}
		Consider the solution $\u(\tau,\mathfrak{s},\omega,\cdot)$ of \eqref{WZ_SCBF} for $\tau\in[\mathfrak{s},\mathfrak{s}+T]$, where $T>0$. Assume that the sequence $\{\u_{0,n}\}_{n\in\N}\subset B$. We know that (see the proof of Lemma \ref{LemmaUe}) 
		\begin{equation}\label{PCB1}
			\left\{
			\begin{aligned}
				\{\u(\cdot,\mathfrak{s},\omega,\u_{0,n})\}_{n\in\N} & \text{ is bounded in }\\ \mathrm{L}^{\infty}(\mathfrak{s},\mathfrak{s}+T;\H)\cap\mathrm{L}^2(\mathfrak{s}, \mathfrak{s}+T;&\V)\cap\mathrm{L}^{r+1}(\mathfrak{s},\mathfrak{s}+T;\widetilde{\L}^{r+1}).
			\end{aligned}
			\right.
		\end{equation}
		From \eqref{S2} along with \eqref{PCB1}, we obtain
		\begin{align}\label{PCB2}
			\{S(\cdot,\cdot,\u(\tau,\mathfrak{s},\omega,\u_{0,n}))\mathcal{Z}_{\delta}(\vartheta_{\tau}\omega)\}_{n\in\N}  \text{ is bounded in } \mathrm{L}^{2}(\mathfrak{s},\mathfrak{s}+T;\H).
		\end{align}
		We also have
		\begin{align}\label{PCB3}
			\{\A(\u(\tau,\mathfrak{s},\omega,\u_{0,n}))\}_{n\in\N} \text{ and }   \{\B(\u(\tau,\mathfrak{s},\omega,\u_{0,n}))\}_{n\in\N} \text{ are bounded in } \mathrm{L}^2(\mathfrak{s},\mathfrak{s}+T;\V'),
		\end{align}
		and
		\begin{align}\label{PCB4}
			\{\mathcal{C}(\u(\tau,\mathfrak{s},\omega,\u_{0,n}))\}_{n\in\N} \text{ is bounded in } \mathrm{L}^{\frac{r+1}{r}}(\mathfrak{s},\mathfrak{s}+T;\widetilde{\L}^{\frac{r+1}{r}}).
		\end{align}
		It follows from \eqref{PCB1}-\eqref{PCB4} and \eqref{WZ_SCBF} that
		\begin{align*}
			\left\{\frac{\d}{\d s}(\u(\tau,\mathfrak{s},\omega,\u_{0,n}))\right\}_{n\in\N} \text{ is bounded in } \mathrm{L}^2(\mathfrak{s},\mathfrak{s}+T;\V')+\mathrm{L}^{\frac{r+1}{r}}(\mathfrak{s},\mathfrak{s}+T;\widetilde{\L}^{\frac{r+1}{r}}).
		\end{align*}
		Since $\mathrm{L}^2(\mathfrak{s},\mathfrak{s}+T;\V')+\mathrm{L}^{\frac{r+1}{r}}(\mathfrak{s},\mathfrak{s}+T;\widetilde{\L}^{\frac{r+1}{r}})\subset\mathrm{L}^{\frac{r+1}{r}}(\mathfrak{s}+T;\V'+\wi\L^{\frac{r+1}{r}})$, the above sequence is bounded in $\mathrm{L}^{\frac{r+1}{r}}(\mathfrak{s}+T;\V'+\wi\L^{\frac{r+1}{r}})$. Note also that $\V\cap\wi\L^{\frac{r+1}{r}}\subset\V\subset\H\subset\V'\subset \V'+\wi\L^{\frac{r+1}{r}}$ and the embedding of $\V\subset\H$ is compact. 	By the \emph{Aubin-Lions compactness lemma}, there exists a subsequence (keeping as it is) and $\v\in\mathrm{L}^2(\mathfrak{s},\mathfrak{s}+T;\H)$ such that 
		\begin{align}\label{PCB6}
			\u(\cdot,\mathfrak{s},\omega,\u_{0,n})\to\v(\cdot) \ \text{ strongly in }\  \mathrm{L}^{2}(\mathfrak{s},\mathfrak{s}+T;\H).
		\end{align}
		Along a further subsequence (again not relabeling), we infer from \eqref{PCB6} that
		\begin{align}\label{PCB7}
			\u(\tau,\mathfrak{s},\omega,\u_{0,n})\to\v(\tau) \text{ in }  \H \ \text{ for almost all }\  \tau\in(\mathfrak{s},\mathfrak{s}+T).
		\end{align}
		Since $\mathfrak{s}<t<T$, we obtain from \eqref{PCB7} that there exists $\tau\in(\mathfrak{s},t)$ such that \eqref{PCB7} holds true for this particular $\tau$. Then  by Lemma \ref{Continuity}, we obtain
		\begin{align*}
			\u(t,\mathfrak{s},\omega,\u_{0,n})=\u(t,\tau,\omega,\u(\tau,\mathfrak{s},\omega,\u_{0,n}))\to \u(t,\tau,\omega,\v(\tau)),
		\end{align*}
		which completes the proof.
	\end{proof}

	In fact, Lemma \ref{PCB} helps us to prove the $\mathfrak{D}$-pullback asymptotic compactness of $\Phi$ in $\H$ on bounded domains.
	\begin{corollary}\label{Asymptotic_B}
		For $d=2$ with $r\geq1$, $d=3$ with $r>3$ and $d=r=3$ with $2\beta\mu\geq1$, assume that $\f\in\mathrm{L}^2_{\emph{loc}}(\R;\V')$ satisfies \eqref{forcing1} and Assumption \ref{NDT1} is fulfilled. Then for every $0<\delta\leq1$, $\mathfrak{s}\in \R,$ $\omega\in \Omega,$ $D=\{D(\mathfrak{s},\omega):\mathfrak{s}\in \R,\omega\in \Omega\}\in \mathfrak{D}$ and $t_n\to \infty,$ $\u_{0,n}\in D(\mathfrak{s}-t_n, \vartheta_{-t_{n}}\omega)$, the sequence $\Phi(t_n,\mathfrak{s}-t_n,\vartheta_{-t_n}\omega,\u_{0,n})$ or $\u(\mathfrak{s},\mathfrak{s}-t_n,\vartheta_{-\mathfrak{s}}\omega,\u_{0,n})$ of solutions of the system \eqref{WZ_SCBF} has a convergent subsequence in $\H$.
	\end{corollary}
	\begin{proof}
		From Lemma \ref{LemmaUe} with $\tau=\mathfrak{s}-1$, we have that there exists $\mathcal{T}=\mathcal{T}(\delta,\mathfrak{s},\omega,D)>0$ such that for all $t\geq\mathcal{T}$ and $\u_{\mathfrak{s}-t}\in D(\mathfrak{s}-t, \vartheta_{-t}\omega)$,
		\begin{align}\label{AB1}
			\u(\mathfrak{s}-1,\mathfrak{s}-t,\vartheta_{-\mathfrak{s}}\omega,\u_{\mathfrak{s}-t})\in\H.
		\end{align}
		Since $t_n\to\infty$ and $\u_{0,n}\in D(\mathfrak{s}-t_n, \vartheta_{-t_{n}}\omega)$, from \eqref{AB1}, we infer that there exists $N_1=N(\delta,\mathfrak{s},\omega,D)>0$ such that 
		\begin{align}\label{AB2}
			\{\u(\mathfrak{s}-1,\mathfrak{s}-t_n,\vartheta_{-\mathfrak{s}}\omega,\u_{0,n})\}_{n\geq N_1}\subset\H.
		\end{align}
		Hence, by \eqref{AB2} and Lemma \ref{PCB}, we conclude that the sequence $$\u(\mathfrak{s},\mathfrak{s}-t_n,\vartheta_{-\mathfrak{s}}\omega,\u_{0,n})=\u(\mathfrak{s},\mathfrak{s}-1,\vartheta_{-\mathfrak{s}}\omega,\u(\mathfrak{s}-1,\mathfrak{s}-t_n,\vartheta_{-\mathfrak{s}}\omega,\u_{0,n}))$$ has a convergent subsequence in $\H$, which completes the proof.
	\end{proof}

	\subsection{Existence of a unique random $\mathfrak{D}$-pullback attractor}
	In this subsection, we start with the result on the existence of a $\mathfrak{D}$-pullback absorbing set in $\H$ for the system \eqref{WZ_SCBF}. Then, we prove the main result of this section, that is, the existence of a unique random $\mathfrak{D}$-pullback attractor for the system \eqref{WZ_SCBF}.
	\begin{lemma}\label{PAS}
		For $d=2$ with $r\geq1$, $d=3$ with $r>3$ and $d=r=3$ with $2\beta\mu\geq1$, assume that $\f\in\mathrm{L}^2_{\emph{loc}}(\R;\V')$ satisfies \eqref{forcing2} and Assumption \ref{NDT1} is fulfilled. Then there exists a closed measurable $\mathfrak{D}$-pullback absorbing set $\mathcal{K}=\{\mathcal{K}(\mathfrak{s},\omega):\mathfrak{s}\in\R, \omega\in\Omega\}\in\mathfrak{D}$ for the  continuous cocycle $\Phi$ associated with the system \eqref{WZ_SCBF}.
	\end{lemma}
	\begin{proof}
		Let us denote, for given $\mathfrak{s}\in\R$ and $\omega\in\Omega$
		\begin{align}
			\mathcal{K}(\mathfrak{s},\omega)=\{\u\in\H:\|\u\|^2_{\H}\leq\mathcal{L}(\mathfrak{s},\omega)\},
		\end{align}
		where
		\begin{align}
			\mathcal{L}(\mathfrak{s},\omega)&=\frac{4}{\min\{\mu,\alpha\}} \int_{-\infty}^{0} e^{\int_{0}^{\xi}\left(\alpha-2\kappa e^{\sigma (\zeta+\mathfrak{s})}\mathcal{Z}_{\delta}(\vartheta_{\zeta}\omega)\right)\d\zeta} \|\f(\cdot,\xi+\mathfrak{s})\|^2_{\V'}\d \xi\nonumber\\&\quad+\int_{-\infty}^{0} e^{\int_{0}^{\xi}\left(\alpha-2\kappa e^{\sigma (\zeta+\mathfrak{s})}\mathcal{Z}_{\delta}(\vartheta_{\zeta}\omega)\right)\d\zeta}\bigg\{2s_3e^{\sigma (\xi+\mathfrak{s})}\left|\mathcal{Z}_{\delta}(\vartheta_{\xi}\omega)\right|\nonumber\\&\qquad+2s_6 e^{2\sigma (\xi+\mathfrak{s})}\left|\mathcal{Z}_{\delta}(\vartheta_{\xi}\omega)\right|^2+s_7\left[e^{\sigma(\xi+\mathfrak{s})}\left|\mathcal{Z}_{\delta}(\vartheta_{\xi}\omega)\right|\right]^{\frac{2}{1-s_5}}\bigg\}\d\xi.
		\end{align}
		Since $\mathcal{L}(\mathfrak{s},\cdot):\Omega\to\R$ is $(\mathscr{F},\mathscr{B}(\R))$-measurable for every $\mathfrak{s}\in\R$, $\mathcal{K}(\mathfrak{s},\cdot):\Omega\to2^{\H}$ is a measurable set-valued mapping. Furthermore, we have from Lemma \ref{LemmaUe} that for each $\mathfrak{s}\in\R$, $\omega\in\Omega$ and $D\in\mathfrak{D}$, there exists $\mathcal{T}=\mathcal{T}(\delta,\mathfrak{s},\omega,D)>0$ such that for all $t\geq\mathcal{T}$,
		\begin{align}\label{PAS1}
			\Phi(t,\mathfrak{s}-t,\vartheta_{-t}\omega,D(\mathfrak{s}-t,\vartheta_{-t}\omega))=\u(\mathfrak{s},\mathfrak{s}-t,\vartheta_{-\mathfrak{s}}\omega,D(\mathfrak{s}-t,\vartheta_{-t}\omega))\subseteq\mathcal{K}(\mathfrak{s},\omega).
		\end{align}
		Now, in order to complete the proof, we only need to prove that $\mathcal{K}\in\mathfrak{D}$, that is, for every $c>0$, $\mathfrak{s}\in\R$ and $\omega\in\Omega$ $$\lim_{t\to-\infty}e^{ct}\|\mathcal{K}(\mathfrak{s} +t,\vartheta_{t}\omega)\|^2_{\H}=0.$$ 
		For every $c>0$, $\mathfrak{s}\in\R$ and $\omega\in\Omega$,
		\begin{align}\label{PAS2}
			&	\lim_{t\to-\infty}e^{ct}\|\mathcal{K}(\mathfrak{s} +t,\vartheta_{t}\omega)\|^2_{\H}\nonumber\\&=\lim_{t\to-\infty}e^{ct}\mathcal{L}(\mathfrak{s} +t,\vartheta_{t}\omega)\nonumber\\&=\lim_{t\to-\infty}\frac{4e^{ct}}{\min\{\mu,\alpha\}} \int_{-\infty}^{0} e^{\int_{0}^{\xi}\left(\alpha-2\kappa e^{\sigma (\zeta+\mathfrak{s}+t)}\mathcal{Z}_{\delta}(\vartheta_{\zeta+t}\omega)\right)\d\zeta} \|\f(\cdot,\xi+\mathfrak{s}+t)\|^2_{\V'}\d \xi\nonumber\\&\quad+\lim_{t\to-\infty}e^{ct}\int_{-\infty}^{0} e^{\int_{0}^{\xi}\left(\alpha-2\kappa e^{\sigma (\zeta+\mathfrak{s}+t)}\mathcal{Z}_{\delta}(\vartheta_{\zeta+t}\omega)\right)\d\zeta}\bigg\{2s_3e^{\sigma (\xi+\mathfrak{s}+t)}\left|\mathcal{Z}_{\delta}(\vartheta_{\xi+t}\omega)\right|\nonumber\\&\qquad\qquad\qquad\qquad+2s_6 e^{2\sigma (\xi+\mathfrak{s}+t)}\left|\mathcal{Z}_{\delta}(\vartheta_{\xi+t}\omega)\right|^2+s_7\left[e^{\sigma(\xi+\mathfrak{s}+t)}\left|\mathcal{Z}_{\delta}(\vartheta_{\xi+t}\omega)\right|\right]^{\frac{2}{1-s_5}}\bigg\}\d\xi.
		\end{align}
		Taking into account of \eqref{N1}-\eqref{N2} and following the same steps applied in the proof of Lemma 2.6,  \cite{GGW}, we obtain that there exists a $\mathcal{T}_1=\mathcal{T}_1(\omega)>0$ such that for $\xi\leq0$, $\sigma>0$, $\kappa>0$, $0<\delta\leq1$ with $\sigma\delta\neq1$ and $t\leq-\mathcal{T}_1$,
		\begin{align}\label{PAS3}
			-2\kappa \int_{0}^{\xi}e^{\sigma (\zeta+\mathfrak{s}+t)}\mathcal{Z}_{\delta}(\vartheta_{\zeta+t}\omega)\d\zeta\leq c_0 \left(\frac{1}{\sigma}+2-\xi-2t\right),
		\end{align}
		where $c_0=\min\left\{\alpha-\gamma,\frac{c}{4\kappa}\right\}$. 
		
		Let $c_1=\min\left\{\frac{c}{2},\gamma+\sigma\right\}.$ It implies from \eqref{PAS2}-\eqref{PAS3} that, $\sigma>0$, $\kappa>0$, $0<\delta\leq1$ with $\sigma\delta\neq1$ and $t\leq-\mathcal{T}_1$,
		\begin{align}\label{PAS4}
			&\lim_{t\to-\infty}e^{ct}\|\mathcal{K}(\mathfrak{s} +t,\vartheta_{t}\omega)\|^2_{\H}\nonumber\\&\leq e^{\frac{c_0(2\sigma+1)}{\sigma}}\lim_{t\to-\infty}e^{\frac{c}{2}t}\int_{-\infty}^{0} e^{\gamma\xi}\bigg\{\frac{4}{\min\{\mu,\alpha\}}\|\f(\cdot,\xi+\mathfrak{s}+t)\|^2_{\V'} +2s_3e^{\sigma (\xi+\mathfrak{s}+t)}\left|\mathcal{Z}_{\delta}(\vartheta_{\xi+t}\omega)\right|\nonumber\\&\qquad\qquad\qquad\qquad+2s_6 e^{2\sigma (\xi+\mathfrak{s}+t)}\left|\mathcal{Z}_{\delta}(\vartheta_{\xi+t}\omega)\right|^2+s_7\left[e^{\sigma(\xi+\mathfrak{s}+t)}\left|\mathcal{Z}_{\delta}(\vartheta_{\xi+t}\omega)\right|\right]^{\frac{2}{1-s_5}}\bigg\}\d\xi\nonumber\\&\leq4\frac{e^{\frac{c_0(2\sigma+1)}{\sigma}-\frac{c}{2}\mathfrak{s}}}{\min\{\mu,\alpha\}}\lim_{t\to-\infty}e^{\frac{c}{2}t}\int_{-\infty}^{0} e^{\gamma\xi}\|\f(\cdot,\xi+t)\|^2_{\V'}\d\xi\nonumber\\&\quad+ e^{\frac{c_0(2\sigma+1)}{\sigma}}\lim_{t\to-\infty}\int_{-\infty}^{0}e^{c_1t}  \bigg\{2s_3e^{(\gamma+\sigma)\xi}e^{\sigma (\mathfrak{s}+t)}\left|\mathcal{Z}_{\delta}(\vartheta_{\xi+t}\omega)\right|+2s_6 e^{(\gamma+2\sigma)\xi}e^{2\sigma (\mathfrak{s}+t)}\left|\mathcal{Z}_{\delta}(\vartheta_{\xi+t}\omega)\right|^2\nonumber\\&\qquad\qquad\qquad\qquad+s_7e^{\left(\gamma+\frac{2\sigma}{1-s_5}\right)\xi}\left[e^{\sigma(\mathfrak{s}+t)}\left|\mathcal{Z}_{\delta}(\vartheta_{\xi+t}\omega)\right|\right]^{\frac{2}{1-s_5}}\bigg\}\d\xi\nonumber\\&\leq4\frac{e^{\frac{c_0(2\sigma+1)}{\sigma}-\frac{c}{2}\mathfrak{s}}}{\min\{\mu,\alpha\}}\lim_{t\to-\infty}e^{\frac{c}{2}t}\int_{-\infty}^{0} e^{\gamma\xi}\|\f(\cdot,\xi+t)\|^2_{\V'}\d\xi\nonumber\\&\quad+ e^{\frac{c_0(2\sigma+1)}{\sigma}}\lim_{t\to-\infty}\bigg[2s_3e^{\sigma (\mathfrak{s}+t)}\int_{-\infty}^{0}e^{c_1(\xi+t)}  \left|\mathcal{Z}_{\delta}(\vartheta_{\xi+t}\omega)\right|\d\xi\nonumber\\&\quad+2s_6e^{2\sigma (\mathfrak{s}+t)}\int_{-\infty}^{0}e^{c_1(\xi+t)}\left|\mathcal{Z}_{\delta}(\vartheta_{\xi+t}\omega)\right|^2\d\xi+s_7e^{\frac{2\sigma}{1-s_5} (\mathfrak{s}+t)}\int_{-\infty}^{0}e^{c_1(\xi+t)}\left|\mathcal{Z}_{\delta}(\vartheta_{\xi+t}\omega)\right|^{\frac{2}{1-s_5}}\d\xi\bigg]\nonumber\\&\leq4\frac{e^{\frac{c_0(2\sigma+1)}{\sigma}-\frac{c}{2}\mathfrak{s}}}{\min\{\mu,\alpha\}}\lim_{t\to-\infty}e^{\frac{c}{2}t}\int_{-\infty}^{0} e^{\gamma\xi}\|\f(\cdot,\xi+t)\|^2_{\V'}\d\xi\nonumber\\&\quad+ e^{\frac{c_0(2\sigma+1)}{\sigma}}\lim_{t\to-\infty}\bigg[2s_3e^{\sigma (\mathfrak{s}+t)}\int_{-\infty}^{t}e^{c_1\xi}  \left|\mathcal{Z}_{\delta}(\vartheta_{\xi}\omega)\right|\d\xi+2s_6e^{2\sigma (\mathfrak{s}+t)}\int_{-\infty}^{t}e^{c_1\xi}\left|\mathcal{Z}_{\delta}(\vartheta_{\xi}\omega)\right|^2\d\xi\nonumber\\&\qquad\qquad\qquad\qquad\qquad+s_7e^{\frac{2\sigma}{1-s_5} (\mathfrak{s}+t)}\int_{-\infty}^{t}e^{c_1\xi}\left|\mathcal{Z}_{\delta}(\vartheta_{\xi}\omega)\right|^{\frac{2}{1-s_5}}\d\xi\bigg].
		\end{align}
		Using \eqref{N3}, one can easily find that 
		\begin{align}\label{PAS5}
			&\int_{-\infty}^{0}e^{c_1\xi}  \left\{\left|\mathcal{Z}_{\delta}(\vartheta_{\xi}\omega)\right|+\left|\mathcal{Z}_{\delta}(\vartheta_{\xi}\omega)\right|^2+\left|\mathcal{Z}_{\delta}(\vartheta_{\xi}\omega)\right|^{\frac{2}{1-s_5}}\right\}\d\xi <\infty.
		\end{align}
		Taking \eqref{forcing2} and \eqref{PAS5} into account, along with \eqref{PAS4}, we find that 
		\begin{align*}
			\lim_{t\to-\infty}e^{ct}\|\mathcal{K}(\mathfrak{s} +t,\vartheta_{t}\omega)\|^2_{\H}=0.
		\end{align*}
		Also, in the case when $\sigma=0$ or $\kappa=0$ or $\sigma\delta=1$, the above convergence holds true. In fact, in these particular cases, the proof is easier.
	\end{proof}
	
	\begin{theorem}\label{WZ_RA_B}
		For $d=2$ with $r\geq1$, $d=3$ with $r>3$ and $d=r=3$ with $2\beta\mu>1$, assume that $\f\in\mathrm{L}^2_{\emph{loc}}(\R;\V')$ satisfies \eqref{forcing2} and Assumption \ref{NDT1} is fulfilled. Then there exists a unique random $\mathfrak{D}$-pullback attractor $$\mathscr{A}=\{\mathscr{A}(\mathfrak{s},\omega):\mathfrak{s}\in\R, \omega\in\Omega\}\in\mathfrak{D},$$ for the  continuous cocycle $\Phi$ associated with the system \eqref{WZ_SCBF} in $\H$.
	\end{theorem}
	\begin{proof}
		The proof follows from Corollary \ref{Asymptotic_B}, Lemma \ref{PAS} and the abstract theory given in \cite{SandN_Wang} (see Theorem 2.23 in \cite{SandN_Wang}).
	\end{proof}
	
	\section{Random pullback attractors for Wong-Zakai approximations: whole domain $\R^d$} \label{sec4}\setcounter{equation}{0}
	In this section, we prove the existence of a unique random $\mathfrak{D}$-pullback attractor for the system  \eqref{WZ_SCBF} on $\R^d$ with nonlinear diffusion term satisfying the assumptions given below.
	\subsection{Nonlinear diffusion term}
	In this work, we prove the existence of random pullback attractors for Wong-Zakai approximations of SCBF equations on $\R^d$ under two different assumptions on the nonlinear diffusion term $S(t,x,\u)$ appearing in \eqref{WZ_SCBF}. Those two different assumptions are as follows:
	\begin{assumption}\label{NDT2}
		We assume that the nonlinear diffusion term $$S(t,x,\u)=e^{\sigma t}\left[\kappa\u+\mathcal{S}(\u)+\h(x)\right],$$ where $\sigma\geq0, \kappa\geq0$ and $\h\in\H$. Also, $\mathcal{S}:\V\to\H$ is a continuous function satisfying \eqref{S2}-\eqref{S3} and
		\begin{align}\label{S1}
			\left(\mathcal{S}(\u),\u\right)=0,  \  \text{ for all }\  \u\in\V.
		\end{align}
		Note that the condition \eqref{S4} is much weaker than \eqref{S1} and for $s_3=s_4=0$ in \eqref{S4}, one can obtain \eqref{S1}.
	\end{assumption}
	\begin{remark}
		One can take $\h\in\widetilde{\L}^{\frac{r+1}{r}}$ also.
	\end{remark}
	\begin{assumption}\label{NDT3}
		Let $S:\R\times\R^d\times\R^d\to\R^d$ be a continuous function such that for all $t\in\R$ and $x,y\in\R^d$,
		\begin{align}
			|S(t,x,y)|&\leq \mathcal{S}_1(t,x)|y|^{q-1}+\mathcal{S}_2(t,x),\label{GS1}
		\end{align}
		where $1\leq q<r+1$, $\mathcal{S}_1\in\mathrm{L}^{\infty}_{\emph{loc}}(\R;\L^{\frac{r+1}{r+1-q}}(\R^d))$ and $\mathcal{S}_2\in\mathrm{L}^{\infty}_{\emph{loc}}(\R;\L^{\frac{r+1}{r}}(\R^d))$. Furthermore, we assume that $S(t,x,y)$ is locally Lipschitz continuous with respect to the third variable.
	\end{assumption}
	\begin{remark}
		One can take $\mathcal{S}_2\in\mathrm{L}^{\infty}_{\emph{loc}}(\R;\L^{2}(\R^d))$ also.
	\end{remark}
	\begin{example}
		Let us discuss some examples for such nonlinear diffusion terms, which satisfy the above Assumptions. 
		\begin{itemize}
			\item [(1)] Let $\mathcal{S}:\V\to\H$ be a nonlinear operator defined by $\mathcal{S}(\u)=\B(\g_2,\u)$ for all $\u\in\V$, where $\g_2$ is a fixed element of $\D(\A)$ and $\mathcal{S}$ satisfies \eqref{S2}-\eqref{S3} and \eqref{S1}, see \cite{GGW}. Hence $S(t,x,\u)=e^{\sigma t}\left[\kappa\u+\B(\g_2,\u)+\h(x)\right]$ satisfies Assumption \ref{NDT2} but not Assumption \ref{NDT3}.
			\item [(2)] Take $S(t,x,\u)=\mathcal{S}_1(t,x)|\textbf{u}|$, where $\mathcal{S}_1\in\mathrm{L}^{\infty}_{\emph{loc}}(\R;\L^{\frac{r+1}{r-1}}(\R^d))$ for $r>1$. Hence $S(t,x,\u)$ satisfies Assumption \ref{NDT3} but not Assumption \ref{NDT2}.
		\end{itemize}
		The above examples show that Assumptions \ref{NDT2} and \ref{NDT3} cover different classes of functions.
	\end{example}
	
	\subsection{Random pullback attractors under Assumption \ref{NDT2}}
	In this subsection, we prove the existence of a unique random $\mathfrak{D}$-pullback attractor under Assumptions \ref{NDT2} and \ref{DNFT1} on the nonlinear diffusion term $S(t,x,\u)$ and deterministic non-autonomous forcing term $\f(x,t)$, respectively.
	
	\begin{lemma}\label{ContinuityUB1}
		For $d=2$ with $r\geq1$, $d=3$ with $r>3$ and $d=r=3$ with $2\beta\mu>1$, assume that $\f\in \mathrm{L}^2_{\emph{loc}}(\mathbb{R};\V')$ and Assumption \ref{NDT2} is fulfilled. Then, the solution of \eqref{WZ_SCBF} is continuous in initial data $\u_{\mathfrak{s}}(x).$
	\end{lemma}
	\begin{proof}
		See Lemma \ref{Continuity}.
	\end{proof}
	\subsubsection{Uniform estimates and $\mathfrak{D}$-pullback asymptotic compactness of solutions}
	We start by proving uniform estimates for the solutions of \eqref{WZ_SCBF}. To prove the asymptotic compactness, the method of energy equations is introduced in \cite{Ball}. Due to the lack of compact Sobolev embeddings in unbounded domains, we prove $\mathfrak{D}$-pullback asymptotic compactness of solutions of \eqref{WZ_SCBF} on unbounded domains using the idea given in \cite{Ball}. The following lemma is a particular case of Lemma \ref{LemmaUe}.
	\begin{lemma}\label{LemmaUe1}
		For $d=2$ with $r\geq1$, $d=3$ with $r>3$ and $d=r=3$ with $2\beta\mu\geq1$, assume that $\f\in \mathrm{L}^2_{\emph{loc}}(\mathbb{R};\V')$ satisfies \eqref{forcing1} and Assumption \ref{NDT2} is fulfilled. Then for every $0<\delta\leq1$, $\mathfrak{s}\in\R,$ $ \omega\in \Omega$ and $D=\{D(\mathfrak{s},\omega):\mathfrak{s}\in\R, \omega\in\Omega\}\in\mathfrak{D},$ there exists $\mathcal{T}=\mathcal{T}(\delta, \mathfrak{s}, \omega, D)>0$ such that for all $t\geq \mathcal{T}$ and $\tau\geq \mathfrak{s}-t$, the solution $\u$ of the system \eqref{WZ_SCBF} with $\omega$ replaced by $\vartheta_{-\mathfrak{s}}\omega$ satisfies 
		\begin{align}\label{ue^1}
			&\|\u(\tau,\mathfrak{s}-t,\vartheta_{-\mathfrak{s}}\omega,\u_{\mathfrak{s}-t})\|^2_{\H} \nonumber\\&\leq\frac{4}{\min\{\mu,\alpha\}} \int_{-\infty}^{\tau-\mathfrak{s}} e^{\int_{\tau-\mathfrak{s}}^{\xi}\left(\alpha-2\kappa e^{\sigma (\zeta+\mathfrak{s})}\mathcal{Z}_{\delta}(\vartheta_{\zeta}\omega)\right)\d\zeta} \|\f(\cdot,\xi+\mathfrak{s})\|^2_{\V'}\d \xi\nonumber\\&\quad+2s_6\int_{-\infty}^{\tau-\mathfrak{s}} e^{2\sigma (\xi+\mathfrak{s})}e^{\int_{\tau-\mathfrak{s}}^{\xi}\left(\alpha-2\kappa e^{\sigma (\zeta+\mathfrak{s})}\mathcal{Z}_{\delta}(\vartheta_{\zeta}\omega)\right)\d\zeta} \left|\mathcal{Z}_{\delta}(\vartheta_{\xi}\omega)\right|^2\d\xi,
		\end{align}
		where $\u_{\mathfrak{s}-t}\in D(\mathfrak{s}-t,\vartheta_{-t}\omega)$ and $s_6=\frac{2}{\alpha}\|\h\|^2_{\H}$.
	\end{lemma}
	\begin{proof}
		When $s_3=s_4=0$, assumptions \eqref{S1} and \eqref{S4} are the same. Hence, the proof is immediate from Lemma \ref{LemmaUe} by putting $s_3=s_4=0$.
	\end{proof}
	The following Lemma is a direct consequence of Lemma \ref{LemmaUe1}.
	\begin{lemma}\label{LemmaUe2}
		For $d=2$ with $r\geq1$, $d=3$ with $r>3$ and $d=r=3$ with $2\beta\mu\geq1$, assume that $\f\in \mathrm{L}^2_{\emph{loc}}(\mathbb{R};\V')$ satisfies \eqref{forcing1} and Assumption \ref{NDT2} is fulfilled. Then for every $0<\delta\leq1$, $\mathfrak{s}\in\R,$ $ \omega\in \Omega$ and $D=\{D(\mathfrak{s},\omega):\mathfrak{s}\in\R, \omega\in\Omega\}\in\mathfrak{D},$ there exists $\mathcal{T}=\mathcal{T}(\delta,\mathfrak{s},\omega,D)>0$ such that for every $l\geq 0$ and for all $t\geq \mathcal{T}+l$, the solution $\u$ of the system \eqref{WZ_SCBF} with $\omega$ replaced by $\vartheta_{-\mathfrak{s}}\omega$ satisfies 
		\begin{align}\label{ue^2}
			&\|\u(\mathfrak{s}-l,\mathfrak{s}-t,\vartheta_{-\mathfrak{s}}\omega,\u_{\mathfrak{s}-t})\|^2_{\H} \nonumber\\&\leq\frac{4}{\min\{\mu,\alpha\}} \int_{-\infty}^{-l} e^{\int_{-l}^{\xi}\left(\alpha-2\kappa e^{\sigma (\zeta+\mathfrak{s})}\mathcal{Z}_{\delta}(\vartheta_{\zeta}\omega)\right)\d\zeta} \|\f(\cdot,\xi+\mathfrak{s})\|^2_{\V'}\d \xi\nonumber\\&\quad+2s_6\int_{-\infty}^{-l} e^{2\sigma (\xi+\mathfrak{s})}e^{\int_{-l}^{\xi}\left(\alpha-2\kappa e^{\sigma (\zeta+\mathfrak{s})}\mathcal{Z}_{\delta}(\vartheta_{\zeta}\omega)\right)\d\zeta} \left|\mathcal{Z}_{\delta}(\vartheta_{\xi}\omega)\right|^2\d\xi,
		\end{align}
		where $\u_{\mathfrak{s}-t}\in D(\mathfrak{s}-t,\vartheta_{-t}\omega)$ and $s_6=\frac{2}{\alpha}\|\h\|^2_{\H}$.
	\end{lemma} 
	\begin{proof}
		Given $\mathfrak{s}\in \R$ and $l\geq 0$, let $\tau=\mathfrak{s}-l$. Let $\mathcal{T}>0$ be the constant claimed in Lemma \ref{LemmaUe1}. Also, $t\geq \mathcal{T}+l$ implies $t\geq \mathcal{T}$ and $\tau\geq \mathfrak{s}-t$. Thus, the desired result is immediate from Lemma \ref{LemmaUe1}. 
	\end{proof}

	The following Lemma is a result on weak continuity of the solutions of \eqref{WZ_SCBF}, which helps us to prove the $\mathfrak{D}$-pullback asymptotic compactness of the solutions of system \eqref{WZ_SCBF}.
	\begin{lemma}\label{weak}
		For $d=2$ with $r\geq1$, $d=3$ with $r>3$ and $d=r=3$ with $2\beta\mu\geq1$, assume that $\f\in\mathrm{L}^2_{\emph{loc}}(\R;\V')$ and Assumption \ref{NDT2} is fulfilled. Let $0<\delta\leq1$, $\mathfrak{s}\in\R, \omega\in \Omega$ and $\u_{\mathfrak{s}}^0, \u_{\mathfrak{s}}^n\in \H$ for all $n\in\N.$ If $\u_{\mathfrak{s}}^n\xrightharpoonup{w}\u_{\mathfrak{s}}^0$ in $\H$, then the solution $\u$ of the system \eqref{WZ_SCBF} satisfies the following convergences:
		\begin{itemize}
			\item [(i)] $\u(\xi,\mathfrak{s},\omega,\u_{\mathfrak{s}}^n)\xrightharpoonup{w}\u(\xi,\mathfrak{s},\omega,\u_{\mathfrak{s}}^0)$ in $\H$ for all $\xi\geq \mathfrak{s}$.
			\item [(ii)] $\u(\cdot,\mathfrak{s},\omega,\u_{\mathfrak{s}}^n)\xrightharpoonup{w}\u(\cdot,\mathfrak{s},\omega,\u_{\mathfrak{s}}^0)$ in $\mathrm{L}^2((\mathfrak{s},\mathfrak{s}+T);\V)$ for every $T>0$.
			\item [(iii)] $\u(\cdot,\mathfrak{s},\omega,\u_{\mathfrak{s}}^n)\xrightharpoonup{w}\u(\cdot,\mathfrak{s},\omega,\u_{\mathfrak{s}}^0)$ in $\mathrm{L}^{r+1}((\mathfrak{s},\mathfrak{s}+T);\widetilde{\L}^{r+1})$ for every $T>0$.
		\end{itemize}
	\end{lemma}
	\begin{proof}
		Using a standard method as in \cite{KM1} (see Lemmas 5.2 and 5.3 in \cite{KM1}), one can complete the proof.
	\end{proof}
	Now, we establish the $\mathfrak{D}$-pullback asymptotic compactness of the solutions of the system \eqref{WZ_SCBF}.
	\begin{lemma}\label{Asymptotic_UB}
		For $d=2$ with $r\geq1$, $d=3$ with $r>3$ and $d=r=3$ with $2\beta\mu\geq1$, assume that $\f\in\mathrm{L}^2_{\emph{loc}}(\R;\V')$ satisfies \eqref{forcing1} and Assumption \ref{NDT2} is fulfilled. Then for every $0<\delta\leq1$, $\mathfrak{s}\in \R,$ $\omega\in \Omega,$ $D=\{D(\mathfrak{s},\omega):\mathfrak{s}\in \R,\omega\in \Omega\}\in \mathfrak{D}$ and $t_n\to \infty,$ $\u_{0,n}\in D(\mathfrak{s}-t_n, \vartheta_{-t_{n}}\omega)$, the sequence $\Phi(t_n,\mathfrak{s}-t_n,\vartheta_{-t_n}\omega,\u_{0,n})$ or $\u(\mathfrak{s},\mathfrak{s}-t_n,\vartheta_{-\mathfrak{s}}\omega,\u_{0,n})$ of solutions of the system \eqref{WZ_SCBF} has a convergent subsequence in $\H$.
	\end{lemma}
	\begin{proof}
		It follows from Lemma \ref{LemmaUe2} with $l=0$ that there exists $\mathcal{T}=\mathcal{T}(\delta,\mathfrak{s},\omega,D)>0$ such that for all $t\geq \mathcal{T}$,
		\begin{align}\label{ac1}
			\|\u(\mathfrak{s},\mathfrak{s}-t,\vartheta_{-\mathfrak{s}}\omega,\u_{\mathfrak{s}-t})\|^2_{\H} &\leq\frac{4}{\min\{\mu,\alpha\}} \int_{-\infty}^{0} e^{\int_{0}^{\xi}\left(\alpha-2\kappa e^{\sigma (\zeta+\mathfrak{s})}\mathcal{Z}_{\delta}(\vartheta_{\zeta}\omega)\right)\d\zeta} \|\f(\cdot,\xi+\mathfrak{s})\|^2_{\V'}\d \xi\nonumber\\&\quad+2s_6\int_{-\infty}^{0} e^{2\sigma (\xi+\mathfrak{s})}e^{\int_{0}^{\xi}\left(\alpha-2\kappa e^{\sigma (\zeta+\mathfrak{s})}\mathcal{Z}_{\delta}(\vartheta_{\zeta}\omega)\right)\d\zeta} \left|\mathcal{Z}_{\delta}(\vartheta_{\xi}\omega)\right|^2\d\xi\nonumber\\&=:M(\mathfrak{s},\omega),
		\end{align}
		where $\u_{\mathfrak{s}-t}\in D(\mathfrak{s}-t,\vartheta_{-t}\omega).$ Since $t_n\to \infty$, there exists $N_0\in\N$ such that $t_n\geq \mathcal{T}$ for all $n\geq N_0$. As, it is given that $\u_{0,n}\in D(\mathfrak{s}-t_n, \vartheta_{-t_{n}}\omega)$, \eqref{ac1} implies that for all $n\geq N_0$, 
		\begin{align*}
			\|\u(\mathfrak{s},\mathfrak{s}-t_n,\vartheta_{-\mathfrak{s}}\omega,\u_{0,n})\|^2_{\H} \leq M(\mathfrak{s},\omega),
		\end{align*}
		and hence $\{\u(\mathfrak{s},\mathfrak{s}-t_n,\vartheta_{-\mathfrak{s}}\omega,\u_{0,n})\}_{n\geq N_0}\subseteq\H$ is a bounded sequence, which implies that there exists $\tilde{\u}\in \H$ and a subsequence (keeping same label) such that 
		\begin{align}\label{ac2}
			\u(\mathfrak{s},\mathfrak{s}-t_n,\vartheta_{-\mathfrak{s}}\omega,\u_{0,n})\xrightharpoonup{w}\tilde{\u}\  \text{ in }\  \H.
		\end{align}
		By the weak lower semicontinuous property of norms and \eqref{ac2}, we get
		\begin{align}\label{ac3}
			\|\tilde{\u}\|_{\H}\leq\liminf_{n\to\infty}\|\u(\mathfrak{s},\mathfrak{s}-t_n,\vartheta_{-\mathfrak{s}}\omega,\u_{0,n})\|_{\H}.
		\end{align}
		In order to get the desired result, we have to prove that $\u(\mathfrak{s},\mathfrak{s}-t_n,\vartheta_{-\mathfrak{s}}\omega,\u_{0,n})\to\tilde{\u}$ in $\H$ strongly, that is, we only need to show that
		\begin{align}\label{ac4}
			\|\tilde{\u}\|_{\H}\geq\limsup_{n\to\infty}\|\u(\mathfrak{s},\mathfrak{s}-t_n,\vartheta_{-\mathfrak{s}}\omega,\u_{0,n})\|_{\H}.
		\end{align}
		The method of energy equations  introduced in \cite{Ball} will help us to prove \eqref{ac4}. For a given $l\in \N$ ($l\leq t_n$), we can write 
		\begin{align}\label{ac5}
			\u(\mathfrak{s},\mathfrak{s}-t_n,\vartheta_{-\mathfrak{s}}\omega,\u_{0,n})=\u(\mathfrak{s},\mathfrak{s}-l,\vartheta_{-\mathfrak{s}}\omega,\u(\mathfrak{s}-l,\mathfrak{s}-t_n,\vartheta_{-\mathfrak{s}}\omega,\u_{0,n})).
		\end{align} 
		For each $l$, let $N_l$ be sufficiently large such that $t_n\geq \mathcal{T}+l$ for all $n\geq N_l$. By Lemma \ref{LemmaUe2}, we have for $l\geq N_l$,
		\begin{align*}
			&\|\u(\mathfrak{s}-l,\mathfrak{s}-t_n,\vartheta_{-\mathfrak{s}}\omega,\u_{0,n})\|^2_{\H} \nonumber\\&\leq\frac{4}{\min\{\mu,\alpha\}} \int_{-\infty}^{-l} e^{\int_{-l}^{\xi}\left(\alpha-2\kappa e^{\sigma (\zeta+\mathfrak{s})}\mathcal{Z}_{\delta}(\vartheta_{\zeta}\omega)\right)\d\zeta} \|\f(\cdot,\xi+\mathfrak{s})\|^2_{\V'}\d \xi\nonumber\\&\quad+2s_6\int_{-\infty}^{-l} e^{2\sigma (\xi+\mathfrak{s})}e^{\int_{-l}^{\xi}\left(\alpha-2\kappa e^{\sigma (\zeta+\mathfrak{s})}\mathcal{Z}_{\delta}(\vartheta_{\zeta}\omega)\right)\d\zeta} \left|\mathcal{Z}_{\delta}(\vartheta_{\xi}\omega)\right|^2\d\xi,
		\end{align*}
		which infer that the sequence $\{\u(\mathfrak{s}-l,\mathfrak{s}-t_n,\vartheta_{-\mathfrak{s}}\omega,\u_{0,n})\}_{n\geq N_l}$ is bounded in $\H$,  for each $l\in \N$. By the diagonal process, there exists a subsequence (denoting by same symbol for convenience) and an element $\tilde{\u}_{l}\in \H$ for each $l\in\N$ such that 
		\begin{align}\label{ac6}
			\u(\mathfrak{s}-l,\mathfrak{s}-t_n,\vartheta_{-\mathfrak{s}}\omega,\u_{0,n})\xrightharpoonup{w}\tilde{\u}_{l} \ \text{ in }\  \H.
		\end{align}
		From \eqref{ac5}-\eqref{ac6} along with Lemma \ref{weak}, we obtain that for $l\in\N$,
		\begin{align}\label{ac7}
			\u(\mathfrak{s},\mathfrak{s}-t_n,\vartheta_{-\mathfrak{s}}\omega,\u_{0,n})\xrightharpoonup{w}\u(\mathfrak{s},\mathfrak{s}-l,\vartheta_{-\mathfrak{s}}\omega,\tilde{\u}_{l}) \ \text{ in } \ \H,
		\end{align}
		\begin{align}\label{ac8}
			\u(\cdot,\mathfrak{s}-l,\vartheta_{-\mathfrak{s}}\omega,\u(\mathfrak{s}-l,\mathfrak{s}-t_n,\vartheta_{-\mathfrak{s}}\omega,\u_{0,n}))\xrightharpoonup{w}\u(\cdot,\mathfrak{s}-l,\vartheta_{-\mathfrak{s}}\omega,\tilde{\u}_{l}) \text{ in } \mathrm{L}^2((\mathfrak{s}-l,\mathfrak{s});\V),
		\end{align}
		and
		\begin{align}\label{ac8'}
			\u(\cdot,\mathfrak{s}-l,\vartheta_{-\mathfrak{s}}\omega,\u(\mathfrak{s}-l,\mathfrak{s}-t_n,\vartheta_{-\mathfrak{s}}\omega,\u_{0,n}))\xrightharpoonup{w}\u(\cdot,\mathfrak{s}-l,\vartheta_{-\mathfrak{s}}\omega,\tilde{\u}_{l}) \text{ in } \mathrm{L}^{r+1}((\mathfrak{s}-l,\mathfrak{s});\widetilde{\L}^{r+1}).
		\end{align}
		Clearly, \eqref{ac2} and \eqref{ac7} imply that 
		\begin{align}\label{ac9}
			\u(\mathfrak{s},\mathfrak{s}-l,\vartheta_{-\mathfrak{s}}\omega,\tilde{\u}_{l})=\tilde{\u}.
		\end{align}
		Taking into account that $\left(\mathcal{S}(\u),\u\right)=0$, from \eqref{ue0} we have
		\begin{align}\label{ac10}
			&	\frac{\d}{\d t} \|\u\|^2_{\H} + \left(\alpha-2\kappa e^{\sigma t}\mathcal{Z}_{\delta}(\vartheta_{t}\omega)\right)\|\u\|^2_{\H}\nonumber\\&=-2\mu\|\nabla\u\|^2_{\H}-\alpha\|\u\|^2_{\H} - 2\beta\|\u\|^{r+1}_{\wi \L^{r+1}} + 2\left\langle\f,\u\right\rangle +2e^{\sigma t}\mathcal{Z}_{\delta}(\vartheta_{t}\omega)\left(\h,\u\right).
		\end{align}
		It follows by applying variation of constant formula to \eqref{ac10}, that for each $\omega\in \Omega,$ $ \tau\in \R$ and $\mathfrak{s}\geq \tau$,
		\begin{align}\label{ac11}
			\|\u(\mathfrak{s},\tau,\omega,\u_{\tau})\|^2_{\H} &= e^{\int_{\mathfrak{s}}^{\tau}\left(\alpha-2\kappa e^{\sigma\zeta}\mathcal{Z}_{\delta}(\vartheta_{\zeta}\omega)\right)\d\zeta}\|\u_{\tau}\|^2_{\H} \nonumber\\&\quad-2\mu\int_{\tau}^{\mathfrak{s}}e^{\int_{\mathfrak{s}}^{\xi}\left(\alpha-2\kappa e^{\sigma\zeta}\mathcal{Z}_{\delta}(\vartheta_{\zeta}\omega)\right)\d\zeta}\|\nabla\u(\xi,\tau,\omega,\u_{\tau})\|^2_{\H}\d\xi\nonumber\\&\quad-\alpha\int_{\tau}^{\mathfrak{s}}e^{\int_{\mathfrak{s}}^{\xi}\left(\alpha-2\kappa e^{\sigma\zeta}\mathcal{Z}_{\delta}(\vartheta_{\zeta}\omega)\right)\d\zeta}\|\u(\xi,\tau,\omega,\u_{\tau})\|^2_{\H}\d\xi\nonumber\\&\quad-2\beta\int_{\tau}^{\mathfrak{s}}e^{\int_{\mathfrak{s}}^{\xi}\left(\alpha-2\kappa e^{\sigma\zeta}\mathcal{Z}_{\delta}(\vartheta_{\zeta}\omega)\right)\d\zeta}\|\u(\xi,\tau,\omega,\u_{\tau})\|^{r+1}_{\wi \L^{r+1}}\d\xi\nonumber\\&\quad+2\int_{\tau}^{\mathfrak{s}}e^{\int_{\mathfrak{s}}^{\xi}\left(\alpha-2\kappa e^{\sigma\zeta}\mathcal{Z}_{\delta}(\vartheta_{\zeta}\omega)\right)\d\zeta}\left\langle\f(\cdot,\xi),\u(\xi,\tau,\omega,\u_{\tau})\right\rangle\d\xi\nonumber\\&\quad+2\int_{\tau}^{\mathfrak{s}}e^{\sigma\xi}e^{\int_{\mathfrak{s}}^{\xi}\left(\alpha-2\kappa e^{\sigma\zeta}\mathcal{Z}_{\delta}(\vartheta_{\zeta}\omega)\right)\d\zeta}\mathcal{Z}_{\delta}(\vartheta_{\xi}\omega)\left(\h,\u(\xi,\tau,\omega,\u_{\tau})\right)\d\xi.
		\end{align}
		From \eqref{ac9} and \eqref{ac11}, it is immediate that 
		\begin{align}\label{ac12}
			\|\tilde{\u}\|^2_{\H}&=\|\u(\mathfrak{s},\mathfrak{s}-l,\vartheta_{-\mathfrak{s}}\omega,\tilde{\u}_{l})\|^2_{\H} \nonumber\\&= e^{\int_{0}^{-l}\left(\alpha-2\kappa e^{\sigma(\zeta+\mathfrak{s})}\mathcal{Z}_{\delta}(\vartheta_{\zeta}\omega)\right)\d\zeta}\|\tilde{\u}_{l}\|^2_{\H} \nonumber\\&\quad-2\mu\int_{-l}^{0}e^{\int_{0}^{\xi}\left(\alpha-2\kappa e^{\sigma(\zeta+\mathfrak{s})}\mathcal{Z}_{\delta}(\vartheta_{\zeta}\omega)\right)\d\zeta}\|\nabla\u(\xi+\mathfrak{s},\mathfrak{s}-l,\vartheta_{-\mathfrak{s}}\omega,\tilde{\u}_{l})\|^2_{\H}\d\xi\nonumber\\&\quad-\alpha\int_{-l}^{0}e^{\int_{0}^{\xi}\left(\alpha-2\kappa e^{\sigma(\zeta+\mathfrak{s})}\mathcal{Z}_{\delta}(\vartheta_{\zeta}\omega)\right)\d\zeta}\|\u(\xi+\mathfrak{s},\mathfrak{s}-l,\vartheta_{-\mathfrak{s}}\omega,\tilde{\u}_{l})\|^2_{\H}\d\xi\nonumber\\&\quad-2\beta\int_{-l}^{0}e^{\int_{0}^{\xi}\left(\alpha-2\kappa e^{\sigma(\zeta+\mathfrak{s})}\mathcal{Z}_{\delta}(\vartheta_{\zeta}\omega)\right)\d\zeta}\|\u(\xi+\mathfrak{s},\mathfrak{s}-l,\vartheta_{-\mathfrak{s}}\omega,\tilde{\u}_{l})\|^{r+1}_{\wi \L^{r+1}}\d\xi\nonumber\\&\quad+2\int_{-l}^{0}e^{\int_{0}^{\xi}\left(\alpha-2\kappa e^{\sigma(\zeta+\mathfrak{s})}\mathcal{Z}_{\delta}(\vartheta_{\zeta}\omega)\right)\d\zeta}\left\langle\f(\cdot,\xi+\mathfrak{s}),\u(\xi+\mathfrak{s},\mathfrak{s}-l,\vartheta_{-\mathfrak{s}}\omega,\tilde{\u}_{l})\right\rangle\d\xi\nonumber\\&\quad+2\int_{-l}^{0}e^{\sigma(\xi+\mathfrak{s})}e^{\int_{0}^{\xi}\left(\alpha-2\kappa e^{\sigma(\zeta+\mathfrak{s})}\mathcal{Z}_{\delta}(\vartheta_{\zeta}\omega)\right)\d\zeta}\mathcal{Z}_{\delta}(\vartheta_{\xi}\omega)\left(\h,\u(\xi+\mathfrak{s},\mathfrak{s}-l,\vartheta_{-\mathfrak{s}}\omega,\tilde{\u}_{l})\right)\d\xi.
		\end{align}
		Similarly, from \eqref{ac5} and \eqref{ac11}, we obtain
		\begin{align}\label{ac13}
			&\|\u(\mathfrak{s},\mathfrak{s}-t_n,\vartheta_{-\mathfrak{s}}\omega,\u_{0,n})\|^2_{\H}\nonumber\\&=\|\u(\mathfrak{s},\mathfrak{s}-l,\vartheta_{-\mathfrak{s}}\omega,\u(\mathfrak{s}-l,\mathfrak{s}-t_n,\vartheta_{-\mathfrak{s}}\omega,\u_{0,n}))\|^2_{\H} \nonumber\\&= e^{\int_{0}^{-l}\left(\alpha-2\kappa e^{\sigma(\zeta+\mathfrak{s})}\mathcal{Z}_{\delta}(\vartheta_{\zeta}\omega)\right)\d\zeta}\|\u(\mathfrak{s}-l,\mathfrak{s}-t_n,\vartheta_{-\mathfrak{s}}\omega,\u_{0,n})\|^2_{\H} \nonumber\\&\quad-2\mu\int_{-l}^{0}e^{\int_{0}^{-l}\left(\alpha-2\kappa e^{\sigma(\zeta+\mathfrak{s})}\mathcal{Z}_{\delta}(\vartheta_{\zeta}\omega)\right)\d\zeta}\|\nabla\u(\xi+\mathfrak{s},\mathfrak{s}-l,\vartheta_{-\mathfrak{s}}\omega,\u(\mathfrak{s}-l,\mathfrak{s}-t_n,\vartheta_{-\mathfrak{s}}\omega,\u_{0,n}))\|^2_{\H}\d\xi\nonumber\\&\quad-\alpha\int_{-l}^{0}e^{\int_{0}^{-l}\left(\alpha-2\kappa e^{\sigma(\zeta+\mathfrak{s})}\mathcal{Z}_{\delta}(\vartheta_{\zeta}\omega)\right)\d\zeta}\|\u(\xi+\mathfrak{s},\mathfrak{s}-l,\vartheta_{-\mathfrak{s}}\omega,\u(\mathfrak{s}-l,\mathfrak{s}-t_n,\vartheta_{-\mathfrak{s}}\omega,\u_{0,n}))\|^2_{\H}\d\xi\nonumber\\&\quad-2\beta\int_{-l}^{0}e^{\int_{0}^{-l}\left(\alpha-2\kappa e^{\sigma(\zeta+\mathfrak{s})}\mathcal{Z}_{\delta}(\vartheta_{\zeta}\omega)\right)\d\zeta}\|\u(\xi+\mathfrak{s},\mathfrak{s}-l,\vartheta_{-\mathfrak{s}}\omega,\u(\mathfrak{s}-l,\mathfrak{s}-t_n,\vartheta_{-\mathfrak{s}}\omega,\u_{0,n}))\|^{r+1}_{\wi \L^{r+1}}\d\xi\nonumber\\&\quad+2\int_{-l}^{0}e^{\int_{0}^{-l}\left(\alpha-2\kappa e^{\sigma(\zeta+\mathfrak{s})}\mathcal{Z}_{\delta}(\vartheta_{\zeta}\omega)\right)\d\zeta}\nonumber\\&\qquad\qquad\times\left\langle\f(\cdot,\xi+\mathfrak{s}),\u(\xi+\mathfrak{s},\mathfrak{s}-l,\vartheta_{-\mathfrak{s}}\omega,\u(\mathfrak{s}-l,\mathfrak{s}-t_n,\vartheta_{-\mathfrak{s}}\omega,\u_{0,n}))\right\rangle\d\xi\nonumber\\&\quad+2\int_{-l}^{0}e^{\sigma(\xi+\mathfrak{s})}e^{\int_{0}^{-l}\left(\alpha-2\kappa e^{\sigma(\zeta+\mathfrak{s})}\mathcal{Z}_{\delta}(\vartheta_{\zeta}\omega)\right)\d\zeta}\mathcal{Z}_{\delta}(\vartheta_{\xi}\omega)\nonumber\\&\qquad\qquad\times\left(\h,\u(\xi+\mathfrak{s},\mathfrak{s}-l,\vartheta_{-\mathfrak{s}}\omega,\u(\mathfrak{s}-l,\mathfrak{s}-t_n,\vartheta_{-\mathfrak{s}}\omega,\u_{0,n}))\right)\d\xi.
		\end{align}
		Now, we examine the limits of each term of right-hand side of \eqref{ac13} as $n\to \infty$. By \eqref{ue2}, we examine the first term with $\tau=\mathfrak{s}-l$ and $t=t_n$ as follows
		\begin{align}\label{ac14}
			&e^{\int_{0}^{-l}\left(\alpha-2\kappa e^{\sigma(\zeta+\mathfrak{s})}\mathcal{Z}_{\delta}(\vartheta_{\zeta}\omega)\right)\d\zeta}\|\u(\mathfrak{s}-l,\mathfrak{s}-t_n,\vartheta_{-\mathfrak{s}}\omega,\u_{0,n})\|^2_{\H} \nonumber\\&\leq e^{\int_{0}^{-t_n}\left(\alpha-2\kappa e^{\sigma(\zeta+\mathfrak{s})}\mathcal{Z}_{\delta}(\vartheta_{\zeta}\omega)\right)\d\zeta}\|\u_{0,n}\|^2_{\H}\nonumber\\&\quad+\frac{4}{\min\{\mu,\alpha\}}\int_{-\infty}^{-l} e^{\int_{0}^{\xi}\left(\alpha-2\kappa e^{\sigma (\zeta+\mathfrak{s})}\mathcal{Z}_{\delta}(\vartheta_{\zeta}\omega)\right)\d\zeta} \|\f(\cdot,\xi+\mathfrak{s})\|^2_{\V'}\d \xi\nonumber\\&\quad+s_6\int_{-\infty}^{-l} e^{2\sigma (\xi+\mathfrak{s})} e^{\int_{0}^{\xi}\left(\alpha-2\kappa e^{\sigma (\zeta+\mathfrak{s})}\mathcal{Z}_{\delta}(\vartheta_{\zeta}\omega)\right)\d\zeta}\left|\mathcal{Z}_{\delta}(\vartheta_{\xi}\omega)\right|^2\d\xi.
		\end{align}
		Since $\u_{0,n}\in D(\mathfrak{s}-t_n,\vartheta_{-t_n}\omega)$, we have
		\begin{align}\label{ac15}
			e^{\int_{0}^{-t_n}\left(\alpha-2\kappa e^{\sigma(\zeta+\mathfrak{s})}\mathcal{Z}_{\delta}(\vartheta_{\zeta}\omega)\right)\d\zeta}\|\u_{0,n}\|^2_{\H}\leq e^{\int_{0}^{-t_n}\left(\alpha-2\kappa e^{\sigma(\zeta+\mathfrak{s})}\mathcal{Z}_{\delta}(\vartheta_{\zeta}\omega)\right)\d\zeta}\|D(\mathfrak{s}-t_n,\vartheta_{-t_n}\omega)\|^2_{\H}\to 0,
		\end{align}
		as $n\to \infty$. Combining \eqref{ac15} along with \eqref{ac14}, we have
		\begin{align}\label{ac16}
			&\limsup_{n\to\infty}e^{\int_{0}^{-l}\left(\alpha-2\kappa e^{\sigma(\zeta+\mathfrak{s})}\mathcal{Z}_{\delta}(\vartheta_{\zeta}\omega)\right)\d\zeta}\|\u(\mathfrak{s}-l,\mathfrak{s}-t_n,\vartheta_{-\mathfrak{s}}\omega,\u_{0,n})\|^2_{\H}\nonumber\\& \leq\frac{4}{\min\{\mu,\alpha\}}\int_{-\infty}^{-l} e^{\int_{0}^{\xi}\left(\alpha-2\kappa e^{\sigma (\zeta+\mathfrak{s})}\mathcal{Z}_{\delta}(\vartheta_{\zeta}\omega)\right)\d\zeta} \|\f(\cdot,\xi+\mathfrak{s})\|^2_{\V'}\d \xi\nonumber\\&\quad+s_6\int_{-\infty}^{-l} e^{2\sigma (\xi+\mathfrak{s})} e^{\int_{0}^{\xi}\left(\alpha-2\kappa e^{\sigma (\zeta+\mathfrak{s})}\mathcal{Z}_{\delta}(\vartheta_{\zeta}\omega)\right)\d\zeta}\left|\mathcal{Z}_{\delta}(\vartheta_{\xi}\omega)\right|^2\d\xi.
		\end{align} 
		From \eqref{ac8}, we get
		\begin{align}\label{ac17}
			&	\lim_{n\to\infty} 2\int_{-l}^{0}e^{\int_{0}^{-l}\left(\alpha-2\kappa e^{\sigma(\zeta+\mathfrak{s})}\mathcal{Z}_{\delta}(\vartheta_{\zeta}\omega)\right)\d\zeta}\nonumber\\&\qquad\qquad\times\left\langle\f(\cdot,\xi+\mathfrak{s}),\u(\xi+\mathfrak{s},\mathfrak{s}-l,\vartheta_{-\mathfrak{s}}\omega,\u(\mathfrak{s}-l,\mathfrak{s}-t_n,\vartheta_{-\mathfrak{s}}\omega,\u_{0,n}))\right\rangle\d\xi \nonumber\\&= 2\int_{-l}^{0}e^{\int_{0}^{-l}\left(\alpha-2\kappa e^{\sigma(\zeta+\mathfrak{s})}\mathcal{Z}_{\delta}(\vartheta_{\zeta}\omega)\right)\d\zeta}\left\langle\f(\cdot,\xi+\mathfrak{s}),\u(\xi+\mathfrak{s},\mathfrak{s}-l,\vartheta_{-\mathfrak{s}}\omega,\tilde{\u}_{l})\right\rangle\d\xi,
		\end{align}
		and
		\begin{align}\label{ac18}
			&	\lim_{n\to\infty} 2\int_{-l}^{0}e^{\sigma(\xi+\mathfrak{s})}e^{\int_{0}^{-l}\left(\alpha-2\kappa e^{\sigma(\zeta+\mathfrak{s})}\mathcal{Z}_{\delta}(\vartheta_{\zeta}\omega)\right)\d\zeta}\mathcal{Z}_{\delta}(\vartheta_{\xi}\omega)\nonumber\\&\qquad\qquad\times\left(\h,\u(\xi+\mathfrak{s},\mathfrak{s}-l,\vartheta_{-\mathfrak{s}}\omega,\u(\mathfrak{s}-l,\mathfrak{s}-t_n,\vartheta_{-\mathfrak{s}}\omega,\u_{0,n}))\right)\d\xi \nonumber\\&= 2\int_{-l}^{0}e^{\sigma(\xi+\mathfrak{s})}e^{\int_{0}^{-l}\left(\alpha-2\kappa e^{\sigma(\zeta+\mathfrak{s})}\mathcal{Z}_{\delta}(\vartheta_{\zeta}\omega)\right)\d\zeta}\mathcal{Z}_{\delta}(\vartheta_{\xi}\omega)\left(\h,\u(\xi+\mathfrak{s},\mathfrak{s}-l,\vartheta_{-\mathfrak{s}}\omega,\tilde{\u}_l)\right)\d\xi.
		\end{align}
		Using \eqref{ac8} and the weak lower semicontinuity property of norms, we obtain
		\begin{align*}
			&\liminf_{n\to\infty}	\bigg\{2\mu\int_{-l}^{0}e^{\int_{0}^{-l}\left(\alpha-2\kappa e^{\sigma(\zeta+\mathfrak{s})}\mathcal{Z}_{\delta}(\vartheta_{\zeta}\omega)\right)\d\zeta}\nonumber\\&\qquad\qquad\times\|\nabla\u(\xi+\mathfrak{s},\mathfrak{s}-l,\vartheta_{-\mathfrak{s}}\omega,\u(\mathfrak{s}-l,\mathfrak{s}-t_n,\vartheta_{-\mathfrak{s}}\omega,\u_{0,n}))\|^2_{\H}\d\xi\bigg\} \nonumber\\&\geq 2\mu\int_{-l}^{0}e^{\int_{0}^{-l}\left(\alpha-2\kappa e^{\sigma(\zeta+\mathfrak{s})}\mathcal{Z}_{\delta}(\vartheta_{\zeta}\omega)\right)\d\zeta}\|\nabla\u(\xi+\mathfrak{s},\mathfrak{s}-l,\vartheta_{-\mathfrak{s}}\omega,\tilde{\u}_l)\|^2_{\H}\d\xi,
		\end{align*} 
		or
		\begin{align}\label{ac19}
			&\limsup_{n\to\infty}	\bigg\{-2\mu\int_{-l}^{0}e^{\int_{0}^{-l}\left(\alpha-2\kappa e^{\sigma(\zeta+\mathfrak{s})}\mathcal{Z}_{\delta}(\vartheta_{\zeta}\omega)\right)\d\zeta}\nonumber\\&\qquad\qquad\times\|\nabla\u(\xi+\mathfrak{s},\mathfrak{s}-l,\vartheta_{-\mathfrak{s}}\omega,\u(\mathfrak{s}-l,\mathfrak{s}-t_n,\vartheta_{-\mathfrak{s}}\omega,\u_{0,n}))\|^2_{\H}\d\xi\bigg\} \nonumber\\&\leq -2\mu\int_{-l}^{0}e^{\int_{0}^{-l}\left(\alpha-2\kappa e^{\sigma(\zeta+\mathfrak{s})}\mathcal{Z}_{\delta}(\vartheta_{\zeta}\omega)\right)\d\zeta}\|\nabla\u(\xi+\mathfrak{s},\mathfrak{s}-l,\vartheta_{-\mathfrak{s}}\omega,\tilde{\u}_l)\|^2_{\H}\d\xi.
		\end{align} 
		Similarly, using \eqref{ac8}-\eqref{ac8'} and the weak lower semicontinuity property of norms, we get 
		\begin{align}\label{ac20}
			&\limsup_{n\to\infty}	\bigg\{-\alpha\int_{-l}^{0}e^{\int_{0}^{-l}\left(\alpha-2\kappa e^{\sigma(\zeta+\mathfrak{s})}\mathcal{Z}_{\delta}(\vartheta_{\zeta}\omega)\right)\d\zeta}\nonumber\\&\qquad\qquad\times\|\u(\xi+\mathfrak{s},\mathfrak{s}-l,\vartheta_{-\mathfrak{s}}\omega,\u(\mathfrak{s}-l,\mathfrak{s}-t_n,\vartheta_{-\mathfrak{s}}\omega,\u_{0,n}))\|^2_{\H}\d\xi\bigg\} \nonumber\\&\leq -\alpha\int_{-l}^{0}e^{\int_{0}^{-l}\left(\alpha-2\kappa e^{\sigma(\zeta+\mathfrak{s})}\mathcal{Z}_{\delta}(\vartheta_{\zeta}\omega)\right)\d\zeta}\|\u(\xi+\mathfrak{s},\mathfrak{s}-l,\vartheta_{-\mathfrak{s}}\omega,\tilde{\u}_l)\|^2_{\H}\d\xi,
		\end{align} 
		and
		\begin{align}\label{ac21}
			&\limsup_{n\to\infty}	\bigg\{-2\beta\int_{-l}^{0}e^{\int_{0}^{-l}\left(\alpha-2\kappa e^{\sigma(\zeta+\mathfrak{s})}\mathcal{Z}_{\delta}(\vartheta_{\zeta}\omega)\right)\d\zeta}\nonumber\\&\qquad\qquad\times\|\u(\xi+\mathfrak{s},\mathfrak{s}-l,\vartheta_{-\mathfrak{s}}\omega,\u(\mathfrak{s}-l,\mathfrak{s}-t_n,\vartheta_{-\mathfrak{s}}\omega,\u_{0,n}))\|^{r+1}_{\widetilde{\L}^{r+1}}\d\xi\bigg\} \nonumber\\&\leq -2\beta\int_{-l}^{0}e^{\int_{0}^{-l}\left(\alpha-2\kappa e^{\sigma(\zeta+\mathfrak{s})}\mathcal{Z}_{\delta}(\vartheta_{\zeta}\omega)\right)\d\zeta}\|\u(\xi+\mathfrak{s},\mathfrak{s}-l,\vartheta_{-\mathfrak{s}}\omega,\tilde{\u}_l)\|^{r+1}_{\widetilde{\L}^{r+1}}\d\xi.
		\end{align}
		Combining \eqref{ac16}-\eqref{ac21}, and using it in \eqref{ac13}, we find
		\begin{align}\label{ac22}
			&\limsup_{n\to\infty}\|\u(\mathfrak{s},\mathfrak{s}-t_n,\vartheta_{-\mathfrak{s}}\omega,\u_{0,n})\|^2_{\H}\nonumber\\&\leq\frac{4}{\min\{\mu,\alpha\}}\int_{-\infty}^{-l} e^{\int_{0}^{\xi}\left(\alpha-2\kappa e^{\sigma (\zeta+\mathfrak{s})}\mathcal{Z}_{\delta}(\vartheta_{\zeta}\omega)\right)\d\zeta} \|\f(\cdot,\xi+\mathfrak{s})\|^2_{\V'}\d \xi\nonumber\\&\quad+s_6\int_{-\infty}^{-l} e^{2\sigma (\xi+\mathfrak{s})} e^{\int_{0}^{\xi}\left(\alpha-2\kappa e^{\sigma (\zeta+\mathfrak{s})}\mathcal{Z}_{\delta}(\vartheta_{\zeta}\omega)\right)\d\zeta}\left|\mathcal{Z}_{\delta}(\vartheta_{\xi}\omega)\right|^2\d\xi\nonumber\\&\quad-2\mu\int_{-l}^{0}e^{\int_{0}^{-l}\left(\alpha-2\kappa e^{\sigma(\zeta+\mathfrak{s})}\mathcal{Z}_{\delta}(\vartheta_{\zeta}\omega)\right)\d\zeta}\|\nabla\u(\xi+\mathfrak{s},\mathfrak{s}-l,\vartheta_{-\mathfrak{s}}\omega,\tilde{\u}_l)\|^2_{\H}\d\xi\nonumber\\&\quad-\alpha\int_{-l}^{0}e^{\int_{0}^{-l}\left(\alpha-2\kappa e^{\sigma(\zeta+\mathfrak{s})}\mathcal{Z}_{\delta}(\vartheta_{\zeta}\omega)\right)\d\zeta}\|\u(\xi+\mathfrak{s},\mathfrak{s}-l,\vartheta_{-\mathfrak{s}}\omega,\tilde{\u}_l)\|^2_{\H}\d\xi\nonumber\\&\quad-2\beta\int_{-l}^{0}e^{\int_{0}^{-l}\left(\alpha-2\kappa e^{\sigma(\zeta+\mathfrak{s})}\mathcal{Z}_{\delta}(\vartheta_{\zeta}\omega)\right)\d\zeta}\|\u(\xi+\mathfrak{s},\mathfrak{s}-l,\vartheta_{-\mathfrak{s}}\omega,\tilde{\u}_l)\|^{r+1}_{\widetilde{\L}^{r+1}}\d\xi\nonumber\\&\quad+2\int_{-l}^{0}e^{\int_{0}^{-l}\left(\alpha-2\kappa e^{\sigma(\zeta+\mathfrak{s})}\mathcal{Z}_{\delta}(\vartheta_{\zeta}\omega)\right)\d\zeta}\left\langle\f(\cdot,\xi+\mathfrak{s}),\u(\xi+\mathfrak{s},\mathfrak{s}-l,\vartheta_{-\mathfrak{s}}\omega,\tilde{\u}_{l})\right\rangle\d\xi\nonumber\\&\quad+2\int_{-l}^{0}e^{\sigma(\xi+\mathfrak{s})}e^{\int_{0}^{-l}\left(\alpha-2\kappa e^{\sigma(\zeta+\mathfrak{s})}\mathcal{Z}_{\delta}(\vartheta_{\zeta}\omega)\right)\d\zeta}\mathcal{Z}_{\delta}(\vartheta_{\xi}\omega)\left(\h,\u(\xi+\mathfrak{s},\mathfrak{s}-l,\vartheta_{-\mathfrak{s}}\omega,\tilde{\u}_l)\right)\d\xi.
		\end{align}
		Making use of \eqref{ac12} in \eqref{ac22}, we get
		\begin{align}\label{ac23}
			&\limsup_{n\to\infty}\|\u(\mathfrak{s},\mathfrak{s}-t_n,\vartheta_{-\mathfrak{s}}\omega,\u_{0,n})\|^2_{\H}\nonumber\\&\leq\frac{4}{\min\{\mu,\alpha\}}\int_{-\infty}^{-l} e^{\int_{0}^{\xi}\left(\alpha-2\kappa e^{\sigma (\zeta+\mathfrak{s})}\mathcal{Z}_{\delta}(\vartheta_{\zeta}\omega)\right)\d\zeta} \|\f(\cdot,\xi+\mathfrak{s})\|^2_{\V'}\d \xi\nonumber\\&\quad+s_6\int_{-\infty}^{-l} e^{2\sigma (\xi+\mathfrak{s})} e^{\int_{0}^{\xi}\left(\alpha-2\kappa e^{\sigma (\zeta+\mathfrak{s})}\mathcal{Z}_{\delta}(\vartheta_{\zeta}\omega)\right)\d\zeta}\left|\mathcal{Z}_{\delta}(\vartheta_{\xi}\omega)\right|^2\d\xi+\|\tilde{\u}\|^2_{\H}.
		\end{align}
		Finally, passing the limit $l\to \infty$ in \eqref{ac23}, we arrive at \eqref{ac4}, which completes the proof.
	\end{proof}

	\subsubsection{Existence of random $\mathfrak{D}$-pullback attractors}
	In this subsection, we start with the proof of existence of random $\mathfrak{D}$-pullback absorbing set in $\H$ for the system \eqref{WZ_SCBF}. Finally, we prove the main result of this section, that is, the existence of random $\mathfrak{D}$-pullback attractors for the system \eqref{WZ_SCBF}.
	\begin{lemma}\label{PAS'}
		For $d=2$ with $r\geq1$, $d=3$ with $r>3$ and $d=r=3$ with $2\beta\mu\geq1$, assume that $\f\in\mathrm{L}^2_{\emph{loc}}(\R;\V')$ satisfies \eqref{forcing2} and Assumption \ref{NDT2} is fulfilled. Then there exists a closed measurable $\mathfrak{D}$-pullback absorbing set $\widetilde{\mathcal{K}}=\{\widetilde{\mathcal{K}}(\mathfrak{s},\omega):\mathfrak{s}\in\R, \omega\in\Omega\}\in\mathfrak{D}$ for the  continuous cocycle $\Phi$ associated with the system \eqref{WZ_SCBF}.
	\end{lemma}
	\begin{proof}
		When $s_3=s_4=0$, assumptions \eqref{S1} and \eqref{S4} are the same. Hence, the proof is same as the proof of Lemma \ref{PAS} by putting $s_3=s_4=0$.
	\end{proof}
	Now, we are able to provide the main results of this section. 
	\begin{theorem}\label{WZ_RA_UB}
		For $d=2$ with $r\geq1$, $d=3$ with $r>3$ and $d=r=3$ with $2\beta\mu>1$, assume that $\f\in\mathrm{L}^2_{\emph{loc}}(\R;\V')$ satisfies \eqref{forcing2} and Assumption \ref{NDT2} is fulfilled. Then there exists a unique random $\mathfrak{D}$-pullback attractor $$\widetilde{\mathscr{A}}=\{\widetilde{\mathscr{A}}(\mathfrak{s},\omega):\mathfrak{s}\in\R, \omega\in\Omega\}\in\mathfrak{D},$$ for the  the continuous cocycle $\Phi$ associated with the system \eqref{WZ_SCBF} in $\H$.
	\end{theorem}
	\begin{proof}
		The proof follows from Lemma \ref{Asymptotic_UB}, Lemma \ref{PAS'} and the abstract theory given in \cite{SandN_Wang} (Theorem 2.23 in \cite{SandN_Wang}).
	\end{proof}
	\begin{remark}\label{remark1}
		It is remarkable to note that if we replace \eqref{S3} by  
		\begin{align*}
			|\left(\mathcal{S}(\u)-\mathcal{S}(\v),\w\right)|&\leq s_2\|\u-\v\|_{\H}\|\w\|_{\H}, \ \text{ for all }\  \u,\v\in\V \text{ and }\w \in\H,
		\end{align*}
		then we can also include the case $2\beta\mu=1$ for $d=r=3$ in Lemma \ref{Continuity} and hence our main result of section \ref{sec3} and this subsection, that is, Theorems \ref{WZ_RA_B} and \ref{WZ_RA_UB} hold true for the case $d=r=3$ with $2\beta\mu=1$ also.
	\end{remark}
	\subsection{Random pullback attractors under Assumption \ref{NDT3} }\label{subsec4.3}
	In this subsection, we prove the existence of unique random $\mathfrak{D}$-pullback attractor under Assumption \eqref{NDT3} on nonlinear diffusion term $S(t,x,\u)$. In order to prove the results of this subsection, we need the following assumption on non-autonomous forcing term $\f(\cdot,\cdot)$.
	\begin{assumption}\label{DNFT3}
		We assume that external forcing term $\f\in\mathrm{L}^2_{\mathrm{loc}}(\R;\H)$  satisfies
		\begin{itemize}
			\item [(i)] 
			\begin{align}\label{forcing3}
				\int_{-\infty}^{\mathfrak{s}} e^{\alpha\xi}\|\f(\cdot,\xi)\|^2_{\H}\d \xi<\infty, \ \ \text{ for all }\  \mathfrak{s}\in\R.
			\end{align}
			Moreover, \eqref{forcing3} implies that 
			\begin{align}\label{forcing4}
				\lim_{k\to\infty}\int_{-\infty}^{0}\int\limits_{|x|\geq k} e^{\alpha\xi}|\f(x,\xi+\mathfrak{s})|^2\d x\d \xi=0, \ \ \ \text{ for all }\ \mathfrak{s}\in\R.
			\end{align}
			\item [(ii)] for every $c>0$
			\begin{align}\label{forcing5}
				\lim_{\tau\to-\infty}e^{c\tau}\int_{-\infty}^{0} e^{\alpha\xi}\|\f(\cdot,\xi+\tau)\|^2_{\H}\d \xi=0,
			\end{align}
			where $\alpha>0$ is the Darcy coefficient.
		\end{itemize}
	\end{assumption}
	\begin{lemma}\label{ContinuityUB2}
		For $d=2$ with $r\geq1$, $d=3$ with $r>3$ and $d=r=3$ with $2\beta\mu\geq1$, assume that $\f\in \mathrm{L}^2_{\emph{loc}}(\mathbb{R};\H)$ and Assumption \ref{NDT3} is fulfilled. Then, the solution of \eqref{WZ_SCBF} is continuous in initial data $\u_{\mathfrak{s}}(x).$
	\end{lemma}
	\begin{proof}
		Let $\u_{1}(t)$ and $\u_{2}(t)$ be two solutions of \eqref{WZ_SCBF}, then $\mathfrak{X}(t)=\u_{1}(t)-\u_{2}(t)$ with $\mathfrak{X}(\mathfrak{s})=\u_{1,\mathfrak{s}}(x)-\u_{2,\mathfrak{s}}(x)$ satisfies
		\begin{align}\label{Conti9}
			\frac{\d\mathfrak{X}(t)}{\d t}&=-\mu \A\mathfrak{X}(t)-\alpha\mathfrak{X}(t)-\left\{\B\big(\u_{1}(t)\big)-\B\big(\u_{2}(t)\big)\right\} -\beta\left\{\mathcal{C}\big(\u_1(t)\big)-\mathcal{C}\big(\u_2(t)\big)\right\}\nonumber\\&\quad+\left[S(t,x,\u_1(t))-S(t,x,\u_2(t))\right]\mathcal{Z}_{\delta}(\vartheta_t\omega),
		\end{align}
		in $\V'+\widetilde{\L}^{\frac{r+1}{r}}$. Taking the inner product with $\mathfrak{X}(\cdot)$ to the equation \eqref{Conti1}, we obtain
		\begin{align}\label{Conti10}
			\frac{1}{2}\frac{\d}{\d t} \|\mathfrak{X}(t)\|^2_{\H} &=-\mu \|\nabla\mathfrak{X}(t)\|^2_{\H} - \alpha\|\mathfrak{X}(t)\|^2_{\H} -\left\langle\B\big(\u_1(t)\big)-\B\big(\u_2(t)\big), \mathfrak{X}(t)\right\rangle \nonumber\\&\quad-\beta\left\langle\mathcal{C}\big(\u_1(t)\big)-\mathcal{C}\big(\u_2(t)\big),\mathfrak{X}(t)\right\rangle\nonumber\\&\quad + \mathcal{Z}_{\delta}(\vartheta_{t}\omega)\left\langle S(t,x,\u_1(t))-S(t,x,\u_2(t)),\mathfrak{X}(t)\right\rangle ,
		\end{align}
		for a.e. $t\in[\mathfrak{s},\mathfrak{s}+T] \text{ with } T>0$. By the locally Lipschitz continuity of the nonlinear diffusion term (see Assumption \ref{NDT3}), we get
		\begin{align}\label{Conti11}
			&	\mathcal{Z}_{\delta}(\vartheta_{t}\omega)\left\langle S(t,x,\u_1)-S(t,x,\u_2),\mathfrak{X}\right\rangle\nonumber\\&\leq	\left|\mathcal{Z}_{\delta}(\vartheta_{t}\omega)\right|\|S(t,x,\u_1)-S(t,x,\u_2)\|_{\H}\|\mathfrak{X}\|_{\H}\leq C	\left|\mathcal{Z}_{\delta}(\vartheta_{t}\omega)\right|\|\mathfrak{X}\|^2_{\H}.
		\end{align}
		From \eqref{MO_c}, we have
		\begin{align}\label{Conti12}
			-\beta \left\langle\mathcal{C}\big(\u_1\big)-\mathcal{C}\big(\u_2\big),\mathfrak{X}\right\rangle\leq -\frac{\beta}{2}\||\mathfrak{X}||\u_1|^{\frac{r-1}{2}}\|^2_{\H} - \frac{\beta}{2}\||\mathfrak{X}||\u_2|^{\frac{r-1}{2}}\|^2_{\H}.
		\end{align}
		\vskip 2mm
		\noindent
		\textbf{Case I:} \textit{When $d=2$ and $r>1$.} Using \eqref{b1}, \eqref{441} and Lemma \ref{Young}, we obtain
		\begin{align}\label{Conti13}
			\left| \left\langle\B\big(\u_1\big)-\B\big(\u_2\big), \mathfrak{X}\right\rangle\right|&=\left|\left\langle\B\big(\mathfrak{X},\mathfrak{X} \big), \u_1\right\rangle\right|\leq\frac{\mu}{4}\|\nabla\mathfrak{X}\|^2_{\H}+C\|\u_1\|^4_{\widetilde{\L}^4}\|\mathfrak{X}\|^2_{\H}\nonumber\\&\leq\frac{\mu}{4}\|\nabla\mathfrak{X}\|^2_{\H}+C\|\u_1\|^2_{\H}\|\nabla\u_1\|^2_{\H}\|\mathfrak{X}\|^2_{\H}.
		\end{align}
		Making use of \eqref{Conti11}-\eqref{Conti13} in \eqref{Conti10}, we get
		\begin{align}\label{Conti14}
			&	\frac{\d}{\d t} \|\mathfrak{X}(t)\|^2_{\H} \leq C \|\u_1(t)\|^2_{\H}\|\nabla\u_1(t)\|^2_{\H}\|\mathfrak{X}(t)\|^2_{\H},\text{ for a.e. } t\in[\mathfrak{s},\mathfrak{s}+T]. 
		\end{align}
		\vskip 2mm
		\noindent
		\textbf{Case II:} \textit{When $d= 3$ and $r\geq3$ ($r>3$ with any $\beta,\mu>0$ and $r=3$ with $2\beta\mu\geq1$).} The nonlinear term $\left|\left\langle\B\big(\u_1\big)-\B\big(\u_2\big), \mathfrak{X}\right\rangle\right|$ can be estimated using \eqref{441}, Lemmas \ref{Holder} and \ref{Young} as 
		\begin{align}\label{Conti16}
			\left|\left\langle\B\big(\u_1\big)-\B\big(\u_2\big), \mathfrak{X}\right\rangle\right|\leq\begin{cases}
				\frac{1}{2\beta}\|\nabla\mathfrak{X}\|^2_{\H}+\frac{\beta}{2}\||\u_1||\mathfrak{X}|\|^2_{\H},  &\text{ for } r=3,\\ \frac{\mu}{4}\|\nabla\mathfrak{X}\|_{\H}^2+\frac{\beta}{2}\||\mathfrak{X}||\u_1|^{\frac{r-1}{2}}\|^2_{\H}+C\|\mathfrak{X}\|^2_{\H} \ \ \ &\text{ for } r>3.
			\end{cases} 
		\end{align}
		Combining \eqref{Conti10}-\eqref{Conti12} and \eqref{Conti16}, we get for $r\geq3$ ($r>3$ with any $\beta,\mu>0$ and $r=3$ with $2\beta\mu\geq1$),
		\begin{align}\label{Conti17}
			&	\frac{\d}{\d t} \|\mathfrak{X}(t)\|^2_{\H}  \leq C\|\mathfrak{X}(t)\|^2_{\H},\ \ \ \text{ for a.e. }\ \ \ t\in[\mathfrak{s},\mathfrak{s}+T].
		\end{align}  Hence, we conclude the proof by applying Gronwall's inequality to \eqref{Conti14} and \eqref{Conti17}. 
	\end{proof}
	Next, we prove the existence of random $\mathfrak{D}$-pullback absorbing set for continuous cocycle $\Phi$.

	\begin{lemma}\label{LemmaUe3}
		For $d=2$ with $r\geq1$, $d=3$ with $r>3$ and $d=r=3$ with $2\beta\mu\geq1$, assume that $\f\in \mathrm{L}^2_{\emph{loc}}(\mathbb{R};\H)$ satisfies \eqref{forcing3} and Assumption \ref{NDT3} is fulfilled. Then for every $0<\delta\leq1$, $\mathfrak{s}\in\R,$ $ \omega\in \Omega$ and $D=\{D(\mathfrak{s},\omega):\mathfrak{s}\in\R, \omega\in\Omega\}\in\mathfrak{D},$ there exists $\mathscr{T}=\mathscr{T}(\delta, \mathfrak{s}, \omega, D)>0$ such that for all $t\geq \mathcal{T}$ and $\tau\geq \mathfrak{s}-t$, the solution $\u$ of the system \eqref{WZ_SCBF} with $\omega$ replaced by $\vartheta_{-\mathfrak{s}}\omega$ satisfies 
		\begin{align}\label{ue^3}
			&	\|\u(\tau,\mathfrak{s}-t,\vartheta_{-\mathfrak{s}}\omega,\u_{\mathfrak{s}-t})\|^2_{\H}+\int_{\mathfrak{s}-t}^{\tau}e^{\alpha(\xi-\tau)}\bigg[\|\u(\xi,\mathfrak{s}-t,\vartheta_{-\mathfrak{s}}\omega,\u_{\mathfrak{s}-t})\|^2_{\H} \nonumber\\&+\|\nabla\u(\xi,\mathfrak{s}-t,\vartheta_{-\mathfrak{s}}\omega,\u_{\mathfrak{s}-t})\|^2_{\H}+\|\u(\xi,\mathfrak{s}-t,\vartheta_{-\mathfrak{s}}\omega,\u_{\mathfrak{s}-t})\|^{r+1}_{\wi\L^{r+1}}\bigg]\d\xi \nonumber\\&\leq \widetilde{M} \int_{-\infty}^{\tau-\mathfrak{s}}e^{\alpha(\xi+\mathfrak{s}-\tau)}\left[\|\f(\cdot,\xi+\mathfrak{s})\|^2_{\H}+\left|\mathcal{Z}_{\delta}(\vartheta_{\xi}\omega)\right|^{\frac{r+1}{r+1-q}}+\left|\mathcal{Z}_{\delta}(\vartheta_{\xi}\omega)\right|^{\frac{r+1}{r}}\right]\d\xi,
		\end{align}
		where $\u_{\mathfrak{s}-t}\in D(\mathfrak{s}-t,\vartheta_{-t}\omega)$ and $\widetilde{M}$ is positive constant independent of $\tau, \mathfrak{s}, \omega$ and $D$.
	\end{lemma}
	\begin{proof}
		From the first equation of the system \eqref{WZ_SCBF}, we obtain
		\begin{align}\label{ue12}
			&	\frac{1}{2}\frac{\d}{\d t} \|\u\|^2_{\H} +\mu\|\nabla\u\|^2_{\H} + \alpha\|\u\|^2_{\H} + \beta\|\u\|^{r+1}_{\wi \L^{r+1}}\nonumber\\&= (\f,\u) +\mathcal{Z}_{\delta}(\vartheta_{t}\omega)\int_{\mathbb{R}^d}S(t,x,\u)\u\d x\nonumber\\&\leq \|\f\|_{\H}\|\u\|_{\H}+\left|\mathcal{Z}_{\delta}(\vartheta_{t}\omega)\right|\int_{\mathbb{R}^d}\left(\mathcal{S}_1(t,x)|\u|^q+\mathcal{S}_2(t,x)|\u|\right)\d x\nonumber\\&\leq\frac{\alpha}{4}\|\u\|^2_{\H}+\frac{\beta}{2}\|\u\|^{r+1}_{\widetilde{\L}^{r+1}}+M_1\|\f\|^2_{\H}+M_2\left|\mathcal{Z}_{\delta}(\vartheta_{t}\omega)\right|^{\frac{r+1}{r+1-q}}\|\mathcal{S}_1(t)\|^{\frac{r+1}{r+1-q}}_{\L^{\frac{r+1}{r+1-q}}(\R^d)}\nonumber\\&\quad+M_3\left|\mathcal{Z}_{\delta}(\vartheta_{t}\omega)\right|^{\frac{r+1}{r}}\|\mathcal{S}_2(t)\|^{\frac{r+1}{r}}_{\L^{\frac{r+1}{r}}(\R^d)},
		\end{align}
		where we have used \eqref{b0}, \eqref{GS1}, Lemmas \ref{Holder} and \ref{Young}, and $M_1, M_2$ and $M_3$ are positive constants independent of $\tau, \mathfrak{s}, \omega$ and $D$. From \eqref{ue12}, we find 
		\begin{align}\label{ue13}
			&\frac{\d}{\d t} \|\u\|^2_{\H}+ \alpha\|\u\|^2_{\H}+\min\left\{\frac{\alpha}{2},2\mu,\beta\right\}\left[\|\u\|^2_{\H}+\|\nabla\u\|^2_{\H} + \|\u\|^{r+1}_{\wi \L^{r+1}}\right]\nonumber\\& \leq M_1\|\f\|^2_{\H}+M_2\left|\mathcal{Z}_{\delta}(\vartheta_{t}\omega)\right|^{\frac{r+1}{r+1-q}}\|\mathcal{S}_1(t)\|^{\frac{r+1}{r+1-q}}_{\L^{\frac{r+1}{r+1-q}}(\R^d)}+M_3\left|\mathcal{Z}_{\delta}(\vartheta_{t}\omega)\right|^{\frac{r+1}{r}}\|\mathcal{S}_2(t)\|^{\frac{r+1}{r}}_{\L^{\frac{r+1}{r}}(\R^d)}.
		\end{align}
		Applying variation of constant formula to \eqref{ue13} and replacing $\omega$ by $\vartheta_{-\mathfrak{s}}\omega$, we have
		\begin{align}\label{ue14}
			&	\|\u(\tau,\mathfrak{s}-t,\vartheta_{-\mathfrak{s}}\omega,\u_{\mathfrak{s}-t})\|^2_{\H}+\min\left\{\frac{\alpha}{2},2\mu,\beta\right\}\int_{\mathfrak{s}-t}^{\tau}e^{\alpha(\xi-\tau)}\bigg[\|\u(\xi,\mathfrak{s}-t,\vartheta_{-\mathfrak{s}}\omega,\u_{\mathfrak{s}-t})\|^2_{\H} \nonumber\\&+\|\nabla\u(\xi,\mathfrak{s}-t,\vartheta_{-\mathfrak{s}}\omega,\u_{\mathfrak{s}-t})\|^2_{\H}+\|\u(\xi,\mathfrak{s}-t,\vartheta_{-\mathfrak{s}}\omega,\u_{\mathfrak{s}-t})\|^{r+1}_{\wi\L^{r+1}}\bigg]\d\xi \nonumber\\&\leq e^{\alpha(\mathfrak{s}-t-\tau)}\|\u_{\mathfrak{s}-t}\|^2_{\H}+M_4 \int_{\mathfrak{s}-t}^{\tau}e^{\alpha(\xi-\tau)}\left[\|\f(\cdot,\xi)\|^2_{\H}+\left|\mathcal{Z}_{\delta}(\vartheta_{\xi-\mathfrak{s}}\omega)\right|^{\frac{r+1}{r+1-q}}+\left|\mathcal{Z}_{\delta}(\vartheta_{\xi-\mathfrak{s}}\omega)\right|^{\frac{r+1}{r}}\right]\d\xi\nonumber\\&\leq e^{\alpha(\mathfrak{s}-t-\tau)}\|\u_{\mathfrak{s}-t}\|^2_{\H}+M_4 \int_{-\infty}^{\tau-\mathfrak{s}}e^{\alpha(\xi+\mathfrak{s}-\tau)}\left[\|\f(\cdot,\xi+\mathfrak{s})\|^2_{\H}+\left|\mathcal{Z}_{\delta}(\vartheta_{\xi}\omega)\right|^{\frac{r+1}{r+1-q}}+\left|\mathcal{Z}_{\delta}(\vartheta_{\xi}\omega)\right|^{\frac{r+1}{r}}\right]\d\xi,
		\end{align}
		where $M_4$ is a positive constant independent of $\tau, \mathfrak{s}, \omega$ and $D$. Second term of the right hand side of \eqref{ue14} is finite due to \eqref{N3} and \eqref{forcing3}. Since $\u_{\mathfrak{s}-t}\in D(\mathfrak{s}-t,\vartheta_{-t}\omega)$ and $D\in\mathfrak{D}$, we have
		\begin{align}\label{ue15}
			e^{\alpha(\mathfrak{s}-t-\tau)}\|\u_{\mathfrak{s}-t}\|^2_{\H}\leq e^{\alpha(\mathfrak{s}-t-\tau)}\|D(\mathfrak{s}-t,\vartheta_{-t}\omega)\|^2_{\H}\to0,
		\end{align}
		as $t\to\infty$, there exists $\mathscr{T}=\mathscr{T}(\delta,\mathfrak{s},\omega,D)>0$ such that for all $t\geq\mathscr{T}$,
		\begin{align}\label{ue16}
			e^{\alpha(\mathfrak{s}-t-\tau)}\|\u_{\mathfrak{s}-t}\|^2_{\H}\leq \int_{-\infty}^{\tau-\mathfrak{s}}e^{\alpha(\xi+\mathfrak{s}-\tau)}\left[\left|\mathcal{Z}_{\delta}(\vartheta_{\xi}\omega)\right|^{\frac{r+1}{r+1-q}}+\left|\mathcal{Z}_{\delta}(\vartheta_{\xi}\omega)\right|^{\frac{r+1}{r}}\right]\d\xi.
		\end{align}
		From \eqref{ue14} along with \eqref{ue16}, one can complete the proof.
	\end{proof}
	
	\begin{lemma}\label{PAS_GA}
		For $d=2$ with $r\geq1$, $d=3$ with $r>3$ and $d=r=3$ with $2\beta\mu\geq1$, assume that $\f\in\mathrm{L}^2_{\emph{loc}}(\R;\H)$ satisfies \eqref{forcing5} and Assumption \ref{NDT3} is fulfilled. Then the continuous cocycle $\Phi$ associated with the system \eqref{WZ_SCBF} possesses a closed measurable $\mathfrak{D}$-pullback absorbing set $\widehat{\mathcal{K}}=\{\widehat{\mathcal{K}}(\mathfrak{s},\omega):\mathfrak{s}\in\R, \omega\in\Omega\}\in\mathfrak{D}$ defined  for each $\mathfrak{s}\in\R$ and $\omega\in\Omega$
		\begin{align}
			\widehat{\mathcal{K}}(\mathfrak{s},\omega)=\{\u\in\H:\|\u\|^2_{\H}\leq\widehat{\mathcal{L}}(\mathfrak{s},\omega)\},
		\end{align}
		where
		\begin{align}
			\widehat{\mathcal{L}}(\mathfrak{s},\omega)=\wi M \int_{-\infty}^{0}e^{\alpha\xi}\left[\|\f(\cdot,\xi+\mathfrak{s})\|^2_{\H}+\left|\mathcal{Z}_{\delta}(\vartheta_{\xi}\omega)\right|^{\frac{r+1}{r+1-q}}+\left|\mathcal{Z}_{\delta}(\vartheta_{\xi}\omega)\right|^{\frac{r+1}{r}}\right]\d\xi,
		\end{align}
		with $\wi M$ appearing in \eqref{ue^3}.
	\end{lemma}
	\begin{proof}
		Since $\widehat{\mathcal{L}}(\mathfrak{s},\cdot):\Omega\to\R$ is $(\mathscr{F},\mathscr{B}(\R))$-measurable for every $\mathfrak{s}\in\R$, $\widehat{\mathcal{K}}(\mathfrak{s},\cdot):\Omega\to2^{\H}$ is a measurable set-valued mapping. Moreover, we have from Lemma \ref{LemmaUe3} that for each $\mathfrak{s}\in\R$, $\omega\in\Omega$ and $D\in\mathfrak{D}$, there exists $\mathscr{T}=\mathscr{T}(\delta,\mathfrak{s},\omega,D)>0$ such that for all $t\geq\mathscr{T}$,
		\begin{align}\label{PAS6}
			\Phi(t,\mathfrak{s}-t,\vartheta_{-t}\omega,D(\mathfrak{s}-t,\vartheta_{-t}\omega))=\u(\mathfrak{s},\mathfrak{s}-t,\vartheta_{-\mathfrak{s}}\omega,D(\mathfrak{s}-t,\vartheta_{-t}\omega))\subseteq\widehat{\mathcal{K}}(\mathfrak{s},\omega).
		\end{align}
		It only remains to show that $\widehat{\mathcal{K}}\in\mathfrak{D}$, that is, for every $c>0$, $\mathfrak{s}\in\R$ and $\omega\in\Omega$ $$\lim_{t\to-\infty}e^{ct}\|\widehat{\mathcal{K}}(\mathfrak{s} +t,\vartheta_{t}\omega)\|^2_{\H}=0.$$ 
		For every $c>0$, $\mathfrak{s}\in\R$ and $\omega\in\Omega$,
		\begin{align}\label{PAS7}
			&	\lim_{t\to-\infty}e^{ct}\|\widehat{\mathcal{K}}(\mathfrak{s} +t,\vartheta_{t}\omega)\|^2_{\H}=\lim_{t\to-\infty}e^{ct}\widehat{\mathcal{L}}(\mathfrak{s} +t,\vartheta_{t}\omega)\nonumber\\&=\widetilde{M}\lim_{t\to-\infty}e^{ct} \int_{-\infty}^{0}e^{\alpha\xi}\left[\|\f(\cdot,\xi+\mathfrak{s}+t)\|^2_{\H}+\left|\mathcal{Z}_{\delta}(\vartheta_{\xi+t}\omega)\right|^{\frac{r+1}{r+1-q}}+\left|\mathcal{Z}_{\delta}(\vartheta_{\xi+t}\omega)\right|^{\frac{r+1}{r}}\right]\d\xi.
		\end{align}
		Let $c_2=\min\{c,\alpha\}.$ Consider,
		\begin{align}\label{PAS8}
			&\lim_{t\to-\infty}e^{ct} \int_{-\infty}^{0}e^{\alpha\xi}\left[\left|\mathcal{Z}_{\delta}(\vartheta_{\xi+t}\omega)\right|^{\frac{r+1}{r+1-q}}+\left|\mathcal{Z}_{\delta}(\vartheta_{\xi+t}\omega)\right|^{\frac{r+1}{r}}\right]\d\xi\nonumber\\&\leq\lim_{t\to-\infty} \int_{-\infty}^{0}e^{c_2(\xi+t)}\left[\left|\mathcal{Z}_{\delta}(\vartheta_{\xi+t}\omega)\right|^{\frac{r+1}{r+1-q}}+\left|\mathcal{Z}_{\delta}(\vartheta_{\xi+t}\omega)\right|^{\frac{r+1}{r}}\right]\d\xi\nonumber\\&\leq\lim_{t\to-\infty} \int_{-\infty}^{t}e^{c_2\xi}\left[\left|\mathcal{Z}_{\delta}(\vartheta_{\xi}\omega)\right|^{\frac{r+1}{r+1-q}}+\left|\mathcal{Z}_{\delta}(\vartheta_{\xi}\omega)\right|^{\frac{r+1}{r}}\right]\d\xi=0,
		\end{align}
		where we have used the fact that $\int_{-\infty}^{0}e^{c_2\xi}\left[\left|\mathcal{Z}_{\delta}(\vartheta_{\xi}\omega)\right|^{\frac{r+1}{r+1-q}}+\left|\mathcal{Z}_{\delta}(\vartheta_{\xi}\omega)\right|^{\frac{r+1}{r}}\right]\d\xi<\infty$ due to \eqref{N3}.
		Hence \eqref{forcing5}, \eqref{PAS7} and \eqref{PAS8} imply  
		\begin{align*}
			\lim_{t\to-\infty}e^{ct}\|\widehat{\mathcal{K}}(\mathfrak{s} +t,\vartheta_{t}\omega)\|^2_{\H}=0,
		\end{align*}
		as required.
	\end{proof}
	Next lemma plays a pivotal role to prove the $\mathfrak{D}$-pullback asymptotic compactness of $\Phi$ (see Lemma \ref{Asymptotic_UB_GS}).
	\begin{lemma}\label{largeradius}
		For $d=2$ with $r\geq1, d=3$ with $r>3$ and $d=r=3$ with $2\beta\mu\geq1$, assume that $\f\in\mathrm{L}^2_{\mathrm{loc}}(\R;\H)$ and satisfies \eqref{forcing3}. Then, for any $\u_{\mathfrak{s}-t}\in D(\mathfrak{s}-t,\vartheta_{-t}\omega),$ where $D=\{D(\mathfrak{s},\omega):\mathfrak{s}\in\R, \omega\in \Omega\}\in\mathfrak{D}$, and for any $\eta>0$, $\mathfrak{s}\in\R$, $\omega\in \Omega$ and $0<\delta\leq1$, there exists $\mathscr{T}^*=\mathscr{T}^*(\delta,\mathfrak{s},\omega,D,\eta)\geq 1$ and $P^*=P^*(\delta,\mathfrak{s},\omega,\eta)>0$ such that for all $t\geq \mathscr{T}^*$ and $\tau\in[\mathfrak{s}-1,\mathfrak{s}]$, the solution of \eqref{WZ_SCBF}  with $\omega$ replaced by $\vartheta_{-\mathfrak{s}}\omega$ satisfy
		\begin{align}\label{ep}
			\int_{|x|\geq P^*}|\u(\tau,\mathfrak{s}-t,\vartheta_{-\mathfrak{s}}\omega,\u_{\mathfrak{s}-t}) |^2\d x\leq \eta.
		\end{align}
	\end{lemma}
	\begin{proof}
		Let $\Psi$ be a smooth function such that $0\leq\Psi(s)\leq 1$ for $s\in\R^+$ and 
		\begin{align}\label{Psi}
			\Psi(s)=\begin{cases*}
				0,\quad \text{ for }0\leq s\leq 1,\\
				1, \quad\text{ for } s\geq2 .
			\end{cases*}
		\end{align}
		Then, there exists a positive constant $C$ such that $|\Psi'(s)|\leq C$ for all $s\in\R^+$. Taking the inner product of first equation of \eqref{WZ_SCBF} with $\Psi\left(\frac{|x|^2}{k^2}\right)\u$ in $\H$, we have
		\begin{align}\label{ep1}
			\frac{1}{2} \frac{\d}{\d t} \int_{\R^d}\Psi\left(\frac{|x|^2}{k^2}\right)|\u|^2\d x &= -\mu \int_{\R^d}(\A\u) \Psi\left(\frac{|x|^2}{k^2}\right) \u \d x-\alpha \int_{\R^d}\Psi\left(\frac{|x|^2}{k^2}\right)|\u|^2\d x\nonumber\\&\quad-b\left(\u,\u,\Psi\left(\frac{|x|^2}{k^2}\right)\u\right)-\beta \int_{\R^d}\Psi\left(\frac{|x|^2}{k^2}\right)|\u|^{r+1}\d x\nonumber\\&\quad+ \int_{\R^d}\f(x,t)\Psi\left(\frac{|x|^2}{k^2}\right)\u\d x +\mathcal{Z}_{\delta}(\vartheta_{t}\omega)\int_{\mathbb{R}^d}\Psi\left(\frac{|x|^2}{k^2}\right)S(t,x,\u)\u\d x.
		\end{align}
		We estimate each term on right hand side of \eqref{ep1}. Integration by parts helps to obtain
		\begin{align}\label{ep2}
			-&\mu \int_{\R^d}(\A\u) \Psi\left(\frac{|x|^2}{k^2}\right) \u \d x\nonumber\\&= -\mu \int_{\R^d}|\nabla\u|^2 \Psi\left(\frac{|x|^2}{k^2}\right)  \d x -\mu \int_{\R^d} \Psi'\left(\frac{|x|^2}{k^2}\right)\frac{2}{k^2}(x\cdot\nabla) \u\cdot\u \d x\nonumber\\&= -\mu \int_{\R^d}|\nabla\u|^2 \Psi\left(\frac{|x|^2}{k^2}\right)  \d x -\mu \int\limits_{k\leq|x|\leq \sqrt{2}k}\Psi'\left(\frac{|x|^2}{k^2}\right)\frac{2}{k^2}(x\cdot\nabla) \u\cdot\u \d x\nonumber\\&\leq -\mu \int_{\R^d}|\nabla\u|^2 \Psi\left(\frac{|x|^2}{k^2}\right)  \d x +\frac{2\sqrt{2}\mu}{k} \int\limits_{k\leq|x|\leq \sqrt{2}k}\left|\u\right| \left|\Psi'\left(\frac{|x|^2}{k^2}\right)\right|\left|\nabla \u\right| \d x\nonumber\\&\leq -\mu \int_{\R^d}|\nabla\u|^2 \Psi\left(\frac{|x|^2}{k^2}\right)  \d x +\frac{C}{k} \int_{\R^d}\left|\u\right| \left|\nabla \u\right| \d x\nonumber\\&\leq -\mu \int_{\R^d}|\nabla\u|^2 \Psi\left(\frac{|x|^2}{k^2}\right)  \d x +\frac{C}{k} \left(\|\u\|^2_{\H}+\|\nabla\u\|^2_{\H}\right),
		\end{align}
		and
		\begin{align}\label{ep3}
			-b\left(\u,\u,\Psi\left(\frac{|x|^2}{k^2}\right)\u\right)&=\int_{\R^d} \Psi'\left(\frac{|x|^2}{k^2}\right)\frac{x}{k^2}\cdot\u |\u|^2 \d x=  \int\limits_{k\leq|x|\leq \sqrt{2}k} \Psi'\left(\frac{|x|^2}{k^2}\right)\frac{x}{k^2}\cdot\u |\u|^2 \d x\nonumber\\&\leq \frac{\sqrt{2}}{k} \int\limits_{k\leq|x|\leq \sqrt{2}k} \left|\Psi'\left(\frac{|x|^2}{k^2}\right)\right| |\u|^3 \d x\leq \frac{C}{k}\|\u\|^3_{\wi \L^3}.
		\end{align}
		Using Lemmas \ref{Interpolation} and \ref{Young} in \eqref{ep3}, we get
		\begin{align}
			-b\left(\u,\u,\Psi\left(\frac{|x|^2}{k^2}\right)\u\right)\leq\begin{cases}
				\frac{C}{k}\left[\|\u\|^2_{\H}+\|\u\|^{4}_{\wi \L^{4}}\right],&\ \ \ \text{ for } d=2 \text{ with } r\geq1,\\ \vspace{1mm}
				\frac{C}{k}\left[\|\u\|^2_{\H}+\|\u\|^{r+1}_{\wi \L^{r+1}}\right], &\ \ \ \text{ for } d=3 \text{ with } r\geq3.
			\end{cases}
		\end{align}
		Finally, we estimate the last two terms of \eqref{ep1} as follows
		\begin{align}
			\int_{\R^d}\f(x,t)\Psi\left(\frac{|x|^2}{k^2}\right)\u \d x\leq \frac{\alpha}{2} \int_{\R^d}\Psi\left(\frac{|x|^2}{k^2}\right)|\u|^2\d x +\frac{1}{2\alpha} \int_{\R^d}\Psi\left(\frac{|x|^2}{k^2}\right)|\f(x,t)|^2\d x,
		\end{align}  
		and
		\begin{align}\label{ep4}
			&\mathcal{Z}_{\delta}(\vartheta_{t}\omega)\int_{\mathbb{R}^d}\Psi\left(\frac{|x|^2}{k^2}\right)S(t,x,\u)\u\d x\nonumber\\&\leq\left|\mathcal{Z}_{\delta}(\vartheta_{t}\omega)\right|\int_{\R^d}\Psi\left(\frac{|x|^2}{k^2}\right)\left(\mathcal{S}_1(t,x)|\u|^q+\mathcal{S}_2(t,x)|\u|\right)\d x\nonumber\\&\leq\frac{\beta}{2} \int_{\R^d}\Psi\left(\frac{|x|^2}{k^2}\right)|\u|^{r+1}\d x+C\left|\mathcal{Z}_{\delta}(\vartheta_{t}\omega)\right|^{\frac{r+1}{r+1-q}}\int_{\R^d}\Psi\left(\frac{|x|^2}{k^2}\right)\left|\mathcal{S}_1(t,x)\right|^{\frac{r+1}{r+1-q}}\d x\nonumber\\&\quad+C\left|\mathcal{Z}_{\delta}(\vartheta_{t}\omega)\right|^{\frac{r+1}{r}}\int_{\R^d}\Psi\left(\frac{|x|^2}{k^2}\right)\left|\mathcal{S}_2(t,x)\right|^{\frac{r+1}{r}}\d x.
		\end{align}
		Combining \eqref{ep1}-\eqref{ep4}, we get
		\begin{align}\label{ep5}
			&	\frac{\d}{\d t} \int_{\R^d}\Psi\left(\frac{|x|^2}{k^2}\right)|\u|^2\d x \nonumber\\&\leq -\alpha \int_{\R^d}\Psi\left(\frac{|x|^2}{k^2}\right)|\u|^2\d x+\frac{C}{k} \left(\|\u\|^2_{\V}+\|\u\|^{j}_{\wi \L^{j}}\right)+\frac{1}{\alpha} \int_{|x|\geq k}|\f(x,t)|^2\d x\nonumber\\&\quad+C\left|\mathcal{Z}_{\delta}(\vartheta_{t}\omega)\right|^{\frac{r+1}{r+1-q}}\int_{|x|\geq k}\left|\mathcal{S}_1(t,x)\right|^{\frac{r+1}{r+1-q}}\d x+C\left|\mathcal{Z}_{\delta}(\vartheta_{t}\omega)\right|^{\frac{r+1}{r}}\int_{|x|\geq k}\left|\mathcal{S}_2(t,x)\right|^{\frac{r+1}{r}}\d x,
		\end{align}
		where $j=4$ and $j=r+1$ for $d=2$ and $d=3$, respectively. Applying Gronwall's inequality to \eqref{ep5} on $(\mathfrak{s}-t,\tau)$ where $\tau\in[\mathfrak{s}-1,\mathfrak{s}]$ and replacing $\omega$ by $\vartheta_{-\mathfrak{s}}\omega$, we obtain that, for $\mathfrak{s}\in\R, t\geq 1$ and $\omega\in \Omega$,
		\begin{align}\label{ep6}
			&\int_{\R^d}\Psi\left(\frac{|x|^2}{k^2}\right)|\u(\tau,\mathfrak{s}-t,\vartheta_{-\mathfrak{s}}\omega,\u_{\mathfrak{s}-t})|^2\d x \nonumber\\&\leq e^{\alpha( \mathfrak{s}-t-\tau)}\int_{\R^d}\Psi\left(\frac{|x|^2}{k^2}\right)|\u_{\mathfrak{s}-t}|^2\d x\nonumber\\&\quad+\frac{C}{k}\int_{\mathfrak{s}-t}^{\tau}e^{\alpha(\xi-\tau)} \left(\|\u(\xi,\mathfrak{s}-t,\vartheta_{-\mathfrak{s}}\omega,\u_{\mathfrak{s}-t})\|^2_{\V}+\|u(\xi,\mathfrak{s}-t,\vartheta_{-\mathfrak{s}}\omega,\u_{\mathfrak{s}-t})\|^{j}_{\wi\L^{j}}\right)\d\xi\nonumber\\&\quad+C\int_{-\infty}^{\tau-\mathfrak{s}}e^{\alpha (\xi+\mathfrak{s}-\tau)}\bigg[\int_{|x|\geq k}|\f(x,\xi+\mathfrak{s})|^2\d x+\left|\mathcal{Z}_{\delta}(\vartheta_{\xi}\omega)\right|^{\frac{r+1}{r+1-q}}\int_{|x|\geq k}\left|\mathcal{S}_1(\xi+\mathfrak{s},x)\right|^{\frac{r+1}{r+1-q}}\d x\nonumber\\&\qquad+\left|\mathcal{Z}_{\delta}(\vartheta_{\xi}\omega)\right|^{\frac{r+1}{r}}\int_{|x|\geq k}\left|\mathcal{S}_2(\xi+\mathfrak{s},x)\right|^{\frac{r+1}{r}}\d x\bigg]\d \xi,
		\end{align}
		where $j=4$ and $j=r+1$ for $d=2$ and $d=3$, respectively. Since $\u_{\mathfrak{s}-t}\in D(\mathfrak{s}-t,\vartheta_{-t}\omega)$ and $D\in\mathfrak{D}$, from \eqref{ue15}, we find that, for given $\eta>0,$ there exists $T_1=T_1(\delta,\mathfrak{s},\omega,D,\eta)\geq 1$ such that  for all $t\geq T_1,$
		\begin{align}\label{ep7}
			e^{\alpha (\mathfrak{s}-t-\tau)}\int_{\R^d}\Psi\left(\frac{|x|^2}{k^2}\right)|\u_{\mathfrak{s}-t}|^2\d x\leq\frac{\eta}{4}.
		\end{align}
		Since $\f\in\mathrm{L}^2_{\mathrm{loc}}(\R;\H)$ satisfies \eqref{forcing4}, for a given $\eta>0,$ there exists $P_1=P_1(\delta,\mathfrak{s},\omega,\eta)>0$ such that for all $k\geq P_1$,
		\begin{align}\label{ep8}
			&C\int_{-\infty}^{\tau-\mathfrak{s}}e^{\alpha (\xi+\mathfrak{s}-\tau)} \int_{|x|\geq k}|\f(x,\xi+\mathfrak{s})|^2\d x\d\xi\leq Ce^{\alpha}\int_{-\infty}^{0}e^{\alpha \xi} \int_{|x|\geq k}|\f(x,\xi+\mathfrak{s})|^2\d x\d\xi\leq \frac{\eta}{4}.
		\end{align}
		By \eqref{N3}, we have
		\begin{align*}
			&	C\int_{-\infty}^{\tau-\mathfrak{s}}e^{\alpha (\xi+\mathfrak{s}-\tau)}\bigg[\left|\mathcal{Z}_{\delta}(\vartheta_{\xi}\omega)\right|^{\frac{r+1}{r+1-q}}\int_{\R^d}\left|\mathcal{S}_1(\xi+\mathfrak{s},x)\right|^{\frac{r+1}{r+1-q}}\d x\nonumber\\&\qquad+\left|\mathcal{Z}_{\delta}(\vartheta_{\xi}\omega)\right|^{\frac{r+1}{r}}\int_{\R^d}\left|\mathcal{S}_2(\xi+\mathfrak{s},x)\right|^{\frac{r+1}{r}}\d x\bigg]\d \xi\nonumber\\&\leq C\int_{-\infty}^{\tau-\mathfrak{s}}e^{\alpha (\xi+\mathfrak{s}-\tau)}\bigg[\left|\mathcal{Z}_{\delta}(\vartheta_{\xi}\omega)\right|^{\frac{r+1}{r+1-q}}+\left|\mathcal{Z}_{\delta}(\vartheta_{\xi}\omega)\right|^{\frac{r+1}{r}}\bigg]\d \xi <\infty,
		\end{align*}
		from which it follows that, for given $\eta>0,$ there exists $P_2=P_2(\delta,\mathfrak{s},\omega,\eta)>0$ such that for all $k\geq P_2$,
		\begin{align}\label{ep9}
			&	C\int_{-\infty}^{\tau-\mathfrak{s}}e^{\alpha (\xi+\mathfrak{s}-\tau)}\bigg[\left|\mathcal{Z}_{\delta}(\vartheta_{\xi}\omega)\right|^{\frac{r+1}{r+1-q}}\int_{|x|\geq k}\left|\mathcal{S}_1(\xi+\mathfrak{s},x)\right|^{\frac{r+1}{r+1-q}}\d x\nonumber\\&\qquad+\left|\mathcal{Z}_{\delta}(\vartheta_{\xi}\omega)\right|^{\frac{r+1}{r}}\int_{|x|\geq k}\left|\mathcal{S}_2(\xi+\mathfrak{s},x)\right|^{\frac{r+1}{r}}\d x\bigg]\d \xi\leq\frac{\eta}{4}.
		\end{align}
		Due to \eqref{ue^3} and \eqref{ep8}-\eqref{ep9}, there exists $T_2=T_2(\delta,\mathfrak{s},\omega,D,\eta)>0$ and $P_3=P_3(\delta,\mathfrak{s},\omega,\eta)>0$ such that for all $t\geq T_2$ and $k\geq P_3$,
		\begin{align}\label{ep10}
			\frac{C}{k}\int_{\mathfrak{s}-t}^{\tau}e^{\alpha(\xi-\tau)} \left(\|\u(\xi,\mathfrak{s}-t,\vartheta_{-\mathfrak{s}}\omega,\u_{\mathfrak{s}-t})\|^2_{\V}+\|\u(\xi,\mathfrak{s}-t,\vartheta_{-\mathfrak{s}}\omega,\u_{\mathfrak{s}-t})\|^{r+1}_{\wi\L^{r+1}}\right)\d\xi\leq \frac{\eta}{4},
		\end{align}
		for $d=3$ and
		\begin{align}\label{ep10*}
			\frac{C}{k}\int_{\mathfrak{s}-t}^{\tau}e^{\alpha(\xi-\tau)} \left(\|\u(\xi,\mathfrak{s}-t,\vartheta_{-\mathfrak{s}}\omega,\u_{\mathfrak{s}-t})\|^2_{\V}+\|\u(\xi,\mathfrak{s}-t,\vartheta_{-\mathfrak{s}}\omega,\u_{\mathfrak{s}-t})\|^{4}_{\wi\L^{4}}\right)\d\xi\leq \frac{\eta}{4},
		\end{align}
		where we have used $\|\u\|^4_{\wi\L^4}\leq\|\u\|^2_{\H}\|\nabla\u\|^2_{\H}$, for $d=2$.

		Let $\mathscr{T}^*=\mathscr{T}^*(\delta,\mathfrak{s},\omega,D,\eta)=\max\{T_1,T_2\}, P^*=P^*(\delta,\mathfrak{s},\omega,\eta)=\max\{P_1,P_2,P_3\}$. Then it implies from \eqref{ep6}-\eqref{ep10} that, for all $t\geq \mathscr{T}^*$ and $k\geq P^*$, we obtain that for all $\tau\in[\mathfrak{s}-1,\mathfrak{s}]$,
		\begin{align*}
			\int\limits_{|x|\geq P^*}|\u(\tau,\mathfrak{s}-t,\vartheta_{-\mathfrak{s}}\omega,\u_{\mathfrak{s}-t})|^2\d x\leq \eta,
		\end{align*}
		which completes the proof. 
	\end{proof}
	\begin{lemma}\label{PCB_U}
		For $d=2$ with $r\geq1$, $d=3$ with $r>3$ and $d=r=3$ with $2\beta\mu\geq1$, assume that $\f\in \mathrm{L}^2_{\emph{loc}}(\mathbb{R};\H)$, $\{\u(\cdot,\mathfrak{s},\omega,\u^n)\}_{n\in\N}$ is a bounded sequence of solutions of \eqref{WZ_SCBF} in $\H$ and Assumption \ref{NDT3} is fulfilled. Then for every $0<\delta\leq1, \omega\in\Omega, \mathfrak{s}\in\R$ and $t>\mathfrak{s}$, there exists $\u^0\in\mathrm{L}^2(\mathfrak{s},\mathfrak{s}+T;\H)$ with $T>0$ and a subsequence $\{\u(\cdot,\mathfrak{s},\omega,\u^{n_m})\}_{m\in\N}$ of $\{\u(\cdot,\mathfrak{s},\omega,\u^n)\}_{n\in\N}$ such that $\u(\tau,\mathfrak{s},\omega,\u^{n_m})\to\u^0(\tau)$ in $\L^2(\mathcal{O}_k)$ as $m\to\infty$ for every $k\in\N$ and for almost all $\tau\in(\mathfrak{s},\mathfrak{s}+T)$, where $$\mathcal{O}_k=\{x\in\R^d:|x|<k\}.$$
	\end{lemma}
	\begin{proof}
		Since embedding $\H^1_0(\mathcal{O}_k)\hookrightarrow\L^2(\mathcal{O}_k)$ is compact, the proof can be completed analogously as in  the proof of Lemma \ref{PCB}.
	\end{proof}
	The following theorem proves the $\mathfrak{D}$-pullback asymptotic compactness of $\Phi$.
	\begin{lemma}\label{Asymptotic_UB_GS}
		For $d=2$ with $r\geq1$, $d=3$ with $r>3$ and $d=r=3$ with $2\beta\mu\geq1$, assume that $\f\in\mathrm{L}^2_{\emph{loc}}(\R;\H)$ satisfies \eqref{forcing4} and Assumption \ref{NDT3} is fulfilled. Then for every $0<\delta\leq1$, $\mathfrak{s}\in \R,$ $\omega\in \Omega,$ $D=\{D(\mathfrak{s},\omega):\mathfrak{s}\in \R,\omega\in \Omega\}\in \mathfrak{D}$ and $t_n\to \infty,$ $\u_{0,n}\in D(\mathfrak{s}-t_n, \vartheta_{-t_{n}}\omega)$, the sequence $\Phi(t_n,\mathfrak{s}-t_n,\vartheta_{-t_n}\omega,\u_{0,n})$ or $\u(\mathfrak{s},\mathfrak{s}-t_n,\vartheta_{-\mathfrak{s}}\omega,\u_{0,n})$ of solutions of the system \eqref{WZ_SCBF} has a convergent subsequence in $\H$.
	\end{lemma}
	\begin{proof}
		Lemma \ref{LemmaUe3} implies that there exists $\mathscr{T}=\mathscr{T}(\delta,\mathfrak{s},\omega,D)>0$ and $R(\mathfrak{s},\omega)$ such that for all $t\geq \mathscr{T}$,
		\begin{align}\label{Uac1}
			\|\u(\mathfrak{s}-1,\mathfrak{s}-t,\vartheta_{-\mathfrak{s}}\omega,\u_{\mathfrak{s}-t})\|^2_{\H} \leq R(\mathfrak{s},\omega),
		\end{align}
		where $\u_{\mathfrak{s}-t}\in D(\mathfrak{s}-t,\vartheta_{-t}\omega).$ Since $t_n\to \infty$, there exists $N_2\in\N$ such that $t_n\geq \mathscr{T}$ for all $n\geq N_2$. Since $\u_{0,n}\in D(\mathfrak{s}-t_n, \vartheta_{-t_{n}}\omega)$, \eqref{Uac1} implies that for all $n\geq N_2$, 
		\begin{align*}
			\|\u(\mathfrak{s}-1,\mathfrak{s}-t_n,\vartheta_{-\mathfrak{s}}\omega,\u_{0,n})\|^2_{\H} \leq R(\mathfrak{s},\omega),
		\end{align*}
		and hence 
		\begin{align}\label{Uac2}
			\{\u(\mathfrak{s}-1,\mathfrak{s}-t_n,\vartheta_{-\mathfrak{s}}\omega,\u_{0,n})\}_{n\geq N_2} \text{ is a bounded sequence in }\H.
		\end{align}
		It yields from \eqref{Uac2} and Lemma \ref{PCB_U} that there exists $\tau\in(\mathfrak{s}-1,\mathfrak{s})$, $\u_0$ in $\H$ and a subsequence (not relabeled) such that for every $k\in\N$ as $n\to\infty$
		\begin{align}\label{Uac3}
			\u(\tau,\mathfrak{s}-t_n,\vartheta_{-\mathfrak{s}}\omega,\u_{0,n})=\u(\tau,\mathfrak{s}-1,\vartheta_{-\mathfrak{s}}\omega,\u(\mathfrak{s}-1,\mathfrak{s}-t_n,\vartheta_{-\mathfrak{s}}\omega,\u_{0,n}))\to \u_0 \  \text{ in }\  \L^2(\mathcal{O}_k).
		\end{align}
		Since $\u_0\in\H$, for given $\eta>0$, there exists $K_1=K_1(\delta,\mathfrak{s},\omega,\eta)>0$ such that for all $k\geq K_1$, 
		\begin{align}\label{Uac4}
			\int\limits_{|x|\geq k}|\u_0|^2\d x\leq\frac{\eta}{3}.
		\end{align}
		Also, it follows from Lemma \ref{largeradius} that there exists $N_3=N_3(\delta,\mathfrak{s},\omega,D,\eta)\geq 1$ and $K_2=K_2(\delta,\mathfrak{s},\omega,\eta)\geq K_1$ such that for all $n\geq N_3$ and $k\geq K_2$,
		\begin{align}\label{Uac5}
			\int\limits_{|x|\geq k}|\u(\tau,\mathfrak{s}-t_n,\vartheta_{-\mathfrak{s}}\omega,\u_{0,n})|^2\d x\leq\frac{\eta}{3}.
		\end{align}
		From \eqref{Uac3}, we have that there exists $N_4=N_4(\delta,\mathfrak{s},\omega,D,\eta)>N_3$ such that for all $n\geq N_4$,
		\begin{align}\label{Uac6}
			\int\limits_{|x|< K_2}|\u(\tau,\mathfrak{s}-t_n,\vartheta_{-\mathfrak{s}}\omega,\u_{0,n})-\u_0|^2\d x\leq\frac{\eta}{3}.
		\end{align}
		Finally, Lemma \ref{ContinuityUB2} implies
		\begin{align}\label{Uac7}
			&\|\u(\mathfrak{s},\tau,\vartheta_{-\mathfrak{s}}\omega,\u(\tau,\mathfrak{s}-t_n,\vartheta_{-\mathfrak{s}}\omega,\u_{0,n}))-\u(\mathfrak{s},\tau,\vartheta_{-\mathfrak{s}}\omega,\u_0)\|^2_{\H}\nonumber\\&\leq C\|\u(\tau,\mathfrak{s}-t_n,\vartheta_{-\mathfrak{s}}\omega,\u_{0,n})-\u_0\|^2_{\H}\nonumber\\&\leq C\left[\int\limits_{|x|<K_2}|\u(\tau,\mathfrak{s}-t_n,\vartheta_{-\mathfrak{s}}\omega,\u_{0,n})-\u_0|^2\d x+\int\limits_{|x|\geq K_2}|\u(\tau,\mathfrak{s}-t_n,\vartheta_{-\mathfrak{s}}\omega,\u_{0,n})-\u_0|^2\d x\right]\nonumber\\&\leq C\left[\int\limits_{|x|<K_2}|\u(\tau,\mathfrak{s}-t_n,\vartheta_{-\mathfrak{s}}\omega,\u_{0,n})-\u_0|^2\d x+\hspace{-2mm}\int\limits_{|x|\geq K_2}(|\u(\tau,\mathfrak{s}-t_n,\vartheta_{-\mathfrak{s}}\omega,\u_{0,n})|^2+|\u_0|^2)\d x\right].
		\end{align}
		Hence, \eqref{Uac7} along with \eqref{Uac4}-\eqref{Uac6} conclude the proof.
	\end{proof}
	The following theorem is the main results of this section, that is, the existence of a unique random $\mathfrak{D}$-pullback attractor for $\Phi$. 
	\begin{theorem}\label{WZ_RA_UB_GS}
		For $d=2$ with $r\geq1$, $d=3$ with $r>3$ and $d=r=3$ with $2\beta\mu\geq1$, assume that $\f\in\mathrm{L}^2_{\emph{loc}}(\R;\H)$ satisfies \eqref{forcing5} and Assumption \ref{NDT2} is fulfilled. Then there exists a unique random $\mathfrak{D}$-pullback attractor $$\widehat{\mathscr{A}}=\{\widehat{\mathscr{A}}(\mathfrak{s},\omega):\mathfrak{s}\in\R, \omega\in\Omega\}\in\mathfrak{D},$$ for the  the continuous cocycle $\Phi$ associated with the system \eqref{WZ_SCBF} in $\H$.
	\end{theorem}
	\begin{proof}
		Note that Lemmas \ref{PAS_GA} and \ref{Asymptotic_UB_GS} provide the existence of a closed measurable $\mathfrak{D}$-pullback absorbing set for $\Phi$ and asymptotic compactness of $\Phi$, respectively. Hence, Lemmas \ref{PAS_GA} and \ref{Asymptotic_UB_GS} together with the abstract theory given in \cite{SandN_Wang} (Theorem 2.23, \cite{SandN_Wang}) complete the proof.
	\end{proof}

	\section{Convergence of attractors: Additive white noise} \label{sec6}\setcounter{equation}{0}
	In this section, we examine the approximations of solutions of the following stochastic CBF equations with additive white noise. For given $\sigma\geq0$ and $\textbf{g}\in\D(\A)$,
	\begin{equation}\label{SCBF_Add}
		\left\{
		\begin{aligned}
			\frac{\partial \u}{\partial t}+\mu\A\u+\B(\u)+\alpha\u+\beta\mathcal{C}(\u)&=\boldsymbol{f}+e^{\sigma t}\textbf{g}(x)\frac{\d \W}{\d t}, \ \ \ \text{ in } \mathbb{R}^n\times(\mathfrak{s},\infty), \\ 
			\u(x,\mathfrak{s})&=\u_{\mathfrak{s}}(x),\ \ \ \ \ \ \ \ \ \ \ \ \ \ \ \ \ x\in \mathbb{R}^n \text{ and }\mathfrak{s}\in\R.
		\end{aligned}
		\right.
	\end{equation}
	For $\delta>0$, consider the pathwise random equations:
	\begin{equation}\label{WZ_SCBF_Add}
		\left\{
		\begin{aligned}
			\frac{\partial \u_{\delta}}{\partial t}+\mu\A\u_{\delta}+\B(\u_{\delta})+\alpha\u_{\delta}+\beta\mathcal{C}(\u_{\delta})&=\boldsymbol{f}+e^{\sigma t}\textbf{g}(x)\mathcal{Z}_{\delta}(\vartheta_{t}\omega), \ \ \text{ in } \mathbb{R}^n\times(\mathfrak{s},\infty), \\ 
			\u_{\delta}(x,\mathfrak{s})&=\u_{\delta,\mathfrak{s}}(x),\ \ \ \  \ \ \ \ \ \ \ \ \ \ \ \ \ \ \ x\in \mathbb{R}^n \text{ and }\mathfrak{s}\in\R.
		\end{aligned}
		\right.
	\end{equation}
	Throughout this section, we  prove results for 2D SCBF with $r>1$ and 3D SCBF with $r\geq3$ ($r>3$ for any $\beta,\mu>0$ and $r=3$ for $2\beta\mu\geq1$). Since, 2D SCBF equations are linear perturbation of 2D stochastic NSE for $r=1$, one can prove the results of this section for $d=2$ and $r=1$ by using the same arguments as it is done for 2D stochastic NSE on Poincar\'e domains in \cite{GGW} (see Section 3 in \cite{GGW}). In \cite{GGW}, authors imposed an extra condition (Assumption \ref{GA}) on $\textbf{g}(\cdot),$ which was used to prove the existence of $\mathfrak{D}$-pullback absorbing set. We observe that there is no need to impose Assumption \ref{GA} on $\textbf{g}(\cdot)$ for $r>1$. 
	\begin{assumption}\label{DNFT2}
		For this section, we assume that external forcing term $\f$ satisfies the following two assumptions.
		\begin{itemize}
			\item [(i)] 
			\begin{align}\label{forcing6}
				\int_{-\infty}^{\mathfrak{s}} e^{\alpha\xi}\|\f(\cdot,\xi)\|^2_{\V'}\d \xi<\infty, \ \ \text{ for all }\  \mathfrak{s}\in\R.
			\end{align}
			\item [(ii)] for every $c>0$
			\begin{align}\label{forcing7}
				\lim_{\tau\to-\infty}e^{c\tau}\int_{-\infty}^{0} e^{\alpha\xi}\|\f(\cdot,\xi+\tau)\|^2_{\V'}\d \xi=0,
			\end{align}
			where $\alpha>0$ is Darcy coefficient.
		\end{itemize}
	\end{assumption}
	\subsection{Random $\mathfrak{D}$-pullback attractor for SCBF equations with additive white noise}
	Consider, for some $\ell>0$ 
	\begin{align*}
		\y(\vartheta_{t}\omega) =  \int_{-\infty}^{t} e^{-\ell(t-\tau)}\d \W(\tau), \ \ \omega\in \Omega,
	\end{align*} which is the stationary solution of the one dimensional Ornstein-Uhlenbeck equation
	\begin{align*}
		\d\y(\vartheta_t\omega) + \ell\y(\vartheta_t\omega)\d t =\d\W(t).
	\end{align*}
	Let us recall from \cite{FAN} that there exists a $\vartheta$-invariant subset of $\Omega$ (will be denoted by $\Omega$ itself) of full measure such that $\y(\vartheta_t\omega)$ is continuous in $t$ for every $\omega\in \Omega,$ and
	\begin{align}
		\lim_{t\to \pm \infty} \frac{|\y(\vartheta_t\omega)|}{|t|}=0   \text{\ \ and  \ \ }
		\lim_{t\to \pm \infty} \frac{1}{t} \int_{0}^{t} \y(\vartheta_{\xi}\omega)\d\xi =0.\label{Y2}
	\end{align}
	Define
	\begin{align}\label{T_add}
		\v(t,\mathfrak{s},\omega,\v_{\mathfrak{s}})=\u(t,\mathfrak{s},\omega,\u_{\mathfrak{s}})-e^{\sigma t}\textbf{g}(x)\y(\vartheta_{t}\omega).
	\end{align}
	Then, from \eqref{SCBF_Add}, we obtain
	\begin{equation}\label{CSCBF_Add}
		\left\{
		\begin{aligned}
			\frac{\partial \v}{\partial t}+\mu\A\v+\B(\v&+e^{\sigma t}\textbf{g}\y)+\alpha\v+\beta\mathcal{C}(\v+e^{\sigma t}\textbf{g}\y)\\&=\boldsymbol{f}+\left(\ell-\sigma-\alpha\right)e^{\sigma t}\textbf{g}\y-\mu e^{\sigma t}\y\A\textbf{g}, \ \ \ \text{ in } \mathbb{R}^d\times(\mathfrak{s},\infty), \\ 
			\v(x,\mathfrak{s})&=\v_{\mathfrak{s}}(x)=\u_{\mathfrak{s}}(x)-e^{\sigma\mathfrak{s}}\textbf{g}(x)\y(\omega), \quad\qquad x\in \mathbb{R}^d \text{ and }\mathfrak{s}\in\R.
		\end{aligned}
		\right.
	\end{equation}
	For all $\mathfrak{s}\in\R,$ $t>\mathfrak{s},$ and for every $\v_{\mathfrak{s}}\in\H$ and $\omega\in\Omega$, \eqref{CSCBF_Add} has a unique solution $\v(\cdot,\mathfrak{s},\omega,\v_{\mathfrak{s}})\in \mathrm{C}([\mathfrak{s},\mathfrak{s}+T];\H)\cap\mathrm{L}^2(\mathfrak{s}, \mathfrak{s}+T;\V)\cap\mathrm{L}^{r+1}(\mathfrak{s},\mathfrak{s}+T;\widetilde{\L}^{r+1})$. Moreover, $\v(t,\mathfrak{s},\omega,\v_{\mathfrak{s}})$ is continuous with respect to initial data $\v_{\mathfrak{s}}(x)$ (Lemma \ref{Continuity_add}) and $(\mathscr{F},\mathscr{B}(\H))$-measurable in $\omega\in\Omega.$ Define a cocycle $\Phi_0:\R^+\times\R\times\Omega\times\H\to\H$ for the system \eqref{SCBF_Add} such that for given $t\in\R^+, \mathfrak{s}\in\R, \omega\in\Omega$ and $\u_{\mathfrak{s}}\in\H$,
	\begin{align}\label{Phi_0}
		\Phi_0(t,\mathfrak{s},\omega,\u_{\mathfrak{s}}) &=\u(t+\mathfrak{s},\mathfrak{s},\vartheta_{-\mathfrak{s}}\omega,\u_{\mathfrak{s}})=\v(t+\mathfrak{s},\mathfrak{s},\vartheta_{-\mathfrak{s}}\omega,\u_{\mathfrak{s}})+e^{\sigma(t+\mathfrak{s})}\textbf{g}\y(\vartheta_{t}\omega).
	\end{align}
	
	\begin{lemma}\label{Continuity_add}
		For $d=2$ with $r>1$, $d=3$ with $r>3$ and $d=r=3$ with $2\beta\mu\geq1$, assume that $\f\in \mathrm{L}^2_{\emph{loc}}(\mathbb{R};\V')$. Then, the solution of \eqref{CSCBF_Add} is continuous in initial data $\v_{\mathfrak{s}}(x).$
	\end{lemma}
	\begin{proof}
		See the proofs of Theorem 4.8 in \cite{KM} and Theorem 4.10 in \cite{KM3} for $d=2$ and $d=3$, respectively.
	\end{proof}
	Next, we prove the existence of $\mathfrak{D}$-pullback absorbing set of $\Phi_0.$
	\begin{lemma}\label{LemmaUe_add}
		For $d=2$ with $r>1$, $d=3$ with $r>3$ and $d=r=3$ with $2\beta\mu\geq1$, assume that $\f\in \mathrm{L}^2_{\emph{loc}}(\mathbb{R};\V')$ satisfies \eqref{forcing7}. Then $\Phi_0$ possesses a closed measurable $\mathfrak{D}$-pullback absorbing set $\mathcal{K}^1_0=\{\mathcal{K}^1_0(\mathfrak{s},\omega):\mathfrak{s}\in\R, \omega\in\Omega\}\in\mathfrak{D}$ in $\H$ given by
		\begin{align}\label{ue_add}
			\mathcal{K}^1_0(\mathfrak{s},\omega)=\{\u\in\H:\|\u\|^2_{\H}\leq \mathcal{R}^1_0(\mathfrak{s},\omega)\},\ \  \text{ for } d=2, 
		\end{align}
		where $\mathcal{R}^1_0(\mathfrak{s},\omega)$ is defined by
		\begin{align}\label{ue_add1}
			\mathcal{R}^1_0(\mathfrak{s},\omega)&=3\|\textbf{g}\|^2_{\H}\left|e^{\sigma\mathfrak{s}}\y(\omega)\right|^2+2R_5 \int_{-\infty}^{0}e^{\alpha\xi}\bigg[\|\f(\cdot,\xi+\mathfrak{s})\|^2_{\V'}+\left|e^{\sigma (\xi+\mathfrak{s})}\y(\vartheta_{\xi}\omega)\right|^2\nonumber\\&\quad+\left|e^{\sigma (\xi+\mathfrak{s})}\y(\vartheta_{\xi}\omega)\right|^{r+1}+\left|e^{\sigma(\xi+\mathfrak{s})}\y(\vartheta_{\xi}\omega)\right|^{\frac{2(r+1)}{r-1}}\bigg]\d\xi,
		\end{align}
		and a closed measurable $\mathfrak{D}$-pullback absorbing set $\mathcal{K}^2_0=\{\mathcal{K}^2_0(\mathfrak{s},\omega):\mathfrak{s}\in\R, \omega\in\Omega\}\in\mathfrak{D}$ in $\H$ given by
		\begin{align}\label{ue_add2}
			\mathcal{K}^2_0(\mathfrak{s},\omega)=\{\u\in\H:\|\u\|^2_{\H}\leq \mathcal{R}^2_0(\mathfrak{s},\omega)\},\ \  \text{ for } d=3, 
		\end{align}
		where $\mathcal{R}^2_0(\mathfrak{s},\omega)$ is defined by
		\begin{align}\label{ue_add3}
			\mathcal{R}^2_0(\mathfrak{s},\omega)&=3\|\textbf{g}\|^2_{\H}\left|e^{\sigma\mathfrak{s}}\y(\omega)\right|^2+2R_8 \int_{-\infty}^{0}e^{\alpha\xi}\bigg[\|\f(\cdot,\xi+\mathfrak{s})\|^2_{\V'}+\left|e^{\sigma (\xi+\mathfrak{s})}\y(\vartheta_{\xi}\omega)\right|^2\nonumber\\&\quad+\left|e^{\sigma (\xi+\mathfrak{s})}\y(\vartheta_{\xi}\omega)\right|^{r+1}\bigg]\d\xi.
		\end{align}
		Here $R_5$ and $R_8$, both are positive constants do not depend on $\mathfrak{s}$ and $\omega$.
	\end{lemma}
	\begin{proof}
		We infer from \eqref{CSCBF_Add} that
		\begin{align}\label{ue17}
			\frac{1}{2}\frac{\d}{\d t}\|\v\|^2_{\H}=&-\mu\|\nabla\v\|^2_{\H}-\alpha\|\v\|^2_{\H}-b(\v+e^{\sigma t}\textbf{g}\y,\v+e^{\sigma t}\textbf{g}\y,\v)-\beta\left\langle\mathcal{C}(\v+e^{\sigma t}\textbf{g}\y),\v\right\rangle\nonumber\\&+\left\langle\f,\v\right\rangle+e^{\sigma t}\y\left(\left(\ell-\sigma-\alpha\right)\textbf{g}-\mu \A\textbf{g},\v\right)\nonumber\\=&-\mu\|\nabla\v\|^2_{\H}-\alpha\|\v\|^2_{\H}-\beta\|\v+e^{\sigma t}\textbf{g}\y\|^{r+1}_{\wi\L^{r+1}}+b(\v+e^{\sigma t}\textbf{g}\y,\v+e^{\sigma t}\textbf{g}\y,e^{\sigma t}\textbf{g}\y)\nonumber\\&+\beta\left\langle\mathcal{C}(\v+e^{\sigma t}\textbf{g}\y),e^{\sigma t}\textbf{g}\y\right\rangle+\left\langle\f,\v\right\rangle+e^{\sigma t}\y\left(\left(\ell-\sigma-\alpha\right)\textbf{g}-\mu \A\textbf{g},\v\right).
		\end{align}
		Applying Lemmas \ref{Holder} and \ref{Young}, we get that there exists a constant $R_1,R_2,R_3>0$ such that
		\begin{align}
			\beta\left\langle\mathcal{C}(\v+e^{\sigma t}\textbf{g}\y),e^{\sigma t}\textbf{g}\y\right\rangle&\leq\beta \left|e^{\sigma t}\y\right|\|\v+e^{\sigma t}\textbf{g}\y\|^{r}_{\wi\L^{r+1}}\|\textbf{g}\|_{\wi\L^{r+1}}\nonumber\\&\leq\frac{\beta}{4}\|\v+e^{\sigma t}\textbf{g}\y\|^{r+1}_{\wi\L^{r+1}} + R_1\left|e^{\sigma t}\y\right|^{r+1},\label{ue18}\\
			\left\langle\f,\v\right\rangle&\leq\|\f\|_{\V'}\|\v\|_{\V}\leq\frac{\min\{\alpha,\mu\}}{6}\|\v\|^2_{\V}+R_2\|\f\|^{2}_{\V'},\label{ue19}\\
			e^{\sigma t}\y\left(\left(\ell-\sigma-\alpha\right)\textbf{g}-\mu \A\textbf{g},\v\right)&= \left(\ell-\sigma-\alpha\right)e^{\sigma t}\y\left(\textbf{g},\v\right)+\mu e^{\sigma t}\y\left(\nabla\textbf{g},\nabla\v\right)\nonumber\\&\leq\frac{\alpha}{6}\|\v\|^2_{\H}+\frac{\mu}{6}\|\nabla\v\|^2_{\H}+R_3\left|e^{\sigma t}\y\right|^2.\label{ue20}
		\end{align}
		\vskip 2mm
		\noindent
		\textbf{Case I:} \textit{When $d=2$ and $r>1$.}
		Using \eqref{b0}, Lemmas \ref{Holder} and \ref{Young}, and Sobolev's inequality, we obtain that there exists a constant $R_4>0$ such that
		\begin{align}\label{ue21}
			&\left|b(\v+e^{\sigma t}\textbf{g}\y,\v+e^{\sigma t}\textbf{g}\y,e^{\sigma t}\textbf{g}\y)\right|\nonumber\\&= \left|b(\v+e^{\sigma t}\textbf{g}\y,e^{\sigma t}\textbf{g}\y,\v)\right|\nonumber\\&\leq \left|e^{\sigma t}\y\right|\|\v+e^{\sigma t}\textbf{g}\y\|_{\wi\L^{r+1}}\|\nabla\textbf{g}\|_{\H}\|\v\|_{\wi\L^{\frac{2(r+1)}{r-1}}}\nonumber\\&\leq\frac{\beta}{4}\|\v+e^{\sigma t}\textbf{g}\y\|^{r+1}_{\wi\L^{r+1}} +\frac{\min\{\alpha,\mu\}}{6}\|\v\|^{2}_{\V}+R_4\left|e^{\sigma t}\y\right|^{\frac{2(r+1)}{r-1}}.
		\end{align}
		Combining \eqref{ue17}-\eqref{ue21}, we get
		\begin{align}\label{ue22}
			\frac{\d}{\d t}\|\v\|^2_{\H}+\alpha\|\v\|^2_{\H}\leq R_5\left[\|\f\|^{2}_{\V'}+\left|e^{\sigma t}\y\right|^2+\left|e^{\sigma t}\y\right|^{r+1}+\left|e^{\sigma t}\y\right|^{\frac{2(r+1)}{r-1}}\right],
		\end{align}
		where $R_5=\max\{2R_1,2R_2,2R_3,2R_4\}$. Applying variation of constant formula to \eqref{ue22} over $(\mathfrak{s}-t,\tau)$ with $t\geq0$, $\mathfrak{s}\in\R$ and $\tau\geq\mathfrak{s}-t$, and replacing $\omega$ by $\vartheta_{-\mathfrak{s}}\omega$, we obtain 
		\begin{align}\label{ue23}
			&\|\v(\tau,\mathfrak{s}-t,\vartheta_{-\mathfrak{s}}\omega,\v_{\mathfrak{s}-t})\|^2_{\H} \nonumber\\&\leq e^{\alpha(\mathfrak{s}-t-\tau)}\|\v_{\mathfrak{s}-t}\|^2_{\H}+R_5 \int_{-\infty}^{\tau-\mathfrak{s}}e^{\alpha(\xi+\mathfrak{s}-\tau)}\bigg[\|\f(\cdot,\xi+\mathfrak{s})\|^2_{\V'}\nonumber\\&\qquad+\left|e^{\sigma (\xi+\mathfrak{s})}\y(\vartheta_{\xi}\omega)\right|^2+\left|e^{\sigma (\xi+\mathfrak{s})}\y(\vartheta_{\xi}\omega)\right|^{r+1}+\left|e^{\sigma(\xi+\mathfrak{s})}\y(\vartheta_{\xi}\omega)\right|^{\frac{2(r+1)}{r-1}}\bigg]\d\xi.
		\end{align}
		From \eqref{Phi_0}, we have
		\begin{align*}
			\u(\mathfrak{s},\mathfrak{s}-t,\vartheta_{-\mathfrak{s}}\omega,\u_{\mathfrak{s}-t})=\v(\mathfrak{s},\mathfrak{s}-t,\vartheta_{-\mathfrak{s}}\omega,\v_{\mathfrak{s}-t})+e^{\sigma\mathfrak{s}}\y(\omega)\textbf{g},
		\end{align*}
		with $\v_{\mathfrak{s}-t}=\u_{\mathfrak{s}-t}-e^{\sigma(\mathfrak{s}-t)}\y(\vartheta_{-t}\omega)\textbf{g}$, $\u_{\mathfrak{s}-t}\in D(\mathfrak{s}-t,\vartheta_{-t}\omega)$ and $D\in\mathfrak{D}$, which together with \eqref{ue23} gives
		\begin{align}\label{ue24}
			&	\|\u(\tau,\mathfrak{s}-t,\vartheta_{-\mathfrak{s}}\omega,\u_{\mathfrak{s}-t})\|^2_{\H}\nonumber\\&\leq2\|\v(\tau,\mathfrak{s}-t,\vartheta_{-\mathfrak{s}}\omega,\v_{\mathfrak{s}-t})\|^2_{\H}+2\|\textbf{g}\|^2_{\H}\left|e^{\sigma\mathfrak{s}}\y(\omega)\right|^2 \nonumber\\&\leq 4e^{-\alpha t}\left(\|\u_{\mathfrak{s}-t}\|^2_{\H}+\|\textbf{g}\|^2_{\H}\left|e^{\sigma(\mathfrak{s}-t)}\y(\vartheta_{-t}\omega)\right|^2\right)+2\|\textbf{g}\|^2_{\H}\left|e^{\sigma\mathfrak{s}}\y(\omega)\right|^2+2R_5 \int_{-\infty}^{0}e^{\alpha\xi}\nonumber\\&\quad\times\bigg[\|\f(\cdot,\xi+\mathfrak{s})\|^2_{\V'}+\left|e^{\sigma (\xi+\mathfrak{s})}\y(\vartheta_{\xi}\omega)\right|^2+\left|e^{\sigma (\xi+\mathfrak{s})}\y(\vartheta_{\xi}\omega)\right|^{r+1}+\left|e^{\sigma(\xi+\mathfrak{s})}\y(\vartheta_{\xi}\omega)\right|^{\frac{2(r+1)}{r-1}}\bigg]\d\xi.
		\end{align}
		Using \eqref{forcing6} and \eqref{Y2}, and arguing similarly as in the proofs of Lemmas \ref{LemmaUe3} and \ref{PAS_GA}, one can  conclude the proof.
		\vskip 2mm
		\noindent
		\textbf{Case II:} \textit{When $d= 3$ and $r\geq3$ ($r>3$ with any $\beta,\mu>0$ and $r=3$ with $2\beta\mu\geq1$).}
		Using \eqref{b0}, Lemmas \ref{Holder}, \ref{Interpolation} and \ref{Young}, we obtain that there exists two constants $R_6,R_7>0$ such that
		\begin{align}\label{ue25}
			&\left|b(\v+e^{\sigma t}\textbf{g}\y,\v+e^{\sigma t}\textbf{g}\y,e^{\sigma t}\textbf{g}\y)\right|\nonumber\\&=\left|b(\v+e^{\sigma t}\textbf{g}\y,e^{\sigma t}\textbf{g}\y,\v+e^{\sigma t}\textbf{g}\y)\right|\nonumber\\&\leq \left|e^{\sigma t}\y\right|\|\v+e^{\sigma t}\textbf{g}\y\|_{\wi\L^{r+1}}\|\nabla\textbf{g}\|_{\H}\|\v+e^{\sigma t}\textbf{g}\y\|_{\wi\L^{\frac{2(r+1)}{r-1}}}\nonumber\\&\leq\left|e^{\sigma t}\y\right|\|\v+e^{\sigma t}\textbf{g}\y\|^{\frac{r+1}{r-1}}_{\wi\L^{r+1}}\|\nabla\textbf{g}\|_{\H}\|\v+e^{\sigma t}\textbf{g}\y\|^{\frac{r-3}{r-1}}_{\H}\nonumber\\&\leq\frac{\beta}{4}\|\v+e^{\sigma t}\textbf{g}\y\|^{r+1}_{\wi\L^{r+1}} +\frac{\alpha}{6}\|\v\|^{2}_{\H}+R_6\left|e^{\sigma t}\y\right|^{2}, \ \ \text{ for } r>3,
		\end{align}
		and 
		\begin{align}\label{ue26}
			&\left|b(\v+e^{\sigma t}\textbf{g}\y,\v+e^{\sigma t}\textbf{g}\y,e^{\sigma t}\textbf{g}\y)\right|\nonumber\\&=\left|b(\v+e^{\sigma t}\textbf{g}\y,e^{\sigma t}\textbf{g}\y,\v+e^{\sigma t}\textbf{g}\y)\right|\nonumber\\&\leq \left|e^{\sigma t}\y\right|\|\v+e^{\sigma t}\textbf{g}\y\|^2_{\wi\L^{4}}\|\nabla\textbf{g}\|_{\H}\nonumber\\&\leq\frac{\beta}{4}\|\v+e^{\sigma t}\textbf{g}\y\|^{4}_{\wi\L^{4}} +R_7\left|e^{\sigma t}\y\right|^{2}, \ \ \text{ for } r=3.
		\end{align}
		Combining \eqref{ue17}-\eqref{ue20} and \eqref{ue25}-\eqref{ue26}, we get
		\begin{align}\label{ue27}
			\frac{\d}{\d t}\|\v\|^2_{\H}+\alpha\|\v\|^2_{\H}\leq R_8\left[\|\f\|^{2}_{\V'}+\left|e^{\sigma t}\y\right|^2+\left|e^{\sigma t}\y\right|^{r+1}\right],
		\end{align}
		where $R_8=\max\{2R_1,2R_2,2(R_3+R_6),2(R_3+R_7)\}$. Applying variation of constant formula to \eqref{ue22} over $(\mathfrak{s}-t,\tau)$ with $t\geq0$, $\mathfrak{s}\in\R$ and $\tau\geq\mathfrak{s}-t$, and replacing $\omega$ by $\vartheta_{-\mathfrak{s}}\omega$, we obtain 
		\begin{align}\label{ue28}
			\|\v(\tau,\mathfrak{s}-t,\vartheta_{-\mathfrak{s}}\omega,\v_{\mathfrak{s}-t})\|^2_{\H} &\leq e^{\alpha(\mathfrak{s}-t-\tau)}\|\v_{\mathfrak{s}-t}\|^2_{\H}+R_8 \int_{-\infty}^{\tau-\mathfrak{s}}e^{\alpha(\xi+\mathfrak{s}-\tau)}\bigg[\|\f(\cdot,\xi+\mathfrak{s})\|^2_{\V'}\nonumber\\&\quad+\left|e^{\sigma (\xi+\mathfrak{s})}\y(\vartheta_{\xi}\omega)\right|^2+\left|e^{\sigma (\xi+\mathfrak{s})}\y(\vartheta_{\xi}\omega)\right|^{r+1}\bigg]\d\xi,
		\end{align}
		or
		\begin{align}\label{ue29}
			&	\|\u(\tau,\mathfrak{s}-t,\vartheta_{-\mathfrak{s}}\omega,\u_{\mathfrak{s}-t})\|^2_{\H}\nonumber\\&\leq2\|\v(\tau,\mathfrak{s}-t,\vartheta_{-\mathfrak{s}}\omega,\v_{\mathfrak{s}-t})\|^2_{\H}+2\|\textbf{g}\|^2_{\H}\left|e^{\sigma\mathfrak{s}}\y(\omega)\right|^2 \nonumber\\&\leq 4e^{-\alpha t}\left(\|\u_{\mathfrak{s}-t}\|^2_{\H}+\|\textbf{g}\|^2_{\H}\left|e^{\sigma(\mathfrak{s}-t)}\y(\vartheta_{-t}\omega)\right|^2\right)+2\|\textbf{g}\|^2_{\H}\left|e^{\sigma\mathfrak{s}}\y(\omega)\right|^2\nonumber\\&\quad+2R_8 \int_{-\infty}^{0}e^{\alpha\xi}\bigg[\|\f(\cdot,\xi+\mathfrak{s})\|^2_{\V'}+\left|e^{\sigma (\xi+\mathfrak{s})}\y(\vartheta_{\xi}\omega)\right|^2+\left|e^{\sigma (\xi+\mathfrak{s})}\y(\vartheta_{\xi}\omega)\right|^{r+1}\bigg]\d\xi,
		\end{align}
		where $\u_{\mathfrak{s}-t}\in D(\mathfrak{s}-t,\vartheta_{-t}\omega)$ and $D\in\mathfrak{D}$. Using \eqref{forcing6} and \eqref{Y2}, and arguing similarly as in the proofs of Lemmas \ref{LemmaUe3} and \ref{PAS_GA}, we conclude the proof.
	\end{proof}
	The next two Corollaries \ref{convergence_b} and \ref{convergence_c} help us to prove the asymptotic compactness of $\Phi_0$ in this subsection (see Lemma \ref{Asymptotic_UB_Add}) as well as the uniform compactness of random pullback attractors in next subsection (see Lemma \ref{precompact}). 
	\begin{corollary}\label{convergence_b}
		Let $\{\v_m\}_{m\in \mathbb{N}}\subset\mathrm{L}^{\infty}(0, T; \H)\cap\mathrm{L}^2(0, T;\V)\cap\mathrm{L}^{r+1}(0, T;\wi\L^{r+1})$ be a bounded sequence and   $\v \in \mathrm{L}^{\infty}(0, T; \H)\cap\mathrm{L}^2(0, T;\V)\cap\mathrm{L}^{r+1}(0, T;\wi\L^{r+1})$  such that 
		\begin{equation}\label{5.26}
			\left\{
			\begin{aligned}
				\v_m&\xrightharpoonup{w}\v \ \text{ in }\ \mathrm{L}^2(0, T;\V),\\
				\v_m&\xrightharpoonup{w}\v \ \text{ in }\  \mathrm{L}^{r+1}(0, T;\wi\L^{r+1}),\\
				\v_m&\to\v\ \text{ in }\ \mathrm{L}^2(0, T; {\L}^2_{\emph{loc}}(\R^d)).
			\end{aligned}
			\right.
		\end{equation}
		Then, for $d=2$ with $r>1$, $d=3$ with $r>3$ and $d=r=3$ with $2\beta\mu\geq1$, 
		\begin{align}\label{Conver}
			\int_{0}^{T} b(\v_m(t), \v_m(t), \uprho(t)\emph{\textbf{g}}) \d t \to \int_{0}^{T} b(\v(t), \v(t), \uprho(t)\emph{\textbf{g}}) \d t\  \ \text{ as }\  \ m\to\infty,
		\end{align} 
		where $\uprho(t)$ is a continuous function of $t$ on $[0,T]$.
	\end{corollary}
	\begin{proof}
		Since $\v_m, \v \in\mathrm{L}^{\infty}(0, T; \H)\cap\mathrm{L}^2(0, T;\V)\cap\mathrm{L}^{r+1}(0, T;\wi\L^{r+1})$ and $\uprho(t)$ is continuous, there exists a constant $\varpi>0$ such that
		\begin{align*}
			&\sup_{t  \in [0, T]} \|\v_m(t)\|_{\H} +\sup_{t  \in [0, T]} \|\v(t)\|_{\H} + \left(\int_{0}^{T}\|\v_m(t)\|^2_{\V} \d t\right)^{\frac{1}{2}} +  \left(\int_{0}^{T}\|\v(t)\|^2_{\V}\d t\right)^{\frac{1}{2}}\\&+\left(\int_{0}^{T}\|\v_m(t)\|^{r+1}_{\wi\L^{r+1}} \d t\right)^{\frac{1}{r+1}} +  \left(\int_{0}^{T}\|\v(t)\|^{r+1}_{\wi\L^{r+1}}\d t\right)^{\frac{1}{r+1}} + \left(\int_{0}^{T}\left|\uprho(t)\right|^2\d t\right)^{\frac{1}{2}} \leq \varpi.
		\end{align*}
		Since $\textbf{g}\in\D(\A)$, for given $\varepsilon>0$, we can say that there exists $k=k(\varepsilon,\textbf{g})>0$ such that
		\begin{align}\label{Conver1}
			\int_{\R^d\backslash\mathcal{O}_k}\left(|\textbf{g}(x)|^2+|\nabla\textbf{g}(x)|^2\right)\d x \leq \frac{\varepsilon^2}{256 \varpi^6}, \ \ \ \text{ for } d=2,
		\end{align}
		and
		\begin{align}\label{Conver1*}
			\int_{\R^d\backslash\mathcal{O}_k}\left(|\textbf{g}(x)|^2+|\nabla\textbf{g}(x)|^2\right)\d x \leq \frac{\varepsilon^2}{64 \varpi^6}, \ \ \ \text{ for }  d=3,
		\end{align}
		where $\mathcal{O}_k=\{x\in\R^d:|x|<k\}$. Let us consider, for $\u,\v,\w\in\V$, 
		$$b(\u,\v,\w)=b_1(\u,\v,\w)+b_2(\u,\v,\w)$$ where
		$$b_1(\u,\v,\w)=\sum_{i,j=1}^d\int_{\mathcal{O}_k}\u_i\frac{\partial \v_j}{\partial x_i}\w_j\d x \ \ \text{ and } \ \  b_2(\u,\v,\w)=\sum_{i,j=1}^d\int_{\R^d\backslash\mathcal{O}_k}\u_i\frac{\partial \v_j}{\partial x_i}\w_j\d x.$$
		Consider,
		\begin{align}\label{Conver2}
			&\int_{0}^{T} b(\v_m(t), \v_m(t), \uprho(t)\textbf{g}) \d t-\int_{0}^{T} b(\v(t), \v(t), \uprho(t)\textbf{g}) \d t \nonumber\\& =\int_{0}^{T} b(\v_m(t)-\v(t), \v_m(t), \uprho(t)\textbf{g}) \d t +\int_{0}^{T} b(\v(t), \v_m(t)-\v(t), \uprho(t)\textbf{g}) \d t\nonumber\\&=\int_{0}^{T} b_1(\v_m(t)-\v(t), \v_m(t), \uprho(t)\textbf{g}) \d t +\int_{0}^{T} b_2(\v_m(t)-\v(t), \v_m(t), \uprho(t)\textbf{g}) \d t \nonumber\\&\qquad+\int_{0}^{T} b_1(\v(t), \v_m(t)-\v(t), \uprho(t)\textbf{g}) \d t +\int_{0}^{T} b_2(\v(t), \v_m(t)-\v(t), \uprho(t)\textbf{g}) \d t .
		\end{align}
		In order to prove \eqref{Conver}, it is enough to show that each term on the right hand side of \eqref{Conver2} converges to $0$ as $m\to\infty$.
		\vskip 2mm
		\noindent
		\textbf{Case I:} \textit{$d=2$ and $r>1$.}
		By \eqref{b1} and \eqref{Conver1}, we obtain
		\begin{align}\label{Conver3}
			&\left|\int_{0}^{T} b_2( \v_m(t)-\v(t),\v_m(t), \uprho(t)\textbf{g}) \d t\right|\nonumber\\&\leq\frac{\varepsilon}{8\varpi^3}\int_{0}^{T}\left|\uprho(t)\right|\left[\|\v_m(t)\|_{\H}\|\nabla\v_m(t)\|_{\H}+\|\v(t)\|^{\frac{1}{2}}_{\H}\|\nabla\v(t)\|^{\frac{1}{2}}_{\H}\|\v_m(t)\|^{\frac{1}{2}}_{\H}\|\nabla\v_m(t)\|^{\frac{1}{2}}_{\H}\right]\d t \nonumber\\&\leq\frac{\varepsilon}{8\varpi^3}\sup_{t  \in [0, T]} \|\v_m(t)\|_{\H}\left(\int_{0}^{T}\left|\uprho(t)\right|^2\d t\right)^{\frac{1}{2}}\left(\int_{0}^{T}\|\nabla\v_m(t)\|^2_{\H}\d t\right)^{\frac{1}{2}}+ \frac{\varepsilon}{8\varpi^3}\sup_{t  \in [0, T]} \|\v(t)\|^{\frac{1}{2}}_{\H}\nonumber\\&\quad\times\sup_{t  \in [0, T]} \|\v_m(t)\|^{\frac{1}{2}}_{\H}\left(\int_{0}^{T}\left|\uprho(t)\right|^2\d t\right)^{\frac{1}{2}}\left(\int_{0}^{T}\|\nabla\v(t)\|^2_{\H}\d t\right)^{\frac{1}{4}} \left(\int_{0}^{T}\|\nabla\v_m(t)\|^2_{\H}\d t\right)^{\frac{1}{4}}\nonumber\\&\leq \frac{\varepsilon}{4}.
		\end{align}
		By \eqref{b1} and continuity of $\uprho(t)$, we have
		\begin{align}\label{Conver4}
			&\left|\int_{0}^{T} b_1( \v_m(t)-\v(t),\v_m(t), \uprho(t)\textbf{g}) \d t\right|\nonumber\\&\leq C\|\nabla\textbf{g}\|_{\H}\int_{0}^{T}\left|\uprho(t)\right|\|\v_m(t)-\v(t)\|^{\frac{1}{2}}_{\L^2(\mathcal{O}_k)}\|\nabla\left(\v_m(t)-\v(t)\right)\|^{\frac{1}{2}}_{\H}\|\v_m(t)\|^{\frac{1}{2}}_{\H}\|\nabla\v_m(t)\|^{\frac{1}{2}}_{\H}\d t\nonumber\\&\leq CT^{\frac{1}{4}}\sup_{t  \in [0, T]}\left[\left|\uprho(t)\right|\|\v_m(t)\|^{\frac{1}{2}}_{\H}\right]\left(\int_{0}^{T}\|\v_m(t)-\v(t)\|^2_{\L^2(\mathcal{O}_k)}\d t\right)^{\frac{1}{4}}\bigg[\left(\int_{0}^{T}\|\nabla\v_m(t)\|^2_{\H}\d t\right)^{\frac{1}{2}}\nonumber\\&\qquad\qquad+\left(\int_{0}^{T}\|\nabla\v(t)\|^2_{\H}\d t\right)^{\frac{1}{4}}\left(\int_{0}^{T}\|\nabla\v_m(t)\|^2_{\H}\d t\right)^{\frac{1}{4}}\bigg]\nonumber\\&\to0  \ \text{ as }\   m\to \infty.
		\end{align}
		\vskip 2mm
		\noindent
		\textbf{Case II:} \textit{$d= 3$ and $r\geq3$ ($r>3$ with any $\beta,\mu>0$ and $r=3$ with $2\beta\mu\geq1$).} By \eqref{Conver1*}, Lemmas \ref{Holder} and \ref{Interpolation}, we have
		\begin{align}\label{Conver9}
			&\left|\int_{0}^{T} b_2( \v_m(t)-\v(t),\v_m(t), \uprho(t)\textbf{g}) \d t\right|\nonumber\\&\leq\frac{\varepsilon}{8\varpi^3}\int_{0}^{T}\left|\uprho(t)\right|\left[\|\v_m(t)\|_{\wi\L^{r+1}}\|\v_m(t)\|_{\wi\L^{\frac{2(r+1)}{r-1}}}+\|\v(t)\|_{\wi\L^{r+1}}\|\v_m(t)\|_{\wi\L^{\frac{2(r+1)}{r-1}}}\right]\d t\nonumber\\&\leq\frac{\varepsilon}{8\varpi^3}\int_{0}^{T}\left|\uprho(t)\right|\left[\|\v_m(t)\|^{\frac{r+1}{r-1}}_{\wi\L^{r+1}}\|\v_m(t)\|^{\frac{r-3}{r-1}}_{\H}+\|\v(t)\|_{\wi\L^{r+1}}\|\v_m(t)\|^{\frac{2}{r-1}}_{\wi\L^{r+1}}\|\v_m(t)\|^{\frac{r-3}{r-1}}_{\H}\right]\d t \nonumber\\&\leq\frac{\varepsilon}{8\varpi^3}\left(\int_{0}^{T}\left|\uprho(t)\right|^2\d t\right)^{\frac{1}{2}}\left(\int_{0}^{T}\|\v_m(t)\|^2_{\H}\d t\right)^{\frac{r-3}{2(r-1)}}\bigg[\left(\int_{0}^{T}\|\v_m(t)\|^{r+1}_{\wi\L^{r+1}}\d t\right)^{\frac{1}{r-1}}\nonumber\\&\quad+ \left(\int_{0}^{T}\|\v(t)\|^{r+1}_{\wi\L^{r+1}}\d t\right)^{\frac{1}{r+1}}\left(\int_{0}^{T}\|\v_m(t)\|^{r+1}_{\wi\L^{r+1}}\d t\right)^{\frac{2}{(r+1)(r-1)}}\bigg]\nonumber\\&\leq \frac{\varepsilon}{4},
		\end{align}
		for $r>3$ and
		\begin{align}\label{Conver10}
			&\left|\int_{0}^{T} b_2( \v_m(t)-\v(t),\v_m(t), \uprho(t)\textbf{g}) \d t\right|\nonumber\\&\leq\frac{\varepsilon}{8\varpi^3}\int_{0}^{T}\left|\uprho(t)\right|\|\v_m(t)\|^2_{\wi\L^{4}}+ \frac{\varepsilon}{8\varpi^3}\int_{0}^{T}\left|\uprho(t)\right|\|\v(t)\|_{\wi\L^{4}}\|\v_m(t)\|_{\wi\L^4}\d t\nonumber\\&\leq\frac{\varepsilon}{8\varpi^3}\left(\int_{0}^{T}\left|\uprho(t)\right|^2\d t\right)^{\frac{1}{2}}\left[\left(\int_{0}^{T}\|\v_m(t)\|^{4}_{\wi\L^{4}}\d t\right)^{\frac{1}{2}}+\left(\int_{0}^{T}\|\v(t)\|^{4}_{\wi\L^{4}}\d t\right)^{\frac{1}{4}}\left(\int_{0}^{T}\|\v_m(t)\|^{4}_{\wi\L^{4}}\d t\right)^{\frac{1}{4}}\right]\nonumber\\&\leq \frac{\varepsilon}{4}, 
		\end{align}
		for $r=3$.	By \eqref{b1} and continuity of $\uprho(t)$, we have
		\begin{align}\label{Conver11}
			&\left|\int_{0}^{T} b_1( \v_m(t)-\v(t),\v_m(t), \uprho(t)\textbf{g}) \d t\right|\nonumber\\&\leq C\|\nabla\textbf{g}\|_{\H}\int_{0}^{T}\left|\uprho(t)\right|\|\v_m(t)-\v(t)\|^{\frac{1}{4}}_{\L^2(\mathcal{O}_k)}\|\nabla\left(\v_m(t)-\v(t)\right)\|^{\frac{3}{4}}_{\H}\|\v_m(t)\|^{\frac{1}{4}}_{\H}\|\nabla\v_m(t)\|^{\frac{3}{4}}_{\H}\d t\nonumber\\&\leq CT^{\frac{1}{8}}\sup_{t  \in [0, T]}\left[\left|\uprho(t)\right|\cdot\|\v_m(t)\|^{\frac{1}{4}}_{\H}\right]\left(\int_{0}^{T}\|\v_m(t)-\v(t)\|^2_{\L^2(\mathcal{O}_k)}\d t\right)^{\frac{1}{8}}\bigg[\left(\int_{0}^{T}\|\nabla\v_m(t)\|^2_{\H}\d t\right)^{\frac{3}{4}}\nonumber\\&\quad+\left(\int_{0}^{T}\|\nabla\v(t)\|^2_{\H}\d t\right)^{\frac{3}{8}}\left(\int_{0}^{T}\|\nabla\v_m(t)\|^2_{\H}\d t\right)^{\frac{3}{8}}\bigg]\nonumber\\&\to0  \ \text{ as }\   m\to \infty.
		\end{align}
		Calculations similar to \eqref{Conver3} (for $d=2$ and $r>1$), \eqref{Conver9} (for $d=3$ and $r>3$) and \eqref{Conver10} (for $d=r=3$), and \eqref{Conver4} (for $d=2$ and $r>1$) and \eqref{Conver11} (for $d=3$ and $r\geq3$), we obtain 
		\begin{align}\label{Conver5}
			&\left|\int_{0}^{T} b_2( \v(t),\v_m(t)-\v(t), \uprho(t)\textbf{g}) \d t\right|\leq \frac{\varepsilon}{4},
		\end{align}
		and
		\begin{align}\label{Conver6}
			\left|\int_{0}^{T} b_1(\v(t), \v_m(t)-\v(t), \uprho(t)\textbf{g}) \d t\right|\to0  \ \text{ as }\   m\to \infty,
		\end{align}
		respectively. For given $\varepsilon>0$, we infer from \eqref{Conver4}  (for $d=2$ and $r>1$) and \eqref{Conver11} (for $d=3$ and $r\geq3$), and \eqref{Conver6} that there exists $M(\varepsilon)\in\N$ such that
		\begin{align}\label{Conver7}
			\left|\int_{0}^{T} b_1( \v_m(t)-\v(t),\v_m(t), \uprho(t)\textbf{g}) \d t\right| \leq \frac{\varepsilon}{4},
		\end{align}
		and
		\begin{align}\label{Conver8}
			\left|\int_{0}^{T} b_1(\v(t), \v_m(t)-\v(t), \uprho(t)\textbf{g}) \d t\right|\leq\frac{\varepsilon}{4},
		\end{align}
		respectively, for all $m\geq M(\varepsilon).$ Hence, \eqref{Conver2}, \eqref{Conver3} (for $d=2$ and $r>1$), \eqref{Conver9} (for $d=3$ and $r>3$), \eqref{Conver10} (for $d=r=3$), \eqref{Conver5} and \eqref{Conver7}-\eqref{Conver8} imply that, for a given $\varepsilon>0$, there exists $M(\varepsilon)\in\N$ such that
		\begin{align*}
			\left|\int_{0}^{T} b(\v_m(t), \v_m(t), \uprho(t)\textbf{g}) \d t-\int_{0}^{T} b(\v(t), \v(t), \uprho(t)\textbf{g}) \d t \right|\leq \varepsilon,
		\end{align*}
		for all $m\geq M(\varepsilon),$ which completes the proof.
	\end{proof}

	\begin{corollary}\label{convergence_c}
		Let $\{\v_m\}_{m\in \mathbb{N}}\subset\mathrm{L}^{\infty}(0, T; \H)\cap\mathrm{L}^2(0, T;\V)\cap\mathrm{L}^{r+1}(0, T;\wi\L^{r+1})$ be a bounded sequence and  $\v \in \mathrm{L}^{\infty}(0, T; \H)\cap\mathrm{L}^2(0, T;\V)\cap\mathrm{L}^{r+1}(0, T;\wi\L^{r+1})$ such that \eqref{5.26} be satisfied. Then, for $d=2$ with $r>1$, $d=3$ with $r>3$ and $d=r=3$ with $2\beta\mu\geq1$, 
		\begin{align*}
			\int_{0}^{T} \big\langle\mathcal{C}(\v_m(t)) ,\uprho(t)\emph{\textbf{g}}  \big\rangle \d t \to \int_{0}^{T} \big\langle\mathcal{C}(\v(t)) ,\uprho(t)\emph{\textbf{g}}  \big\rangle \d t,\ \text{ as }\ m\to\infty,
		\end{align*}
		where $\uprho(t)$ is a continuous function of $t$ on $[0,T]$.
	\end{corollary}
	\begin{proof}
		Since $\v_m, \v \in\mathrm{L}^{\infty}(0, T; \H)\cap\mathrm{L}^2(0, T;\V)\cap\mathrm{L}^{r+1}(0, T;\wi\L^{r+1})$ and $\uprho(t)$ is continuous, there exists a constant $\varpi>0$ such that
		\begin{align*}
			\left(\int_{0}^{T}\|\v_m(t)\|^{r+1}_{\wi\L^{r+1}} \d t\right)^{\frac{1}{r+1}} +  \left(\int_{0}^{T}\|\v(t)\|^{r+1}_{\wi\L^{r+1}}\d t\right)^{\frac{1}{r+1}} + \left(\int_{0}^{T}\left|\uprho(t)\right|^{r+1}\d t\right)^{\frac{1}{r+1}} \leq \varpi.
		\end{align*}
		Since $\textbf{g}\in\D(\A)$, Sobolev's embeddings ($\V\subset\widetilde{\L}^{r+1}$ for $d=2$ and $\D(\A)\subset\widetilde{\L}^{r+1}$ for $d=3$) imply that $\textbf{g}\in\widetilde{\L}^{r+1}$, for all $r\geq1$. Hence, for given $\varepsilon>0$, there exists $k=k(\varepsilon,\textbf{g})>0$ such that
		\begin{align}\label{Conver12}
			\left(\int_{\R^d\backslash\mathcal{O}_k}|\textbf{g}(x)|^{r+1}\d x\right)^{\frac{1}{r+1}} \leq \frac{\varepsilon}{4 \varpi^{r+1}},
		\end{align}
		where $\mathcal{O}_k=\{x\in\R^d:|x|<k\}$. Consider,
		\begin{align}\label{Conver13}
			&\int_{0}^{T} \big\langle\mathcal{C}(\v_m(t)) ,\uprho(t)\textbf{g}  \big\rangle \d t- \int_{0}^{T} \big\langle\mathcal{C}(\v(t)) ,\uprho(t)\textbf{g}  \big\rangle \d t\nonumber\\&=\int_{0}^{T}\uprho(t) \left[\int_{\mathbb{R}^d}\left(\left|\v_m(x,t)\right|^{r-1}\v_m(x,t)-\left|\v(x,t)\right|^{r-1}\v(x,t)\right)\textbf{g}(x)  \d x\right] \d t\nonumber\\&=\int_{0}^{T}\uprho(t) \left[\int_{\mathcal{O}_k}\left(\left|\v_m(x,t)\right|^{r-1}\v_m(x,t)-\left|\v(x,t)\right|^{r-1}\v(x,t)\right)\textbf{g}(x)  \d x\right] \d t\nonumber\\&\quad+\int_{0}^{T}\uprho(t) \left[\int_{\mathbb{R}^d\backslash\mathcal{O}_k}\left(\left|\v_m(x,t)\right|^{r-1}\v_m(x,t)-\left|\v(x,t)\right|^{r-1}\v(x,t)\right)\textbf{g}(x)  \d x\right] \d t.
		\end{align}
		From  Lemma \ref{Holder} and \eqref{Conver12}, we infer that 
		\begin{align}\label{Conver14}
			&\left|\int_{0}^{T}\uprho(t) \left[\int_{\mathbb{R}^d\backslash\mathcal{O}_k}\left(\left|\v_m(x,t)\right|^{r-1}\v_m(x,t)-\left|\v(x,t)\right|^{r-1}\v(x,t)\right)\textbf{g}(x)  \d x\right] \d t\right|\nonumber\\&\leq\frac{\varepsilon}{4\varpi^{r+1}}\int_{0}^{T}\left|\uprho(t)\right| \|\v_m(t)\|^r_{\wi\L^{r+1}} \d t +\frac{\varepsilon}{4\varpi^{r+1}}\int_{0}^{T}\left|\uprho(t)\right| \|\v(t)\|^r_{\wi\L^{r+1}} \d t\nonumber\\&\leq\frac{\varepsilon}{4\varpi^{r+1}} \left(\int_{0}^{T}\left|\uprho(t)\right|^{r+1}\d t\right)^{\frac{1}{r+1}}\left[\left(\int_{0}^{T}\|\v_m(t)\|^{r+1}_{\wi\L^{r+1}}\d t\right)^{\frac{r}{r+1}}+\left(\int_{0}^{T}\|\v(t)\|^{r+1}_{\wi\L^{r+1}}\d t\right)^{\frac{r}{r+1}}\right]\nonumber\\&\leq\frac{\varepsilon}{2}.
		\end{align}
		Using Taylor's formula, \eqref{29} and Lemma \ref{Holder}, we achieve
		\begin{align*}
			&\left|\int_{0}^{T}\uprho(t) \left[\int_{\mathcal{O}_k}\left(\left|\v_m(x,t)\right|^{r-1}\v_m(x,t)-\left|\v(x,t)\right|^{r-1}\v(x,t)\right)\textbf{g}(x)  \d x\right] \d t\right|\nonumber\\&\leq C\int_{0}^{T}\left|\uprho(t)\right|\left[\int_{\mathcal{O}_k}\left(\left|\v_m(x,t)\right|^{r-1}+\left|\v(x,t)\right|^{r-1}\right)\left|\v_m(x,t)-\v(x,t)\right|\left|\textbf{g}(x)\right|  \d x\right] \d t\nonumber\\&\leq C\int_{0}^{T}\left|\uprho(t)\right|\left(\|\v_m(t)\|^{r-1}_{\wi\L^{r+1}}+\|\v(t)\|^{r-1}_{\wi\L^{r+1}}\right)\|\v_m(t)-\v(t)\|_{\L^2(\mathcal{O}_k)} \|\textbf{g}\|_{\wi\L^{\frac{2(r+1)}{3-r}}}\d t\nonumber\\&\leq CT^{\frac{3-r}{2(r+1)}}\sup_{t  \in [0, T]}\left|\uprho(t)\right|\left(\|\v_m\|^{r-1}_{\mathrm{L}^{r+1}(0,T;\widetilde{\L}^{r+1})}+\|\v\|^{r-1}_{\mathrm{L}^{r+1}(0,T;\widetilde{\L}^{r+1})}\right)\|\v_m-\v\|_{\mathrm{L}^2(0,T;\L^2(\mathcal{O}_k))}\nonumber\\&
			\to 0  \ \text{ as }\   m\to \infty, 
		\end{align*}
		for $ 1<r<3$,	and
		\begin{align*}
			&\left|\int_{0}^{T}\uprho(t) \left[\int_{\mathcal{O}_k}\left(\left|\v_m(x,t)\right|^{r-1}\v_m(x,t)-\left|\v(x,t)\right|^{r-1}\v(x,t)\right)\textbf{g}(x)  \d x\right] \d t\right|\nonumber\\&\leq C\int_{0}^{T}\left|\uprho(t)\right|\left[\int_{\mathcal{O}_k}\left(\left|\v_m(x,t)\right|^{r-1}+\left|\v(x,t)\right|^{r-1}\right)\left|\v_m(x,t)-\v(x,t)\right|\left|\textbf{g}(x)\right|  \d x\right] \d t\nonumber\\&\leq C\int_{0}^{T}\left|\uprho(t)\right|\left(\|\v_m(t)\|^{r-1}_{\wi\L^{r+1}}+\|\v(t)\|^{r-1}_{\wi\L^{r+1}}\right)\|\v_m(t)-\v(t)\|^{\frac{1}{r-1}}_{\L^{2}(\mathcal{O}_k)}\nonumber\\&\qquad\qquad\times\|\v_m(t)-\v(t)\|^{\frac{r-2}{r-1}}_{\wi\L^{r+1}} \|\textbf{g}\|_{\wi\L^{2(r+1)}}\d t\nonumber\\&\leq CT^{\frac{1}{2(r+1)}}\sup_{t  \in [0, T]}\left|\uprho(t)\right|\left(\|\v_m\|^{r-1}_{\mathrm{L}^{r+1}(0,T;\widetilde{\L}^{r+1})}+\|\v\|^{r-1}_{\mathrm{L}^{r+1}(0,T;\widetilde{\L}^{r+1})}\right)\nonumber\\&\qquad\qquad\times\|\v_m-\v\|^{\frac{1}{r-1}}_{\mathrm{L}^2(0,T;\L^2(\mathcal{O}_k))}\|\v_m-\v\|^{\frac{r-2}{r-1}}_{\mathrm{L}^{r+1}(0,T;\wi\L^{r+1})}\nonumber\\&
			\to 0  \ \text{ as }\   m\to \infty,
		\end{align*}
		for $r\geq 3$,	which implies that, for given $\varepsilon>0$, there exists $M(\varepsilon)\in\N$ such that
		\begin{align}\label{Conver15}
			\left|\int_{0}^{T}\uprho(t) \left[\int_{\mathcal{O}_k}\left(\left|\v_m(x,t)\right|^{r-1}\v_m(x,t)-\left|\v(x,t)\right|^{r-1}\v(x,t)\right)\textbf{g}(x)  \d x\right] \d t\right|\leq\frac{\varepsilon}{2}.
		\end{align}
		Finally, from \eqref{Conver13}-\eqref{Conver15}, we infer that, for a given $\varepsilon>0$, there exists $M(\varepsilon)\in\N$ such that
		\begin{align*}
			\left|\int_{0}^{T} \big\langle\mathcal{C}(\v_m(t)) ,\uprho(t)\textbf{g}  \big\rangle \d t- \int_{0}^{T} \big\langle\mathcal{C}(\v(t)) ,\uprho(t)\textbf{g}  \big\rangle \d t\right|\leq\varepsilon,
		\end{align*}
		for all $m\geq M(\varepsilon),$ which completes the proof.
	\end{proof}
	
	\begin{lemma}\label{weak_add}
		For $d=2$ with $r>1$, $d=3$ with $r>3$ and $d=r=3$ with $2\beta\mu\geq1$, assume that $\f\in\mathrm{L}^2_{\emph{loc}}(\R;\V')$. Let $\mathfrak{s}\in\R, \omega\in \Omega$ and $\v_{\mathfrak{s}}^0, \v_{\mathfrak{s}}^n\in \H$ for all $n\in\N.$ If $\v_{\mathfrak{s}}^n\xrightharpoonup{w}\v_{\mathfrak{s}}^0$ in $\H$, then the solution $\v$ of the system \eqref{CSCBF_Add} satisfies the following convergences:
		\begin{itemize}
			\item [(i)] $\v(\xi,\mathfrak{s},\omega,\v_{\mathfrak{s}}^n)\xrightharpoonup{w}\v(\xi,\mathfrak{s},\omega,\v_{\mathfrak{s}}^0)$ in $\H$ for all $\xi\geq \mathfrak{s}$.
			\item [(ii)] $\v(\cdot,\mathfrak{s},\omega,\v_{\mathfrak{s}}^n)\xrightharpoonup{w}\v(\cdot,\mathfrak{s},\omega,\v_{\mathfrak{s}}^0)$ in $\mathrm{L}^2((\mathfrak{s},\mathfrak{s}+T);\V)$ for every $T>0$.
			\item [(iii)] $\v(\cdot,\mathfrak{s},\omega,\v_{\mathfrak{s}}^n)\xrightharpoonup{w}\v(\cdot,\mathfrak{s},\omega,\v_{\mathfrak{s}}^0)$ in $\mathrm{L}^{r+1}((\mathfrak{s},\mathfrak{s}+T);\widetilde{\L}^{r+1})$ for every $T>0$.
			\item [(iv)] $\v(\cdot,\mathfrak{s},\omega,\v_{\mathfrak{s}}^n)\to\v(\cdot,\mathfrak{s},\omega,\v_{\mathfrak{s}}^0)$ in $\mathrm{L}^2((\mathfrak{s},\mathfrak{s}+T);\L^2(\mathcal{O}_k))$ for every $T>0$ and $k>0$,
			where $\mathcal{O}_k=\{x\in\R^d:|x|<k\}$.
		\end{itemize}
	\end{lemma}
	\begin{proof}
		Using the standard method as in \cite{KM1} (see Lemmas 5.2 and 5.3 in \cite{KM1}), one can complete the proof. Note that, due to compactness of Sobolev embedding in bounded domains, we get the  strong convergence in $\mathrm{L}^2((\mathfrak{s},\mathfrak{s}+T);\L^2(\mathcal{O}_k))$.
	\end{proof}
	\begin{lemma}\label{Asymptotic_UB_Add}
		For $d=2$ with $r>1$, $d=3$ with $r>3$ and $d=r=3$ with $2\beta\mu\geq1$, assume that $\f\in\mathrm{L}^2_{\emph{loc}}(\R;\V')$ satisfies \eqref{forcing6}. Then for every $\mathfrak{s}\in \R,$ $\omega\in \Omega,$ $D=\{D(\mathfrak{s},\omega):\mathfrak{s}\in \R,\omega\in \Omega\}\in \mathfrak{D}$ and $t_n\to \infty,$ $\u_{0,n}\in D(\mathfrak{s}-t_n, \vartheta_{-t_{n}}\omega)$, the sequence $\Phi_0(t_n,\mathfrak{s}-t_n,\vartheta_{-t_n}\omega,\u_{0,n})$ or $\u(\mathfrak{s},\mathfrak{s}-t_n,\vartheta_{-\mathfrak{s}}\omega,\u_{0,n})$ of solutions of the system \eqref{SCBF_Add} has a convergent subsequence in $\H$.
	\end{lemma}
	\begin{proof}
		We infer from \eqref{ue24} (for $d=2$) and \eqref{ue29} (for $d=3$) with $\tau=\mathfrak{s}$ that there exists $\mathfrak{T}=\mathfrak{T}(\mathfrak{s},\omega,D)>0$ such that for all $t\geq \mathfrak{T}$, $\u(\mathfrak{s},\mathfrak{s}-t,\vartheta_{-\mathfrak{s}}\omega,\u_{\mathfrak{s}-t})\in\H$ with $\u_{\mathfrak{s}-t}\in D(\mathfrak{s}-t,\vartheta_{-t}\omega).$ Since $t_n\to \infty$, there exists $N_5\in\N$ such that $t_n\geq \mathfrak{T}$ for all $n\geq N_5$. It is given that $\u_{0,n}\in D(\mathfrak{s}-t_n, \vartheta_{-t_{n}}\omega)$, we get  $\{\u(\mathfrak{s},\mathfrak{s}-t_n,\vartheta_{-\mathfrak{s}}\omega,\u_{0,n})\}_{n\geq N_5}$ is a bounded sequence in $\H$. Due to
		\begin{align}\label{Trans}
			\u(\mathfrak{s},\mathfrak{s}-t_n,\vartheta_{-\mathfrak{s}}\omega,\u_{0,n})=\v(\mathfrak{s},\mathfrak{s}-t_n,\vartheta_{-\mathfrak{s}}\omega,\v_{0,n})	+e^{\sigma\mathfrak{s}}\y(\omega)\textbf{g},
		\end{align}
		with $\v_{0,n}=\u_{0,n}-e^{\sigma(\mathfrak{s}-t)}\y(\vartheta_{-t}\omega)\textbf{g}$, we have that $\{\v(\mathfrak{s},\mathfrak{s}-t_n,\vartheta_{-\mathfrak{s}}\omega,\v_{0,n})\}_{n\geq N_5}$ is also a bounded sequence in $\H$, which implies that there exists $\tilde{\v}\in \H$ and a subsequence (not relabeling) such that 
		\begin{align}\label{ac2add}
			\v(\mathfrak{s},\mathfrak{s}-t_n,\vartheta_{-\mathfrak{s}}\omega,\v_{0,n})\xrightharpoonup{w}\tilde{\v} \ \text{ in }\  \H,
		\end{align}
		which gives
		\begin{align}\label{ac3add}
			\|\tilde{\v}\|_{\H}\leq\liminf_{n\to\infty}\|\v(\mathfrak{s},\mathfrak{s}-t_n,\vartheta_{-\mathfrak{s}}\omega,\v_{0,n})\|_{\H}.
		\end{align}
		For our purpose, we need to show that $\v(\mathfrak{s},\mathfrak{s}-t_n,\vartheta_{-\mathfrak{s}}\omega,\v_{0,n})\to\tilde{\v}$ in $\H$ strongly. Along with the above expression, it is enough  to prove that
		\begin{align}\label{ac4add}
			\|\tilde{\v}\|_{\H}\geq\limsup_{n\to\infty}\|\v(\mathfrak{s},\mathfrak{s}-t_n,\vartheta_{-\mathfrak{s}}\omega,\v_{0,n})\|_{\H}.
		\end{align}
		Now, for a given $j\in \N$ ($j\leq t_n$), we have 
		\begin{align}\label{ac5add}
			\v(\mathfrak{s},\mathfrak{s}-t_n,\vartheta_{-\mathfrak{s}}\omega,\v_{0,n})=\v(\mathfrak{s},\mathfrak{s}-j,\vartheta_{-\mathfrak{s}}\omega,\v(\mathfrak{s}-j,\mathfrak{s}-t_n,\vartheta_{-\mathfrak{s}}\omega,\v_{0,n})).
		\end{align} 
		For each $j$, let $N_j$ be sufficiently large such that $t_n\geq \mathfrak{T}+j$ for all $n\geq N_j$. From \eqref{ue23} (for $d=2$) and \eqref{ue28} (for $d=3$) with $\tau=\mathfrak{s}-j$, we have that the sequence $\{\v(\mathfrak{s}-j,\mathfrak{s}-t_n,\vartheta_{-\mathfrak{s}}\omega,\v_{0,n})\}_{n\geq N_j}$ is bounded in $\H$, for each $j\in \N$. By the diagonal process, there exists a subsequence (denoting by same label) and $\tilde{\u}_{j}\in \H$ for each $j\in\N$ such that 
		\begin{align}\label{ac6add}
			\v(\mathfrak{s}-j,\mathfrak{s}-t_n,\vartheta_{-\mathfrak{s}}\omega,\v_{0,n})\xrightharpoonup{w}\tilde{\v}_{j}\  \text{ in } \ \H.
		\end{align}
		From \eqref{ac5add}-\eqref{ac6add} together with Lemma \ref{weak_add}, we have that for $j\in\N$,
		\begin{align}\label{ac7add}
			\v(\mathfrak{s},\mathfrak{s}-t_n,\vartheta_{-\mathfrak{s}}\omega,\v_{0,n})\xrightharpoonup{w}\v(\mathfrak{s},\mathfrak{s}-j,\vartheta_{-\mathfrak{s}}\omega,\tilde{\v}_{j}) \text{ in } \H,
		\end{align}
		\begin{align}\label{ac8add}
			\v(\cdot,\mathfrak{s}-j,\vartheta_{-\mathfrak{s}}\omega,\v(\mathfrak{s}-j,\mathfrak{s}-t_n,\vartheta_{-\mathfrak{s}}\omega,\v_{0,n}))\xrightharpoonup{w}\v(\cdot,\mathfrak{s}-l,\vartheta_{-\mathfrak{s}}\omega,\tilde{\v}_{j}) \text{ in } \mathrm{L}^2((\mathfrak{s}-j,\mathfrak{s});\V),
		\end{align}
		\begin{align}\label{ac8'add}
			\v(\cdot,\mathfrak{s}-j,\vartheta_{-\mathfrak{s}}\omega,\v(\mathfrak{s}-j,\mathfrak{s}-t_n,\vartheta_{-\mathfrak{s}}\omega,\v_{0,n}))\xrightharpoonup{w}\v(\cdot,\mathfrak{s}-j,\vartheta_{-\mathfrak{s}}\omega,\tilde{\v}_{j}) \text{ in } \mathrm{L}^{r+1}((\mathfrak{s}-j,\mathfrak{s});\widetilde{\L}^{r+1}),
		\end{align}
		and
		\begin{align}\label{ac1add}
			\v(\cdot,\mathfrak{s}-j,\vartheta_{-\mathfrak{s}}\omega,\v(\mathfrak{s}-j,\mathfrak{s}-t_n,\vartheta_{-\mathfrak{s}}\omega,\v_{0,n}))\to\v(\cdot,\mathfrak{s}-l,\vartheta_{-\mathfrak{s}}\omega,\tilde{\v}_{j}) \text{ in } \mathrm{L}^2((\mathfrak{s}-j,\mathfrak{s});\L^2(\mathcal{O}_k)),
		\end{align}
		where $\mathcal{O}_k=\{x\in\R^d:|x|<k\}$. Clearly, \eqref{ac2add} and \eqref{ac7add} imply that 
		\begin{align}\label{ac9add}
			\v(\mathfrak{s},\mathfrak{s}-j,\vartheta_{-\mathfrak{s}}\omega,\tilde{\v}_{j})=\tilde{\v}.
		\end{align}
		From \eqref{ue17} we have
		\begin{align}\label{ac10add}
			&	\frac{\d}{\d t} \|\v\|^2_{\H} + \alpha\|\v\|^2_{\H}\nonumber\\&=-2\mu\|\nabla\v\|^2_{\H}-\alpha\|\v\|^2_{\H} - 2\beta\|\v+e^{\sigma t}\textbf{g}\y\|^{r+1}_{\wi \L^{r+1}} +2\beta\left\langle\mathcal{C}(\v+e^{\sigma t}\textbf{g}\y),e^{\sigma t}\textbf{g}\y\right\rangle\nonumber\\&\quad+2b(\v+e^{\sigma t}\textbf{g}\y,\v+e^{\sigma t}\textbf{g}\y,e^{\sigma t}\textbf{g}\y)+2\left\langle\f,\v\right\rangle+2e^{\sigma t}\y\left(\left(\ell-\sigma-\alpha\right)\textbf{g}-\mu \A\textbf{g},\v\right).
		\end{align}
		Now, following the similar arguments (that is, the method of energy equation) as in Lemma \ref{Asymptotic_UB} (or see Lemma \ref{precompact}) together with the convergence in \eqref{ac7add}-\eqref{ac1add}, Corollary \ref{convergence_b} and Corollary \ref{convergence_c}, we obtain that $\v(\mathfrak{s},\mathfrak{s}-t_n,\vartheta_{-\mathfrak{s}}\omega,\v_{0,n})\to\tilde{\v}$ strongly in $\H$. Hence, from \eqref{Trans}, we get that $\u(\mathfrak{s},\mathfrak{s}-t_n,\vartheta_{-\mathfrak{s}}\omega,\u_{0,n})$ has a convergent subsequence in $\H$, as required.
	\end{proof}
	The following theorem demonstrates the existence of a unique random $\mathfrak{D}$-pullback attractor for the system \eqref{SCBF_Add}.
	\begin{theorem}\label{RA_add}
		For $d=2$ with $r>1$, $d=3$ with $r>3$ and $d=r=3$ with $2\beta\mu\geq1$, assume that $\f\in\mathrm{L}^2_{\emph{loc}}(\R;\V')$ satisfies \eqref{forcing7}. Then there exists a unique random $\mathfrak{D}$-pullback attractor $$\mathscr{A}_0=\{\mathscr{A}_0(\mathfrak{s},\omega):\mathfrak{s}\in\R, \omega\in\Omega\}\in\mathfrak{D},$$ for the  the continuous cocycle $\Phi_0$ associated with the system \eqref{SCBF_Add}, in $\H$.
	\end{theorem}
	\begin{proof}
		The proof follows from Lemma \ref{Asymptotic_UB}, Lemma \ref{PAS'} and the abstract theory available in \cite{SandN_Wang} (Theorem 2.23 in \cite{SandN_Wang}).
	\end{proof}
	\subsection{Upper semicontinuity of random pullback attractors for SCBF equations with additive colored noise}
	Consider the random equation driven by colored noise 
	\begin{align}\label{z1}
		\frac{\d\z_{\delta}}{\d t}=-\ell \z_{\delta} + \mathcal{Z}_{\delta}(\vartheta_t\omega).
	\end{align}
	From \eqref{N3}, it is clear that for all $\omega\in\Omega$, the integral 
	\begin{align}\label{z2}
		\z_{\delta}(\omega)=\int_{-\infty}^{0}e^{\ell \tau}\mathcal{Z}_{\delta}(\vartheta_{\tau}\omega)\d\tau
	\end{align}
	is finite, and hence $\z_{\delta}:\Omega\to\R$ is a well defined random variable. Let us recall some properties of $\z_{\delta}$ from \cite{GGW}.
	\begin{lemma}[Lemma 3.2, \cite{GGW}]
		Let $\z_{\delta}$ be the random variable given by \eqref{z2}. Then the mapping 
		\begin{align}\label{z3}
			(t,\omega)\mapsto\z_{\delta}(\vartheta_{t}\omega)=\int_{-\infty}^{t}e^{-\ell( t-\tau)}\mathcal{Z}_{\delta}(\vartheta_{\tau}\omega)\d\tau
		\end{align}
		is a stationary solution of \eqref{z1} with continuous paths. Moreover, $\mathbb{E}[\z_{\delta}]=0$ and for every $\omega\in\Omega$,
		\begin{align}\label{z4}
			\lim_{\delta\to0}\z_{\delta}(\vartheta_{t}\omega)=\y(\vartheta_{t}\omega) \text{ uniformly on } [\mathfrak{s},\mathfrak{s}+T]\text{ with } \mathfrak{s}\in\R\text{ and } T>0;
		\end{align}
		\begin{align}\label{z5}
			\lim_{t\to\pm\infty}\frac{\left|\z_{\delta}(\vartheta_{t}\omega)\right|}{|t|}=0 \text{ uniformly for }  0<\delta\leq\tilde{\ell};
		\end{align}
		where $\tilde{\ell}=\min\{1,\frac{1}{2\ell}\}$.
	\end{lemma}
	Define
	\begin{align}\label{WZ_T_add}
		\v_{\delta}(t,\mathfrak{s},\omega,\v_{\delta,\mathfrak{s}})=\u_{\delta}(t,\mathfrak{s},\omega,\u_{\delta,\mathfrak{s}})-e^{\sigma t}\textbf{g}(x)\z_{\delta}(\vartheta_{t}\omega).
	\end{align}
	Then, from \eqref{WZ_SCBF_Add}, we obtain
	\begin{equation}\label{WZ_CSCBF_Add}
		\left\{
		\begin{aligned}
			\frac{\partial \v_{\delta}}{\partial t}+\mu\A\v_{\delta}+\B(\v_{\delta}&+e^{\sigma t}\textbf{g}\z_{\delta})+\alpha\v_{\delta}+\beta\mathcal{C}(\v_{\delta}+e^{\sigma t}\textbf{g}\z_{\delta})\\&=\boldsymbol{f}+\left(\ell-\sigma-\alpha\right)e^{\sigma t}\textbf{g}\z_{\delta}-\mu e^{\sigma t}\z_{\delta}\A\textbf{g}, \ \ \ \text{ in } \mathbb{R}^d\times(\mathfrak{s},\infty), \\ 
			\v_{\delta}(x,\mathfrak{s})&=\v_{\delta,\mathfrak{s}}(x)=\u_{\delta,\mathfrak{s}}(x)-e^{\sigma\mathfrak{s}}\textbf{g}(x)\z_{\delta}(\omega), \quad\quad \ \  x\in \mathbb{R}^d \text{ and }\mathfrak{s}\in\R.
		\end{aligned}
		\right.
	\end{equation}
	For all $\mathfrak{s}\in\R,$ $t>\mathfrak{s},$ and for every $\v_{\delta,\mathfrak{s}}\in\H$ and $\omega\in\Omega$, \eqref{WZ_CSCBF_Add} has a unique solution $\v_{\delta}(\cdot,\mathfrak{s},\omega,\v_{\delta,\mathfrak{s}})\in \mathrm{C}([\mathfrak{s},\mathfrak{s}+T];\H)\cap\mathrm{L}^2(\mathfrak{s}, \mathfrak{s}+T;\V)\cap\mathrm{L}^{r+1}(\mathfrak{s},\mathfrak{s}+T;\widetilde{\L}^{r+1})$. Moreover, $\v_{\delta}(t,\mathfrak{s},\omega,\v_{\delta,\mathfrak{s}})$ is continuous with respect to the initial data $\v_{\delta,\mathfrak{s}}(x)$ and $(\mathscr{F},\mathscr{B}(\H))$-measurable in $\omega\in\Omega.$ Define a cocycle $\Phi_{\delta}:\R^+\times\R\times\Omega\times\H\to\H$ for the system \eqref{WZ_SCBF_Add} such that for given $t\in\R^+, \mathfrak{s}\in\R, \omega\in\Omega$ and $\u_{\delta,\mathfrak{s}}\in\H$,
	\begin{align}\label{Phi_d}
		\Phi_{\delta}(t,\mathfrak{s},\omega,\u_{\delta,\mathfrak{s}}) &=\u_{\delta}(t+\mathfrak{s},\mathfrak{s},\vartheta_{-\mathfrak{s}}\omega,\u_{\delta,\mathfrak{s}})=\v_{\delta}(t+\mathfrak{s},\mathfrak{s},\vartheta_{-\mathfrak{s}}\omega,\u_{\delta,\mathfrak{s}})+e^{\sigma(t+\mathfrak{s})}\textbf{g}\z_{\delta}(\vartheta_{t}\omega).
	\end{align}
	\begin{lemma}\label{LemmaUe_add1}
		For $d=2$ with $r>1$, $d=3$ with $r>3$ and $d=r=3$ with $2\beta\mu\geq1$, assume that $\f\in \mathrm{L}^2_{\emph{loc}}(\mathbb{R};\V')$ satisfies \eqref{forcing6}. Then $\Phi_\delta$ possesses a closed measurable $\mathfrak{D}$-pullback absorbing set $\mathcal{K}^1_\delta=\{\mathcal{K}^1_\delta(\mathfrak{s},\omega):\mathfrak{s}\in\R, \omega\in\Omega\}\in\mathfrak{D}$ in $\H$ given by
		\begin{align}\label{ue_add4}
			\mathcal{K}^1_\delta(\mathfrak{s},\omega)=\{\u\in\H:\|\u\|^2_{\H}\leq \mathcal{R}^1_\delta(\mathfrak{s},\omega)\},\ \  \text{ for } d=2, 
		\end{align}
		where $\mathcal{R}^1_\delta(\mathfrak{s},\omega)$ is defined by
		\begin{align}\label{ue_add5}
			\mathcal{R}^1_\delta(\mathfrak{s},\omega)&=3\|\textbf{g}\|^2_{\H}\left|e^{\sigma\mathfrak{s}}\z_{\delta}(\omega)\right|^2+2R_{5} \int_{-\infty}^{0}e^{\alpha\xi}\bigg[\|\f(\cdot,\xi+\mathfrak{s})\|^2_{\V'}+\left|e^{\sigma (\xi+\mathfrak{s})}\z_{\delta}(\vartheta_{\xi}\omega)\right|^2\nonumber\\&\quad+\left|e^{\sigma (\xi+\mathfrak{s})}\z_{\delta}(\vartheta_{\xi}\omega)\right|^{r+1}+\left|e^{\sigma(\xi+\mathfrak{s})}\z_{\delta}(\vartheta_{\xi}\omega)\right|^{\frac{2(r+1)}{r-1}}\bigg]\d\xi,
		\end{align}
		and a closed measurable $\mathfrak{D}$-pullback absorbing set $\mathcal{K}^2_\delta=\{\mathcal{K}^2_\delta(\mathfrak{s},\omega):\mathfrak{s}\in\R, \omega\in\Omega\}\in\mathfrak{D}$ in $\H$ given by
		\begin{align}\label{ue_add6}
			\mathcal{K}^2_\delta(\mathfrak{s},\omega)=\{\u\in\H:\|\u\|^2_{\H}\leq \mathcal{R}^2_\delta(\mathfrak{s},\omega)\},\ \  \text{ for } d=3, 
		\end{align}
		where $\mathcal{R}^2_\delta(\mathfrak{s},\omega)$ is defined by
		\begin{align}\label{ue_add7}
			\mathcal{R}^2_\delta(\mathfrak{s},\omega)&=3\|\textbf{g}\|^2_{\H}\left|e^{\sigma\mathfrak{s}}\z_{\delta}(\omega)\right|^2+2R_{8} \int_{-\infty}^{0}e^{\alpha\xi}\bigg[\|\f(\cdot,\xi+\mathfrak{s})\|^2_{\V'}+\left|e^{\sigma (\xi+\mathfrak{s})}\z_{\delta}(\vartheta_{\xi}\omega)\right|^2\nonumber\\&\quad+\left|e^{\sigma (\xi+\mathfrak{s})}\z_{\delta}(\vartheta_{\xi}\omega)\right|^{r+1}\bigg]\d\xi.
		\end{align}
		Here $R_{5}$ and $R_{8}$ are same as in Lemma \ref{LemmaUe_add}. Furthermore, for every $\mathfrak{s}\in\R$, $\omega\in\Omega$ and $i\in\{1,2\}$,
		\begin{align}\label{ue_add8}
			\lim_{\delta\to0}\mathcal{R}^{i}_\delta(\mathfrak{s},\omega)=\mathcal{R}^{i}_0(\mathfrak{s},\omega),
		\end{align}
		where $\mathcal{R}^{1}_0(\mathfrak{s},\omega)$ and $\mathcal{R}^{2}_0(\mathfrak{s},\omega)$ are given by \eqref{ue_add1} and \eqref{ue_add3}, respectively.
	\end{lemma}
	\begin{proof}
		Since, the systems \eqref{CSCBF_Add} and \eqref{WZ_CSCBF_Add} have similar terms, the existence of closed measurable $\mathfrak{D}$-pullback absorbing sets $\mathcal{K}^{i}_\delta(\mathfrak{s},\omega)$ (for every $i\in\{1,2\}$) is confirmed from the similar calculations as in Lemma  \ref{LemmaUe_add}. Now, it is  left to prove \eqref{ue_add8} only. 
		
		Using \eqref{z5}, we can find $\xi_0<0$ such that for all $0<\delta\leq\tilde{\ell}$,
		\begin{align}\label{ue30}
			\left|\z_{\delta}(\vartheta_{\xi}\omega)\right|\leq|\xi|, \ \ \ \text{ for all } \ \xi\leq\xi_0.
		\end{align}
		From \eqref{ue30}, we get that for all $\xi\leq\xi_0$ and $0<\delta\leq\tilde{\ell}$,
		\begin{align}\label{ue31}
			&	e^{\alpha\xi}\left[\left|e^{\sigma (\xi+\mathfrak{s})}\z_{\delta}(\vartheta_{\xi}\omega)\right|^2+\left|e^{\sigma (\xi+\mathfrak{s})}\z_{\delta}(\vartheta_{\xi}\omega)\right|^{r+1}+\left|e^{\sigma(\xi+\mathfrak{s})}\z_{\delta}(\vartheta_{\xi}\omega)\right|^{\frac{2(r+1)}{r-1}}\right]\nonumber\\&\leq e^{\alpha\xi}\left[\left|e^{\sigma (\xi+\mathfrak{s})}\xi\right|^2+\left|e^{\sigma (\xi+\mathfrak{s})}\xi\right|^{r+1}+\left|e^{\sigma(\xi+\mathfrak{s})}\xi\right|^{\frac{2(r+1)}{r-1}}\right],
		\end{align}
		and 
		\begin{align}\label{ue32}
			\int_{-\infty}^{0}e^{\alpha\xi}\left[\left|e^{\sigma (\xi+\mathfrak{s})}\xi\right|^2+\left|e^{\sigma (\xi+\mathfrak{s})}\xi\right|^{r+1}+\left|e^{\sigma(\xi+\mathfrak{s})}\xi\right|^{\frac{2(r+1)}{r-1}}\right]\d\xi <\infty.
		\end{align}
		Consider,
		\begin{align}\label{ue33}
			&\int_{-\infty}^{0}e^{\alpha\xi}\left[\left|e^{\sigma (\xi+\mathfrak{s})}\z_{\delta}(\vartheta_{\xi}\omega)\right|^2+\left|e^{\sigma (\xi+\mathfrak{s})}\z_{\delta}(\vartheta_{\xi}\omega)\right|^{r+1}+\left|e^{\sigma(\xi+\mathfrak{s})}\z_{\delta}(\vartheta_{\xi}\omega)\right|^{\frac{2(r+1)}{r-1}}\right]\d\xi\nonumber\\&=\int_{\xi_0}^{0}e^{\alpha\xi}\left[\left|e^{\sigma (\xi+\mathfrak{s})}\z_{\delta}(\vartheta_{\xi}\omega)\right|^2+\left|e^{\sigma (\xi+\mathfrak{s})}\z_{\delta}(\vartheta_{\xi}\omega)\right|^{r+1}+\left|e^{\sigma(\xi+\mathfrak{s})}\z_{\delta}(\vartheta_{\xi}\omega)\right|^{\frac{2(r+1)}{r-1}}\right]\d\xi\nonumber\\&\quad+\int_{-\infty}^{\xi_0}e^{\alpha\xi}\left[\left|e^{\sigma (\xi+\mathfrak{s})}\z_{\delta}(\vartheta_{\xi}\omega)\right|^2+\left|e^{\sigma (\xi+\mathfrak{s})}\z_{\delta}(\vartheta_{\xi}\omega)\right|^{r+1}+\left|e^{\sigma(\xi+\mathfrak{s})}\z_{\delta}(\vartheta_{\xi}\omega)\right|^{\frac{2(r+1)}{r-1}}\right]\d\xi.
		\end{align}
		Using \eqref{z4}, \eqref{ue31}, \eqref{ue32} and the Lebesgue Dominated Convergence Theorem, we get 
		\begin{align}\label{ue34}
			&\lim_{\delta\to0}\int_{-\infty}^{\xi_0}e^{\alpha\xi}\left[\left|e^{\sigma (\xi+\mathfrak{s})}\z_{\delta}(\vartheta_{\xi}\omega)\right|^2+\left|e^{\sigma (\xi+\mathfrak{s})}\z_{\delta}(\vartheta_{\xi}\omega)\right|^{r+1}+\left|e^{\sigma(\xi+\mathfrak{s})}\z_{\delta}(\vartheta_{\xi}\omega)\right|^{\frac{2(r+1)}{r-1}}\right]\d\xi\nonumber\\&=\int_{-\infty}^{\xi_0}e^{\alpha\xi}\left[\left|e^{\sigma (\xi+\mathfrak{s})}\y(\vartheta_{\xi}\omega)\right|^2+\left|e^{\sigma (\xi+\mathfrak{s})}\y(\vartheta_{\xi}\omega)\right|^{r+1}+\left|e^{\sigma(\xi+\mathfrak{s})}\y(\vartheta_{\xi}\omega)\right|^{\frac{2(r+1)}{r-1}}\right]\d\xi.
		\end{align}
		From \eqref{z4}, we obtain 
		\begin{align}\label{ue35}
			&\lim_{\delta\to0}\int_{\xi_0}^{0}e^{\alpha\xi}\left[\left|e^{\sigma (\xi+\mathfrak{s})}\z_{\delta}(\vartheta_{\xi}\omega)\right|^2+\left|e^{\sigma (\xi+\mathfrak{s})}\z_{\delta}(\vartheta_{\xi}\omega)\right|^{r+1}+\left|e^{\sigma(\xi+\mathfrak{s})}\z_{\delta}(\vartheta_{\xi}\omega)\right|^{\frac{2(r+1)}{r-1}}\right]\d\xi\nonumber\\&=\int_{\xi_0}^{0}e^{\alpha\xi}\left[\left|e^{\sigma (\xi+\mathfrak{s})}\y(\vartheta_{\xi}\omega)\right|^2+\left|e^{\sigma (\xi+\mathfrak{s})}\y(\vartheta_{\xi}\omega)\right|^{r+1}+\left|e^{\sigma(\xi+\mathfrak{s})}\y(\vartheta_{\xi}\omega)\right|^{\frac{2(r+1)}{r-1}}\right]\d\xi.
		\end{align}
		Hence, from \eqref{z4}, \eqref{ue34} and \eqref{ue35}, we deduce  the desired convergence.
	\end{proof}
	Next result demonstrates the convergence of solution of \eqref{WZ_SCBF_Add} to the solution of \eqref{SCBF_Add} in $\H$ as $\delta\to 0$.
	\begin{lemma}\label{Solu_Conver}
		For $d=2$ with $r>1$, $d=3$ with $r>3$ and $d=r=3$ with $2\beta\mu\geq1$, assume that $\f\in \mathrm{L}^2_{\emph{loc}}(\mathbb{R};\V')$. Suppose that $\{\delta_n\}_{n\in\N}$ be the sequence such that $\delta_n\to0$. Let $\u_{\delta_n}$ and $\u$ be the solutions of \eqref{WZ_SCBF_Add} and \eqref{SCBF_Add} with initial data $\u_{\delta_n,\mathfrak{s}}$ and $\u_{\mathfrak{s}}$, respectively. If $\|\u_{\delta_n,\mathfrak{s}}-\u_{\mathfrak{s}}\|_{\H}\to0$ as $n\to\infty$, then for every $\mathfrak{s}\in\R$, $\omega\in\Omega$ and $t>\mathfrak{s}$,
		\begin{align*}
			\|\u_{\delta_n}(t,\mathfrak{s},\omega,\u_{\delta_n,\mathfrak{s}})-\u(t,\mathfrak{s},\omega,\u_{\mathfrak{s}})\|_{\H} \to0 \  \ \text{ as }\  \  n\to\infty.
		\end{align*}
	\end{lemma}
	\begin{proof}
		Let $\mathfrak{F}=\v_{\delta_n}-\v$, where $\v_{\delta_n}$ and $\v$ are the solutions of \eqref{WZ_CSCBF_Add} and \eqref{CSCBF_Add}, respectively. Also, let $\mathfrak{z}=\z_{\delta_n}-\y$. Then we get from \eqref{WZ_CSCBF_Add} and \eqref{CSCBF_Add} that
		\begin{align}\label{F1}
			\frac{1}{2}\frac{\d}{\d t}\|\mathfrak{F}\|^2_{\H}&=-\mu\|\nabla\mathfrak{F}\|^2_{\H}-\alpha\|\mathfrak{F}\|^2_{\H} -\left\langle \B(\v_{\delta_n}+e^{\sigma t}\textbf{g}\z_{\delta_n})-\B(\v+e^{\sigma t}\textbf{g}\y), \mathfrak{F}\right\rangle\nonumber\\&\quad-\beta\left\langle\mathcal{C}(\v_{\delta_n}+e^{\sigma t}\textbf{g}\z_{\delta_n})-\mathcal{C}(\v+e^{\sigma t}\textbf{g}\y),\mathfrak{F}\right\rangle+e^{\sigma t}\mathfrak{z}\left(\left(\ell-\sigma-\alpha\right)\textbf{g}-\mu \A\textbf{g},\mathfrak{F}\right) \nonumber\\&=-\mu\|\nabla\mathfrak{F}\|^2_{\H}-\alpha\|\mathfrak{F}\|^2_{\H}+e^{\sigma t}\mathfrak{z}\left(\left(\ell-\sigma-\alpha\right)\textbf{g}-\mu \A\textbf{g},\mathfrak{F}\right)\nonumber\\&\quad-\left\langle \B(\v_{\delta_n}+e^{\sigma t}\textbf{g}\z_{\delta_n})-\B(\v+e^{\sigma t}\textbf{g}\y), (\v_{\delta_n}+e^{\sigma t}\textbf{g}\z_{\delta_n})-(\v+e^{\sigma t}\textbf{g}\y)\right\rangle \nonumber\\&\quad-\beta\left\langle\mathcal{C}(\v_{\delta_n}+e^{\sigma t}\textbf{g}\z_{\delta_n})-\mathcal{C}(\v+e^{\sigma t}\textbf{g}\y),(\v_{\delta_n}+e^{\sigma t}\textbf{g}\z_{\delta_n})-(\v+e^{\sigma t}\textbf{g}\y)\right\rangle \nonumber\\&\quad+e^{\sigma t}\mathfrak{z}b(\v_{\delta_n},\v_{\delta_n},\textbf{g})+e^{2\sigma t}\mathfrak{z}\z_{\delta_n}b(\textbf{g},\v_{\delta_n},\textbf{g})-e^{\sigma t}\mathfrak{z}b(\v,\v,\textbf{g})-e^{2\sigma t}\mathfrak{z}\y b(\textbf{g},\v,\textbf{g})\nonumber\\&\quad+\beta e^{\sigma t}\mathfrak{z}\left\langle\mathcal{C}(\v_{\delta_n}+e^{\sigma t}\textbf{g}\z_{\delta_n})-\mathcal{C}(\v+e^{\sigma t}\textbf{g}\y),\textbf{g}\right\rangle.
		\end{align}
		From \eqref{MO_c}, we write
		\begin{align}\label{F2}
			&	-\beta\left\langle\mathcal{C}(\v_{\delta_n}+e^{\sigma t}\textbf{g}\z_{\delta_n})-\mathcal{C}(\v+e^{\sigma t}\textbf{g}\y),(\v_{\delta_n}+e^{\sigma t}\textbf{g}\z_{\delta_n})-(\v+e^{\sigma t}\textbf{g}\y)\right\rangle\nonumber\\&\leq-\frac{\beta}{2}\|\left|\v+e^{\sigma t}\textbf{g}\y\right|^{\frac{r-1}{2}}\left|\mathfrak{F}+e^{\sigma t}\textbf{g}\mathfrak{z}\right|\|^2_{\H}\leq0. 
		\end{align}
		Applying Lemmas \ref{Holder} and \ref{Young}, we obtain
		\begin{align}\label{F3}
			&	\left|\beta e^{\sigma t}\mathfrak{z}\left\langle\mathcal{C}(\v_{\delta_n}+e^{\sigma t}\textbf{g}\z_{\delta_n})-\mathcal{C}(\v+e^{\sigma t}\textbf{g}\y),\textbf{g}\right\rangle\right|\nonumber\\&\leq\beta e^{\sigma t}\left|\mathfrak{z}\right|\left(\|\v_{\delta_n}+e^{\sigma t}\textbf{g}\z_{\delta_n}\|^r_{\wi\L^{r+1}}+\|\v+e^{\sigma t}\textbf{g}\y\|^r_{\wi\L^{r+1}}\right)\|\textbf{g}\|_{\wi\L^{r+1}}\nonumber\\&\leq C e^{\sigma t}\left|\mathfrak{z}\right|\left(\|\v_{\delta_n}\|^r_{\wi\L^{r+1}}+\left|e^{\sigma t}\z_{\delta_n}\right|^r+\|\v\|^r_{\wi\L^{r+1}}+\left|e^{\sigma t}\y\right|^r\right)\nonumber\\&\leq C e^{\sigma t}\left|\mathfrak{z}\right|\left(1+\left|e^{\sigma t}\z_{\delta_n}\right|^r+\left|e^{\sigma t}\y\right|^r+\|\v_{\delta_n}\|^{r+1}_{\wi\L^{r+1}} +\|\v\|^{r+1}_{\wi\L^{r+1}}\right),
		\end{align}
		and
		\begin{align}\label{F4}
			\left|e^{\sigma t}\mathfrak{z}\left(\left(\ell-\sigma-\alpha\right)\textbf{g}-\mu \A\textbf{g},\mathfrak{F}\right)\right|\leq\frac{\alpha}{2}\|\mathfrak{F}\|^2_{\H} + Ce^{2\sigma t}\left|\mathfrak{z}\right|^2.
		\end{align}
		\vskip 2mm
		\noindent
		\textbf{Case I:} \textit{$d= 2$ and $r>1$.}
		Using \eqref{b1}, \eqref{441}, Lemmas \ref{Holder} and \ref{Young}, we get
		\begin{align}\label{F5}
			&	\left|\left\langle \B(\v_{\delta_n}+e^{\sigma t}\textbf{g}\z_{\delta_n})-\B(\v+e^{\sigma t}\textbf{g}\y), (\v_{\delta_n}+e^{\sigma t}\textbf{g}\z_{\delta_n})-(\v+e^{\sigma t}\textbf{g}\y)\right\rangle\right|\nonumber\\&\leq\left|b(\mathfrak{F}+e^{\sigma t}\textbf{g}\mathfrak{z},\mathfrak{F}+e^{\sigma t}\textbf{g}\mathfrak{z},\v+e^{\sigma t}\textbf{g}\y)\right|\nonumber\\&\leq C\|\mathfrak{F}+e^{\sigma t}\textbf{g}\mathfrak{z}\|_{\H}\|\nabla(\mathfrak{F}+e^{\sigma t}\textbf{g}\mathfrak{z})\|_{\H}\|\nabla(\v+e^{\sigma t}\textbf{g}\y)\|_{\H}\nonumber\\&\leq\frac{\mu}{2}\|\nabla\mathfrak{F}\|^2_{\H}+C\left(e^{2\sigma t}\left|\y\right|^2+\|\nabla\v\|^2_{\H}\right)\|\mathfrak{F}\|^2_{\H}+ Ce^{2\sigma t}\left|\mathfrak{z}\right|^2\left(1+e^{2\sigma t}\left|\y\right|^2+\|\nabla\v\|^2_{\H}\right).
		\end{align}
		Again, from \eqref{b1}, Lemmas \ref{Holder} and \ref{Young}, we have
		\begin{align}
			\left|e^{\sigma t}\mathfrak{z}b(\v_{\delta_n},\v_{\delta_n},\textbf{g})\right|&\leq Ce^{\sigma t}\left|\mathfrak{z}\right|\|\v_{\delta_n}\|_{\H}\|\nabla\v_{\delta_n}\|_{\H}\|\nabla\textbf{g}\|_{\H}\leq Ce^{\sigma t}\left|\mathfrak{z}\right|\left(1+\|\v_{\delta_n}\|^2_{\H}\|\nabla\v_{\delta_n}\|^2_{\H}\right),\label{F6}\\
			\left|e^{2\sigma t}\mathfrak{z}\z_{\delta_n}b(\textbf{g},\v_{\delta_n},\textbf{g})\right|&\leq Ce^{2\sigma t}\left|\mathfrak{z}\right|\left|\z_{\delta_n}\right|\|\textbf{g}\|_{\H}\|\nabla\textbf{g}\|_{\H}\|\nabla\v_{\delta_n}\|_{\H}\leq Ce^{2\sigma t}\left|\mathfrak{z}\right|\left|\z_{\delta_n}\right|\left(1+\|\nabla\v_{\delta_n}\|^2_{\H}\right),\label{F7}\\
			\left|e^{\sigma t}\mathfrak{z}b(\v,\v,\textbf{g})\right|&\leq Ce^{\sigma t}\left|\mathfrak{z}\right|\|\v\|_{\H}\|\nabla\v\|_{\H}\|\nabla\textbf{g}\|_{\H}\leq Ce^{\sigma t}\left|\mathfrak{z}\right|\left(1+\|\v\|^2_{\H}\|\nabla\v\|^2_{\H}\right),\label{F8}\\
			\left|e^{2\sigma t}\mathfrak{z}\y b(\textbf{g},\v,\textbf{g})\right|&\leq Ce^{2\sigma t}\left|\mathfrak{z}\right|\left|\y\right|\|\textbf{g}\|_{\H}\|\nabla\textbf{g}\|_{\H}\|\nabla\v\|_{\H}\leq Ce^{2\sigma t}\left|\mathfrak{z}\right|\left|\y\right|\left(1+\|\nabla\v\|^2_{\H}\right).\label{F9}
		\end{align}
		Combining \eqref{F1}-\eqref{F9}, we obtain
		\begin{align}\label{F10}
			\frac{\d}{\d t}\|\mathfrak{F}\|^2_{\H}\leq C\left[Q_1(t)\|\mathfrak{F}\|^2_{\H} + \left|\mathfrak{z}(\vartheta_{t}\omega)\right|Q_2(t)\right],
		\end{align}
		for a.e. $t\in[\mathfrak{s},\mathfrak{s}+T]$, where
		\begin{align*}
			Q_1(t)&=e^{2\sigma t}\left|\y(\vartheta_{t}\omega)\right|^2+\|\nabla\v(t)\|^2_{\H}\  \ \text{ and }\\
			Q_2(t)&=e^{2\sigma t}\left|\mathfrak{z}(\vartheta_{t}\omega)\right|\bigg[1+e^{2\sigma t}\left|\y(\vartheta_{t}\omega)\right|^2+\|\nabla\v(t)\|^2_{\H}\bigg]+e^{\sigma t}\bigg[1+\|\v_{\delta_n}(t)\|^2_{\H}\|\nabla\v_{\delta_n}(t)\|^2_{\H}\nonumber\\&\quad+\|\v(t)\|^2_{\H}\|\nabla\v(t)\|^2_{\H}+\left|e^{\sigma t}\z_{\delta_n}(\vartheta_{t}\omega)\right|^r+\left|e^{\sigma t}\y(\vartheta_{t}\omega)\right|^r+\|\v_{\delta_n}(t)\|^{r+1}_{\wi\L^{r+1}} +\|\v(t)\|^{r+1}_{\wi\L^{r+1}}\bigg]\nonumber\\&\quad+e^{2\sigma t}\left|\z_{\delta_n}(\vartheta_{t}\omega)\right|\bigg[1+\|\nabla\v_{\delta_n}(t)\|^2_{\H}\bigg]+e^{2\sigma t}\left|\y(\vartheta_{t}\omega)\right|\bigg[1+\|\nabla\v(t)\|^2_{\H}\bigg].
		\end{align*}
		Due to continuity of $e^{\sigma t},$ $\y(\vartheta_{t}\omega)$ and $\z_{\delta_n}(\vartheta_{t}\omega)$, and $\v_{\delta_n}, \v\in\mathrm{C}([\mathfrak{s},\mathfrak{s}+T];\H)\cap\mathrm{L}^2(\mathfrak{s}, \mathfrak{s}+T;\V)\cap\mathrm{L}^{r+1}(\mathfrak{s},\mathfrak{s}+T;\widetilde{\L}^{r+1})$, we infer that 
		\begin{align}\label{F11}
			\int_{\mathfrak{s}}^{\mathfrak{s}+T}Q_1(\xi)\d\xi<\infty \ \text{ and }\  \int_{\mathfrak{s}}^{\mathfrak{s}+T}Q_2(\xi)\d\xi<\infty.
		\end{align}
		An application of Gronwall's inequality yields
		\begin{align}\label{F12}
			&\|\v_{\delta_n}(t,\mathfrak{s},\omega,\v_{\delta_n,\mathfrak{s}})-\v(t,\mathfrak{s},\omega,\v_{\mathfrak{s}})\|^2_{\H}\nonumber\\&\leq\left[ \|\v_{\delta_n,\mathfrak{s}}-\v_{\mathfrak{s}}\|^2_{\H}+C\sup_{\xi  \in [\mathfrak{s},\mathfrak{s}+ T]}\left|\z_{\delta_n}(\vartheta_{\xi}\omega)-\y(\vartheta_{\xi}\omega)\right|\int_{\mathfrak{s}}^{\mathfrak{s}+T}Q_1(\xi)\d\xi\right]e^{C\int\limits_{\mathfrak{s}}^{\mathfrak{s}+T}Q_2(\xi)\d\xi}.
		\end{align}
		By \eqref{T_add}, \eqref{WZ_T_add} and \eqref{F12}, we get
		\begin{align}\label{F13}
			&\|\u_{\delta_n}(t,\mathfrak{s},\omega,\u_{\delta_n,\mathfrak{s}})-\u(t,\mathfrak{s},\omega,\u_{\mathfrak{s}})\|^2_{\H}\nonumber\\&\leq2\bigg[ 2\|\u_{\delta_n,\mathfrak{s}}-\u_{\mathfrak{s}}\|^2_{\H} +2e^{\sigma\mathfrak{s}}\left|\z_{\delta_n}(\omega)-\y(\omega)\right|\|\textbf{g}\|^2_{\H}\nonumber\\&\quad+C\sup_{\xi  \in [\mathfrak{s},\mathfrak{s}+ T]}\left|\z_{\delta_n}(\vartheta_{\xi}\omega)-\y(\vartheta_{\xi}\omega)\right|\int\limits_{\mathfrak{s}}^{\mathfrak{s}+T}Q_1(\xi)\d\xi\bigg]e^{C\int\limits_{\mathfrak{s}}^{\mathfrak{s}+T}Q_2(\xi)\d\xi} \nonumber\\&\quad+ 2e^{\sigma t}\left|\z_{\delta_n}(\vartheta_{t}\omega)-\y(\vartheta_{t}\omega)\right|\|\textbf{g}\|^2_{\H},
		\end{align}
		for all $t\in[\mathfrak{s},\mathfrak{s}+T]$. Hence, due to \eqref{z4} and \eqref{F11}, we achieve the required convergence from \eqref{F13}.
		\vskip 2mm
		\noindent
		\textbf{Case II:} \textit{$d= 3$ and $r\geq3$ ($r>3$ with any $\beta,\mu>0$ and $r=3$ with $2\beta\mu\geq1$).}
		Using \eqref{441}, Lemmas \ref{Holder} and \ref{Young}, we get
		\begin{align}\label{F14}
			&	\left|\left\langle \B(\v_{\delta_n}+e^{\sigma t}\textbf{g}\z_{\delta_n})-\B(\v+e^{\sigma t}\textbf{g}\y), (\v_{\delta_n}+e^{\sigma t}\textbf{g}\z_{\delta_n})-(\v+e^{\sigma t}\textbf{g}\y)\right\rangle\right|\nonumber\\&\leq\left|b(\mathfrak{F}+e^{\sigma t}\textbf{g}\mathfrak{z},\mathfrak{F}+e^{\sigma t}\textbf{g}\mathfrak{z},\v+e^{\sigma t}\textbf{g}\y)\right|\nonumber\\&\leq\begin{cases}
				\frac{1}{2\beta}\|\nabla\left(\mathfrak{F}+e^{\sigma t}\textbf{g}\mathfrak{z}\right)\|^2_{\H}+\frac{\beta}{2}\|\left|\v+e^{\sigma t}\textbf{g}\y\right|\left|\mathfrak{F}+e^{\sigma t}\textbf{g}\mathfrak{z}\right|\|^2_{\H}, \text{ for } r=3,\\ \frac{\mu}{4}\|\nabla\left(\mathfrak{F}+e^{\sigma t}\textbf{g}\mathfrak{z}\right)\|^2_{\H}+\frac{\beta}{2}\|\left|\v+e^{\sigma t}\textbf{g}\y\right|^{\frac{r-1}{2}}\left|\mathfrak{F}+e^{\sigma t}\textbf{g}\mathfrak{z}\right|\|^2_{\H}+C\|\mathfrak{F}+e^{\sigma t}\textbf{g}\mathfrak{z}\|^2_{\H}, \text{ for } r>3,	\end{cases}\nonumber\\&\leq\begin{cases}
				\frac{1}{2\beta}\|\nabla\mathfrak{F}\|^2_{\H}+C\left|e^{\sigma t}\mathfrak{z}\right|^2+C\left|e^{\sigma t}\mathfrak{z}\right|\left(1+\|\nabla\v_{\delta_n}\|^2_{\H}+\|\nabla\v\|^2_{\H}\right)+\frac{\beta}{2}\||\v+e^{\sigma t}\textbf{g}\y||\mathfrak{F}+e^{\sigma t}\textbf{g}\mathfrak{z}|\|^2_{\H},\\ \hspace{131mm} \text{ for } r=3,\\
				\frac{\mu}{2}\|\nabla\mathfrak{F}\|^2_{\H}+C\|\mathfrak{F}\|^2_{\H}+ Ce^{2\sigma t}\left|\mathfrak{z}\right|^2+\frac{\beta}{2}\|\left|\v+e^{\sigma t}\textbf{g}\y\right|^{\frac{r-1}{2}}\left|\mathfrak{F}+e^{\sigma t}\textbf{g}\mathfrak{z}\right|\|^2_{\H}, \text{ for } r>3. 
			\end{cases}
		\end{align}
		From \eqref{b1}, Lemmas \ref{Holder} and \ref{Young}, we have
		\begin{align}
			\left|e^{\sigma t}\mathfrak{z}b(\v_{\delta_n},\v_{\delta_n},\textbf{g})\right|&\leq Ce^{\sigma t}\left|\mathfrak{z}\right|\|\v_{\delta_n}\|^{\frac{1}{2}}_{\H}\|\nabla\v_{\delta_n}\|^{\frac{3}{2}}_{\H}\|\nabla\textbf{g}\|_{\H}\leq Ce^{\sigma t}\left|\mathfrak{z}\right|\left(1+\|\v_{\delta_n}\|^{\frac{2}{3}}_{\H}\|\nabla\v_{\delta_n}\|^2_{\H}\right),\label{F16}\\
			\left|e^{2\sigma t}\mathfrak{z}\z_{\delta_n}b(\textbf{g},\v_{\delta_n},\textbf{g})\right|&\leq Ce^{2\sigma t}\left|\mathfrak{z}\right|\left|\z_{\delta_n}\right|\|\textbf{g}\|^{\frac{1}{2}}_{\H}\|\nabla\textbf{g}\|^{\frac{3}{2}}_{\H}\|\nabla\v_{\delta_n}\|_{\H}\leq Ce^{2\sigma t}\left|\mathfrak{z}\right|\left|\z_{\delta_n}\right|\left(1+\|\nabla\v_{\delta_n}\|^2_{\H}\right),\label{F17}\\
			\left|e^{\sigma t}\mathfrak{z}b(\v,\v,\textbf{g})\right|&\leq Ce^{\sigma t}\left|\mathfrak{z}\right|\|\v\|^{\frac{1}{2}}_{\H}\|\nabla\v\|^{\frac{3}{2}}_{\H}\|\nabla\textbf{g}\|_{\H}\leq Ce^{\sigma t}\left|\mathfrak{z}\right|\left(1+\|\v\|^{\frac{2}{3}}_{\H}\|\nabla\v\|^2_{\H}\right),\label{F18}\\
			\left|e^{2\sigma t}\mathfrak{z}\y b(\textbf{g},\v,\textbf{g})\right|&\leq Ce^{2\sigma t}\left|\mathfrak{z}\right|\left|\y\right|\|\textbf{g}\|^{\frac{1}{2}}_{\H}\|\nabla\textbf{g}\|^{\frac{3}{2}}_{\H}\|\nabla\v\|_{\H}\leq Ce^{2\sigma t}\left|\mathfrak{z}\right|\left|\y\right|\left(1+\|\nabla\v\|^2_{\H}\right).\label{F19}
		\end{align}
		Combining \eqref{F1}-\eqref{F4} and \eqref{F14}-\eqref{F19}, we obtain
		\begin{align}\label{F20}
			\frac{\d}{\d t}\|\mathfrak{F}\|^2_{\H}&\leq C\times\begin{cases}
				\left|\mathfrak{z}(\vartheta_{t}\omega)\right|Q_3(t),    &\text{ for } r=3 ,\\ \|\mathfrak{F}\|^2_{\H} + \left|\mathfrak{z}(\vartheta_{t}\omega)\right|Q_4(t), &\text{ for } r>3,
			\end{cases}
		\end{align}
		for a.e. $t\in[\mathfrak{s},\mathfrak{s}+T]$, where
		\begin{align*}
			Q_3(t)&=e^{2\sigma t}\left|\mathfrak{z}(\vartheta_{t}\omega)\right|+e^{\sigma t}\bigg[1+\|\nabla\v_{\delta_n}(t)\|^2_{\H}+\|\nabla\v(t)\|^2_{\H}+\|\v_{\delta_n}(t)\|^{\frac{2}{3}}_{\H}\|\nabla\v_{\delta_n}(t)\|^2_{\H}\nonumber\\&\quad+\|\v(t)\|^{\frac{2}{3}}_{\H}\|\nabla\v(t)\|^2_{\H}+\left|e^{\sigma t}\z_{\delta_n}(\vartheta_{t}\omega)\right|^3+\left|e^{\sigma t}\y(\vartheta_{t}\omega)\right|^3+\|\v_{\delta_n}(t)\|^{4}_{\wi\L^{4}} +\|\v(t)\|^{4}_{\wi\L^{4}}\bigg]\nonumber\\&\quad+e^{2\sigma t}\left|\z_{\delta_n}(\vartheta_{t}\omega)\right|\bigg[1+\|\nabla\v_{\delta_n}(t)\|^2_{\H}\bigg]+e^{2\sigma t}\left|\y(\vartheta_{t}\omega)\right|\bigg[1+\|\nabla\v(t)\|^2_{\H}\bigg]
		\end{align*}
		and \begin{align*}Q_4(t)&=e^{2\sigma t}\left|\mathfrak{z}(\vartheta_{t}\omega)\right|+e^{\sigma t}\bigg[1+\|\v_{\delta_n}(t)\|^{\frac{2}{3}}_{\H}\|\nabla\v_{\delta_n}(t)\|^2_{\H}+\|\v(t)\|^{\frac{2}{3}}_{\H}\|\nabla\v(t)\|^2_{\H}+\left|e^{\sigma t}\z_{\delta_n}(\vartheta_{t}\omega)\right|^r\nonumber\\&\quad+\left|e^{\sigma t}\y(\vartheta_{t}\omega)\right|^r+\|\v_{\delta_n}(t)\|^{r+1}_{\wi\L^{r+1}} +\|\v(t)\|^{r+1}_{\wi\L^{r+1}}\bigg]+e^{2\sigma t}\left|\z_{\delta_n}(\vartheta_{t}\omega)\right|\bigg[1+\|\nabla\v_{\delta_n}(t)\|^2_{\H}\bigg]\nonumber\\&\quad+e^{2\sigma t}\left|\y(\vartheta_{t}\omega)\right|\bigg[1+\|\nabla\v(t)\|^2_{\H}\bigg].
		\end{align*}
		Hence, arguing similarly as in Case I, we obtain the desired convergence.
	\end{proof}
	\begin{lemma}\label{Solu_Conver1}
		For $d=2$ with $r>1$, $d=3$ with $r>3$ and $d=r=3$ with $2\beta\mu\geq1$, assume that $\f\in \mathrm{L}^2_{\emph{loc}}(\mathbb{R};\V')$. Suppose that $\{\delta_n\}_{n\in\N}$ be the sequence such that $\delta_n\to0$. Let $\v_{\delta_n}$ and $\v$ be the solutions of \eqref{WZ_SCBF_Add} and \eqref{SCBF_Add} with initial data $\v_{\delta_n,\mathfrak{s}}$ and $\v_{\mathfrak{s}}$, respectively. If $\v_{\delta_n,\mathfrak{s}}\xrightharpoonup{w}\v_{\mathfrak{s}}$ in $\H$ as $n\to\infty$, then for every $\mathfrak{s}\in\R$ and $\omega\in\Omega$,
		\begin{itemize}
			\item [(i)] $\v_{\delta_n}(\xi,\mathfrak{s},\omega,\v_{\delta_n,\mathfrak{s}})\xrightharpoonup{w}\v(\xi,\mathfrak{s},\omega,\v_{\mathfrak{s}})$ in $\H$ for all $\xi\geq \mathfrak{s}$.
			\item [(ii)] $\v_{\delta_n}(\cdot,\mathfrak{s},\omega,\v_{\delta_n,\mathfrak{s}})\xrightharpoonup{w}\v(\cdot,\mathfrak{s},\omega,\v_{\mathfrak{s}})$ in $\mathrm{L}^2((\mathfrak{s},\mathfrak{s}+T);\V)$ for every $T>0$.
			\item [(iii)] $\v_{\delta_n}(\cdot,\mathfrak{s},\omega,\v_{\delta_n,\mathfrak{s}})\xrightharpoonup{w}\v(\cdot,\mathfrak{s},\omega,\v_{\mathfrak{s}})$ in $\mathrm{L}^{r+1}((\mathfrak{s},\mathfrak{s}+T);\widetilde{\L}^{r+1})$ for every $T>0$.
			\item [(iv)] $\v_{\delta_n}(\cdot,\mathfrak{s},\omega,\v_{\delta_n,\mathfrak{s}})\to\v(\cdot,\mathfrak{s},\omega,\v_{\mathfrak{s}})$ in $\mathrm{L}^2((\mathfrak{s},\mathfrak{s}+T);\L^2(\mathcal{O}_k))$ for every $T>0$ and $k>0$,
			where $\mathcal{O}_k=\{x\in\R^d:|x|<k\}$.
		\end{itemize}
	\end{lemma}
	\begin{proof}
		Proof is similar to Lemma 3.5, \cite{GGW}, and hence we omit it here.
	\end{proof}
	The existence of a unique random $\mathfrak{D}$-pullback attractor (denote here by $\mathscr{A}_\delta$) for continuous cocycle $\Phi_\delta$ follows from Theorem \ref{WZ_RA_UB}. The following Lemma demonstrates the uniform compactness of family of random attractors $\mathscr{A}_{\delta}$.
	\begin{lemma}\label{precompact}
		For $d=2$ with $r>1, d=3$ with $r>3$ and $d=r=3$ with $2\beta\mu\geq1$, assume that $\f\in\mathrm{L}^2_{\mathrm{loc}}(\R;\V')$ and satisfies \eqref{forcing7}. Let $\mathfrak{s}\in\R$ and $\omega\in\Omega$ be fixed. If $\delta_n\to0$ as $n\to\infty$ and $\u_n\in\mathscr{A}_{\delta_n}(\mathfrak{s},\omega)$, then the sequence $\{\u_n\}_{n\in\N}$ has a convergent subsequence in $\H$.
	\end{lemma}
	\begin{proof}
		Since, $\delta_n\to0$, it follows from \eqref{ue_add8} that for every $i\in\{1,2\}$ ($i=1$ and $i=2$ for $d=2$ and $d=3$, respectively), $\mathfrak{s}\in\R$ and $\omega\in\Omega$, there exists $\mathcal{N}_1=\mathcal{N}_1(\mathfrak{s},\omega)$ such that for all $n\geq\mathcal{N}_1$,
		\begin{align}\label{PC2}
			\mathcal{R}^{i}_{\delta_n}(\mathfrak{s},\omega)\leq2\mathcal{R}^{i}_0(\mathfrak{s},\omega).
		\end{align}
		It is given that $\u_n\in\mathscr{A}_{\delta_n}(\mathfrak{s},\omega)$ and due to the property of attractors, we have that $\mathscr{A}_{\delta_n}(\mathfrak{s},\omega)$ is a subset of absorbing set, by Lemma \ref{LemmaUe_add1} and \eqref{PC2} we have, for all $n\geq\mathcal{N}_1$,
		\begin{align}\label{PC3}
			\|\u_n\|^2_{\H}\leq2\mathcal{R}^i_{0}(\mathfrak{s},\omega),
		\end{align}
		where, $i=1$ and $i=2$ for $d=2$ and $d=3$, respectively. It is clear that $\{\u_n:n\in\N\}\subset\H$ and therefore, there exists a subsequence (not relabeling) and $\hat{\u}\in\H$ such that 
		\begin{align}\label{PC4}
			\u_n\xrightharpoonup{w}\hat{\u} \ \text{ in }\  \H.
		\end{align}
		Our next aim to prove that $\u_n\to\hat{\u}$  in  $\H$. Since $\u_n\in\mathscr{A}_{\delta_n}(\mathfrak{s},\omega)$, by the invariance property of $\mathscr{A}_{\delta_n}(\mathfrak{s},\omega)$, for all $l\geq1$, there exists $\u_{n,l}\in\mathscr{A}_{\delta_n}(\mathfrak{s}-l,\vartheta_{-l}\omega)$ such that 
		\begin{align}\label{PC5}
			\u_n=\Phi_{\delta_n}(l,\mathfrak{s}-l,\vartheta_{-l}\omega,\u_{n,l})=\u_{\delta_n}(\mathfrak{s},\mathfrak{s}-l,\vartheta_{-\mathfrak{s}}\omega,\u_{n,l}).
		\end{align}
		Since $\u_{n,l}\in\mathscr{A}_{\delta_n}(\mathfrak{s}-l,\vartheta_{-l}\omega)$ and $\mathscr{A}_{\delta_n}(\mathfrak{s}-l,\vartheta_{-l}\omega)\subseteq\mathcal{K}^i_{\delta_n}(\mathfrak{s}-l,\vartheta_{-l}\omega),$ by Lemma \ref{LemmaUe_add1} and \eqref{PC2} we find that for each $l\geq1$ and $n\geq\mathcal{N}_1(\mathfrak{s}-l,\vartheta_{-l}\omega)$,
		\begin{align}\label{PC6}
			\|\u_{n,l}\|^2_{\H}\leq2\mathcal{R}^i_{0}(\mathfrak{s}-l,\vartheta_{-l}\omega).
		\end{align}
		Due to \eqref{WZ_T_add}, we get 
		\begin{align}\label{PC7}
			\v_{\delta_n}(\mathfrak{s},\mathfrak{s}-l,\vartheta_{-\mathfrak{s}}\omega,\v_{n,l})=\u_{\delta_n}(\mathfrak{s},\mathfrak{s}-l,\vartheta_{-\mathfrak{s}}\omega,\u_{n,l})-e^{\sigma\mathfrak{s}}\textbf{g}(x)\z_{\delta_n}(\omega),
		\end{align}
		where 
		\begin{align}\label{PC8}
			\v_{n,l}=\u_{n,l}-e^{\sigma(\mathfrak{s}-l)}\textbf{g}(x)\z_{\delta_n}(\vartheta_{-l}\omega).
		\end{align}
		From, \eqref{PC5} and \eqref{PC7} we obtain 
		\begin{align}\label{PC9}
			\u_n=\v_{\delta_n}(\mathfrak{s},\mathfrak{s}-l,\vartheta_{-\mathfrak{s}}\omega,\v_{n,l})+e^{\sigma\mathfrak{s}}\textbf{g}(x)\z_{\delta_n}(\omega).
		\end{align}  
		Using \eqref{PC6} and \eqref{PC8}, we find that for $n\geq\mathcal{N}_1(\mathfrak{s}-l,\vartheta_{-l}\omega)$,
		\begin{align}\label{PC10}
			\|\v_{n,l}\|^2_{\H}\leq4\mathcal{R}^i_{0}(\mathfrak{s}-l,\vartheta_{-l}\omega)+2e^{2\sigma(\mathfrak{s}-l)}\left|\textbf{g}(x)\z_{\delta_n}(\vartheta_{-l}\omega)\right|^2.
		\end{align}
		Thanks to \eqref{z4} and \eqref{PC10}, we can find an $\mathcal{N}_2=\mathcal{N}_2(\mathfrak{s},\omega,l)\geq\mathcal{N}_1$ such that for every $l\geq1$ and $n\geq \mathcal{N}_2$,
		\begin{align}\label{PC1}
			\|\v_{n,l}\|^2_{\H}\leq4\mathcal{R}^i_{0}(\mathfrak{s}-l,\vartheta_{-l}\omega)+4e^{2\sigma(\mathfrak{s}-l)}\left|\textbf{g}(x)\right|^2(1+\left|\y(\vartheta_{-l}\omega)\right|^2).
		\end{align}
		By \eqref{PC4} and \eqref{PC9} along with \eqref{z4}, we get 	as $n\to\infty$,
		\begin{align}\label{PC11}
			\v_{\delta_n}(\mathfrak{s},\mathfrak{s}-l,\vartheta_{-\mathfrak{s}}\omega,\v_{n,l}) \xrightharpoonup{w}\hat{\v} \text{ in } \H,\ \text{	where }\ \hat{\v}=\hat{\u}-e^{\sigma\mathfrak{s}}\textbf{g}(x)\y(\omega).
		\end{align}
		It follows from \eqref{PC10} that for each $l\geq1$, the sequence $\{\v_{n,l}\}_{n\in\N}\subset\H$, and hence, we get a subsequence (not relabeling) by a diagonal process such that for every $l\geq1$, there exists $\tilde{\v}_l\in\H$ such that
		\begin{align}\label{PC12}
			\v_{n,l}\xrightharpoonup{w}\tilde{\v}_l  \ \text{ in }\ \H \  \ \text{ as }\  \ n\to\infty.
		\end{align}
		It yields from Lemma \ref{Solu_Conver1} and \eqref{PC12} that as $n\to\infty$,
		\begin{align}
			&\v_{\delta_n}(\mathfrak{s},\mathfrak{s}-l,\vartheta_{-\mathfrak{s}}\omega,\v_{n,l})\xrightharpoonup{w}\v(\mathfrak{s},\mathfrak{s}-l,\vartheta_{-\mathfrak{s}}\omega,\tilde{\v}_{l})\ \text{ in }\ \H, \label{PC13}\\
			&\v_{\delta_n}(\cdot,\mathfrak{s}-l,\vartheta_{-\mathfrak{s}}\omega,\v_{n,l})\xrightharpoonup{w}\v(\cdot,\mathfrak{s}-l,\vartheta_{-\mathfrak{s}}\omega,\tilde{\v}_{l}) \ \text{ in } \ \mathrm{L}^2(\mathfrak{s}-l,\mathfrak{s};\V),\label{PC14}\\
			&	\v_{\delta_n}(\cdot,\mathfrak{s}-l,\vartheta_{-\mathfrak{s}}\omega,\v_{n,l})\xrightharpoonup{w}\v(\cdot,\mathfrak{s}-l,\vartheta_{-\mathfrak{s}}\omega,\tilde{\v}_{l}) \ \text{ in } \ \mathrm{L}^{r+1}(\mathfrak{s}-l,\mathfrak{s};\wi\L^{r+1}),	\label{PC15}\\
			&	\v_{\delta_n}(\cdot,\mathfrak{s}-l,\vartheta_{-\mathfrak{s}}\omega,\v_{n,l})\to\v(\cdot,\mathfrak{s}-l,\vartheta_{-\mathfrak{s}}\omega,\tilde{\v}_{l}) \text{ in }  \mathrm{L}^{2}(\mathfrak{s}-l,\mathfrak{s};\L^{2}(\mathcal{O}_k)),	\label{PC16}
		\end{align}
		where $\mathcal{O}_k=\{x\in\R^d:|x|<k\}$. Now, \eqref{PC11} and \eqref{PC13} imply
		\begin{align}\label{PC17}
			\hat{\v}=\v(\mathfrak{s},\mathfrak{s}-l,\vartheta_{-\mathfrak{s}}\omega,\tilde{\v}_{l}).
		\end{align}
		From \eqref{WZ_CSCBF_Add} together with \eqref{b0}, we obtain
		\begin{align}\label{PC18}
			&\frac{\d}{\d t}\|\v_{\delta_n}\|^2_{\H}+2\alpha\|\v_{\delta_n}\|^2_{\H}+2\mu\|\nabla\v_{\delta_n}\|^2_{\H}+2\beta\|\v_{\delta_n}+e^{\sigma t}\textbf{g}\z_{\delta_n}\|^{r+1}_{\wi\L^{r+1}} \nonumber\\&=2b(\v_{\delta_n}+e^{\sigma t}\textbf{g}\z_{\delta_n},\v_{\delta_n},e^{\sigma t}\textbf{g}\z_{\delta_n})+2\beta\left\langle\mathcal{C}(\v_{\delta_n}+e^{\sigma t}\textbf{g}\z_{\delta_n}),e^{\sigma t}\textbf{g}\z_{\delta_n}\right\rangle+2\left\langle\f(\cdot,t),\v_{\delta_n}\right\rangle\nonumber\\&\quad + 2e^{\sigma t}\z_{\delta_n}\left(\left(\ell-\sigma-\alpha\right)\textbf{g}-\mu\A\textbf{g},\v_{\delta_n}\right). 
		\end{align}
		An application of variation of constant formula with $\omega$ being replaced by $\vartheta_{-\mathfrak{s}}\omega$ yields
		\begin{align}\label{PC19}
			&\|\v_{\delta_n}(\mathfrak{s},\mathfrak{s}-l,\vartheta_{-\mathfrak{s}},\v_{n,l})\|^2_{\H}\nonumber\\&=I_1(n,l)+I_2(n,l)+I_3(n,l)+I_4(n,l)+I_5(n,l)+I_6(n,l)+I_7(n,l),
		\end{align} 
		where 
		\begin{align*}
			I_1(n,l)&=e^{-2\alpha l}\|\v_{n,l}\|^2_{\H},\\
			I_2(n,l)&=-2\mu\int_{-l}^{0}e^{2\alpha\xi}\|\nabla\v_{\delta_n}(\xi+\mathfrak{s},\mathfrak{s}-l,\vartheta_{-\mathfrak{s}},\v_{n,l})\|^2_{\H}\d\xi,\\
			I_3(n,l)&=-2\beta\int_{-l}^{0}e^{2\alpha\xi}\|\v_{\delta_n}(\xi+\mathfrak{s},\mathfrak{s}-l,\vartheta_{-\mathfrak{s}},\v_{n,l})+e^{\sigma (\xi+\mathfrak{s})}\textbf{g}\z_{\delta_n}(\vartheta_{\xi}\omega)\|^{r+1}_{\wi\L^{r+1}}\d\xi,\\
			I_4(n,l)&=	2 \int_{-l}^{0} e^{2\alpha\xi}b(\v_{\delta_n}(\xi+\mathfrak{s},\mathfrak{s}-l,\vartheta_{-\mathfrak{s}},\v_{n,l})+e^{\sigma (\xi+\mathfrak{s})}\textbf{g}\z_{\delta_n}(\vartheta_{\xi}\omega),\nonumber\\&\qquad\qquad\qquad\qquad\v_{\delta_n}(\xi+\mathfrak{s},\mathfrak{s}-l,\vartheta_{-\mathfrak{s}},\v_{n,l}),e^{\sigma (\xi+\mathfrak{s})}\textbf{g}\z_{\delta_n}(\vartheta_{\xi}\omega))\d\xi,\\
			I_5(n,l)&=2\beta \int_{-l}^{0} e^{2\alpha\xi}\left\langle\mathcal{C}(\v_{\delta_n}(\xi+\mathfrak{s},\mathfrak{s}-l,\vartheta_{-\mathfrak{s}},\v_{n,l})+e^{\sigma (\xi+\mathfrak{s})}\textbf{g}\z_{\delta_n}(\vartheta_{\xi}\omega)),e^{\sigma (\xi+\mathfrak{s})}\textbf{g}\z_{\delta_n}(\vartheta_{\xi}\omega)\right\rangle\d\xi,\\
			I_6(n,l)&=2 \int_{-l}^{0} e^{2\alpha\xi}\left\langle\f(\cdot,\xi+\mathfrak{s}),\v_{\delta_n}(\xi+\mathfrak{s},\mathfrak{s}-l,\vartheta_{-\mathfrak{s}},\v_{n,l})\right\rangle\d\xi,
		\end{align*}
		and
		\begin{align*}
			I_7(n,l)&=2 \int_{-l}^{0} e^{2\alpha\xi} e^{\sigma (\xi+\mathfrak{s})}\z_{\delta_n}(\vartheta_{\xi}\omega)\big((\ell-\sigma-\alpha)\textbf{g}-\mu\A\textbf{g},\v_{\delta_n}(\xi+\mathfrak{s},\mathfrak{s}-l,\vartheta_{-\mathfrak{s}},\v_{n,l})\big)\d\xi.
		\end{align*}
		Similar to \eqref{PC19}, by \eqref{CSCBF_Add} and \eqref{PC17}, it can also be obtained that 
		\begin{align}\label{PC19*}
			\|\tilde{\v}\|^2_{\H}&=\|\v(\mathfrak{s},\mathfrak{s}-l,\vartheta_{-\mathfrak{s}}\omega,\tilde{\v}_{l})\|^2_{\H} = e^{-2\alpha l}\|\tilde{\v}_{l}\|^2_{\H} -2\mu\int_{-l}^{0}e^{2\alpha \xi}\|\nabla\v(\xi+\mathfrak{s},\mathfrak{s}-l,\vartheta_{-\mathfrak{s}}\omega,\tilde{\v}_{l})\|^2_{\H}\d\xi\nonumber\\&\quad-2\beta\int_{-l}^{0}e^{2\alpha \xi}\|\v(\xi+\mathfrak{s},\mathfrak{s}-l,\vartheta_{-\mathfrak{s}}\omega,\tilde{\v}_{l})+e^{\sigma (\xi+\mathfrak{s})}\textbf{g}\y(\vartheta_{\xi}\omega)\|^{r+1}_{\wi \L^{r+1}}\d\xi\nonumber\\&\quad+2 \int_{-l}^{0} e^{2\alpha\xi}b(\v(\xi+\mathfrak{s},\mathfrak{s}-l,\vartheta_{-\mathfrak{s}},\tilde{\v}_l)+e^{\sigma (\xi+\mathfrak{s})}\textbf{g}\y(\vartheta_{\xi}\omega),\nonumber\\&\qquad\qquad\qquad\qquad\v(\xi+\mathfrak{s},\mathfrak{s}-l,\vartheta_{-\mathfrak{s}},\tilde{\v}_l),e^{\sigma (\xi+\mathfrak{s})}\textbf{g}\y(\vartheta_{\xi}\omega))\d\xi\nonumber\\&\quad+2\beta \int_{-l}^{0} e^{2\alpha\xi}\left\langle\mathcal{C}(\v(\xi+\mathfrak{s},\mathfrak{s}-l,\vartheta_{-\mathfrak{s}},\tilde{\v}_l)+e^{\sigma (\xi+\mathfrak{s})}\textbf{g}\y(\vartheta_{\xi}\omega)),e^{\sigma (\xi+\mathfrak{s})}\textbf{g}\y(\vartheta_{\xi}\omega)\right\rangle\d\xi\nonumber\\&\quad+2 \int_{-l}^{0} e^{2\alpha\xi}\left\langle\f(\cdot,\xi+\mathfrak{s}),\v(\xi+\mathfrak{s},\mathfrak{s}-l,\vartheta_{-\mathfrak{s}},\tilde{\v}_l)\right\rangle\d\xi\nonumber\\&\quad+2 \int_{-l}^{0} e^{2\alpha\xi} e^{\sigma (\xi+\mathfrak{s})}\y(\vartheta_{\xi}\omega)\big((\ell-\sigma-\alpha)\textbf{g}-\mu\A\textbf{g},\v(\xi+\mathfrak{s},\mathfrak{s}-l,\vartheta_{-\mathfrak{s}},\tilde{\v}_l)\big)\d\xi.
		\end{align}
		From \eqref{PC1}, we have
		\begin{align}\label{PC20}
			\limsup_{n\to\infty}I_1(n,l)\leq 4e^{-2\alpha l}\left[\mathcal{R}^i_{0}(\mathfrak{s}-l,\vartheta_{-l}\omega)+e^{2\sigma(\mathfrak{s}-l)}\left|\textbf{g}(x)\right|^2(1+\left|\y(\vartheta_{-l}\omega)\right|^2)\right].
		\end{align}
		Similarly, from \eqref{PC14}, we obtain
		\begin{align}\label{PC21}
			\limsup_{n\to\infty}I_2(n,l)\leq-2\mu\int_{-l}^{0}e^{2\alpha\xi}\|\nabla\v(\xi+\mathfrak{s},\mathfrak{s}-l,\vartheta_{-\mathfrak{s}},\tilde{\v}_l)\|^2_{\H}\d\xi,
		\end{align}
		and
		\begin{align}\label{PC22}
			\lim_{n\to\infty}I_6(n,l)=2 \int_{-l}^{0} e^{2\alpha\xi}\left\langle\f(\cdot,\xi+\mathfrak{s}),\v(\xi+\mathfrak{s},\mathfrak{s}-l,\vartheta_{-\mathfrak{s}},\tilde{\v}_l)\right\rangle\d\xi.
		\end{align}
		By \eqref{z4} and \eqref{PC14}, we have 
		\begin{align}\label{PC23}
			\lim_{n\to\infty}I_7(n,l)=2 \int_{-l}^{0} e^{2\alpha\xi} e^{\sigma (\xi+\mathfrak{s})}\y(\vartheta_{\xi}\omega)\big((\ell-\sigma-\alpha)\textbf{g}-\mu\A\textbf{g},\v(\xi+\mathfrak{s},\mathfrak{s}-l,\vartheta_{-\mathfrak{s}},\tilde{\v}_l)\big)\d\xi.
		\end{align}
		From \eqref{PC15}, we infer that 
		\begin{align}\label{PC24}
			\limsup_{n\to\infty}I_3(n,l)\leq-2\beta\int_{-l}^{0}e^{2\alpha\xi}\|\v(\xi+\mathfrak{s},\mathfrak{s}-l,\vartheta_{-\mathfrak{s}},\tilde{\v}_l)+e^{\sigma (\xi+\mathfrak{s})}\textbf{g}\y(\vartheta_{\xi}\omega)\|^{r+1}_{\wi\L^{r+1}}\d\xi.
		\end{align}
		Using  \eqref{b0}, we obtain
		\begin{align*}
			I_4(n,l)&= 2 \int_{-l}^{0} e^{2\alpha\xi}b(\v_{\delta_n}(\xi+\mathfrak{s},\mathfrak{s}-l,\vartheta_{-\mathfrak{s}},\v_{n,l})+e^{\sigma (\xi+\mathfrak{s})}\textbf{g}\z_{\delta_n}(\vartheta_{\xi}\omega),\nonumber\\&\qquad\qquad\v_{\delta_n}(\xi+\mathfrak{s},\mathfrak{s}-l,\vartheta_{-\mathfrak{s}},\v_{n,l})+e^{\sigma (\xi+\mathfrak{s})}\textbf{g}\z_{\delta_n}(\vartheta_{\xi}\omega),e^{\sigma (\xi+\mathfrak{s})}\textbf{g}\y(\vartheta_{\xi}\omega))\d\xi\nonumber\\&\quad+ 2 \int_{-l}^{0} e^{2\alpha\xi}\left(\z_{\delta_n}(\vartheta_{\xi}\omega)-\y(\vartheta_{\xi}\omega)\right)\cdot b(\v_{\delta_n}(\xi+\mathfrak{s},\mathfrak{s}-l,\vartheta_{-\mathfrak{s}},\v_{n,l})+e^{\sigma (\xi+\mathfrak{s})}\textbf{g}\z_{\delta_n}(\vartheta_{\xi}\omega),\nonumber\\&\qquad\qquad\qquad\v_{\delta_n}(\xi+\mathfrak{s},\mathfrak{s}-l,\vartheta_{-\mathfrak{s}},\v_{n,l})+e^{\sigma (\xi+\mathfrak{s})}\textbf{g}\z_{\delta_n}(\vartheta_{\xi}\omega),e^{\sigma (\xi+\mathfrak{s})}\textbf{g})\d\xi\nonumber\\&=:I_{8}(n,l)+I_9(n,l).
		\end{align*}
		Making use of \eqref{z4} and Lemma \ref{convergence_b}, we get that $\lim\limits_{n\to\infty}I_9(n,l)=0$, and hence, by Lemma \ref{convergence_b}, we obtain
		\begin{align}\label{PC25}
			\lim_{n\to\infty}I_4(n,l)&=	2 \int_{-l}^{0} e^{2\alpha\xi}b(\v(\xi+\mathfrak{s},\mathfrak{s}-l,\vartheta_{-\mathfrak{s}},\tilde{\v}_l)+e^{\sigma (\xi+\mathfrak{s})}\textbf{g}\y(\vartheta_{\xi}\omega),\nonumber\\&\qquad\qquad\qquad\qquad\v(\xi+\mathfrak{s},\mathfrak{s}-l,\vartheta_{-\mathfrak{s}},\tilde{\v}_l),e^{\sigma (\xi+\mathfrak{s})}\textbf{g}\y(\vartheta_{\xi}\omega))\d\xi.
		\end{align}
		Finally, we have 
		\begin{align*}
			I_5(n,l)=2\beta \int_{-l}^{0} e^{2\alpha\xi}\left\langle\mathcal{C}(\v_{\delta_n}(\xi+\mathfrak{s},\mathfrak{s}-l,\vartheta_{-\mathfrak{s}},\v_{n,l})+e^{\sigma (\xi+\mathfrak{s})}\textbf{g}\z_{\delta_n}(\vartheta_{\xi}\omega)),e^{\sigma (\xi+\mathfrak{s})}\textbf{g}\y(\vartheta_{\xi}\omega)\right\rangle\d\xi\nonumber\\ \quad+2\beta \int_{-l}^{0} e^{2\alpha\xi}\left(\z_{\delta_n}(\vartheta_{\xi}\omega)-\y(\vartheta_{\xi}\omega)\right)\big\langle\mathcal{C}(\v_{\delta_n}(\xi+\mathfrak{s},\mathfrak{s}-l,\vartheta_{-\mathfrak{s}},\v_{n,l})+e^{\sigma (\xi+\mathfrak{s})}\textbf{g}\z_{\delta_n}(\vartheta_{\xi}\omega)),\\e^{\sigma (\xi+\mathfrak{s})}\textbf{g}\big\rangle\d\xi=:I_{10}(n,l)+I_{11}(n,l).\qquad\qquad\qquad\qquad
		\end{align*}
		By \eqref{z4} and Lemma \ref{convergence_c}, we get that $\lim\limits_{n\to\infty}I_{11}(n,l)=0$, and hence, by Lemma \ref{convergence_c}, we deduce
		\begin{align}\label{PC26}
			\lim_{n\to\infty}I_5(n,l)&=2\beta \int_{-l}^{0} e^{2\alpha\xi}\left\langle\mathcal{C}(\v(\xi+\mathfrak{s},\mathfrak{s}-l,\vartheta_{-\mathfrak{s}},\tilde{\v}_l)+e^{\sigma (\xi+\mathfrak{s})}\textbf{g}\y(\vartheta_{\xi}\omega)),e^{\sigma (\xi+\mathfrak{s})}\textbf{g}\y(\vartheta_{\xi}\omega)\right\rangle\d\xi.
		\end{align}
		Combining \eqref{PC19}-\eqref{PC26}, we obtain
		\begin{align}\label{PC27}
			&\limsup_{n\to\infty}\|\v_{\delta_n}(\mathfrak{s},\mathfrak{s}-l,\vartheta_{-\mathfrak{s}},\v_{n,l})\|^2_{\H}\nonumber\\&\leq  4e^{-2\alpha l}\left[\mathcal{R}^i_{0}(\mathfrak{s}-l,\vartheta_{-l}\omega)+e^{2\sigma(\mathfrak{s}-l)}\left|\textbf{g}(x)\right|^2(1+\left|\y(\vartheta_{-l}\omega)\right|^2)\right]+ \|\tilde{\v}\|^2_{\H}- e^{-2\alpha l}\|\tilde{\v}_{l}\|^2_{\H}\nonumber\\&\leq  4e^{-2\alpha l}\left[\mathcal{R}^i_{0}(\mathfrak{s}-l,\vartheta_{-l}\omega)+e^{2\sigma(\mathfrak{s}-l)}\left|\textbf{g}(x)\right|^2(1+\left|\y(\vartheta_{-l}\omega)\right|^2)\right]+ \|\tilde{\v}\|^2_{\H}.
		\end{align}
		Passing $l\to\infty$, we get
		\begin{align*}
			\limsup_{n\to\infty}\|\v_{\delta_n}(\mathfrak{s},\mathfrak{s}-l,\vartheta_{-\mathfrak{s}},\v_{n,l})\|^2_{\H}\leq \|\tilde{\v}\|^2_{\H},
		\end{align*}
		which together with \eqref{PC11} gives
		\begin{align}\label{PC28}
			\v_{\delta_n}(\mathfrak{s},\mathfrak{s}-l,\vartheta_{-\mathfrak{s}}\omega,\v_{n,l}) \to\hat{\v} \text{ in } \H.
		\end{align}
		From \eqref{z4}, \eqref{PC9}, \eqref{PC17} and \eqref{PC28}, we get
		\begin{align*}
			\u_n\to \tilde{\u} \text{ in } \H,
		\end{align*}
		as required, and the proof is completed.  
	\end{proof}
	The next Theorem demonstrates the  upper semicontinuity of random $\mathfrak{D}$-pullback attractors as $\delta\to0$, using the abstract theory given in \cite{non-autoUpperWang} (see Theorem 3.2, \cite{non-autoUpperWang}). 
	\begin{theorem}\label{Main_T_add}
		For $0<\delta\leq 1, d=2$ with $r>1, d=3$ with $r>3$ and $d=r=3$ with $2\beta\mu\geq1$, assume that $\f\in\mathrm{L}^2_{\mathrm{loc}}(\R;\V')$ and \eqref{forcing7} is satisfied. Then for every $\omega\in \Omega$ and $\mathfrak{s}\in\R$,
		\begin{align}\label{U-SC}
			\lim_{\delta\to0}\emph{dist}_{\H}\left(\mathscr{A}_{\delta}(\mathfrak{s},\omega),\mathscr{A}_0(\mathfrak{s},\omega)\right)=0.
		\end{align}
	\end{theorem}
	\begin{proof}
		By Lemma \ref{LemmaUe_add1}, we have for every  $\mathfrak{s}\in\R$ and $\omega\in\Omega$,
		\begin{align}\label{U-SC1}
			\limsup_{\delta\to0}\|\mathcal{K}^i_{\delta}(\mathfrak{s},\omega)\|^2_{\H}\leq\limsup_{\delta\to0}\mathcal{R}^i_{\delta}(\mathfrak{s},\omega)=\mathcal{R}^i_0(\mathfrak{s},\omega).
		\end{align}
		Consider a sequence $\delta\to0$ and $\u_{n,\mathfrak{s}}\to\u_{\mathfrak{s}}$ in $\H$. By Lemma \ref{Solu_Conver} we get that for every $t\geq0, \mathfrak{s}\in\R$ and $\omega\in\Omega$,
		\begin{align}\label{U-SC2}
			\Phi_{\delta}(t,\mathfrak{s},\omega,\u_{n,\mathfrak{s}}) \to \Phi_0(t,\mathfrak{s},\omega,\u_{\mathfrak{s}}) \ \text{ in } \ \H.
		\end{align}
		Hence, by \eqref{U-SC1}, \eqref{U-SC2} and Lemma \ref{precompact} together with Theorem 3.2 in \cite{non-autoUpperWang}, one can conclude the proof.
	\end{proof}

	\section{Convergence of attractors: Multiplicative white noise} \label{sec7}\setcounter{equation}{0}
	In this section, we examine the approximations of solutions of the following stochastic CBF equations with multiplicative white noise,
	\begin{equation}\label{SCBF_Multi}
		\left\{
		\begin{aligned}
			\frac{\partial \u}{\partial t}+\mu\A\u+\B(\u)+\alpha\u+\beta\mathcal{C}(\u)&=\boldsymbol{f}+\u\circ\frac{\d \W}{\d t}, \ \ \ \text{ in } \mathbb{R}^n\times(\mathfrak{s},\infty), \\ 
			\u(x,\mathfrak{s})&=\u_{\mathfrak{s}}(x), \ \ \ \ \ \ \ \ \ \ \ \ \ \ x\in \mathbb{R}^n \text{ and }\mathfrak{s}\in\R.
		\end{aligned}
		\right.
	\end{equation}
	For $\delta>0$, consider the pathwise random equations:
	\begin{equation}\label{WZ_SCBF_Multi}
		\left\{
		\begin{aligned}
			\frac{\partial \u_{\delta}}{\partial t}+\mu\A\u_{\delta}+\B(\u_{\delta})+\alpha\u_{\delta}+\beta\mathcal{C}(\u_{\delta})&=\boldsymbol{f}+\mathcal{Z}_{\delta}(\vartheta_{t}\omega)\u_{\delta}, \ \ \text{ in } \mathbb{R}^n\times(\mathfrak{s},\infty), \\ 
			\u_{\delta}(x,\mathfrak{s})&=\u_{\delta,\mathfrak{s}}(x),\ \ \ \  \ \ \  \ \ \ \ \ \ \ x\in \mathbb{R}^n \text{ and }\mathfrak{s}\in\R.
		\end{aligned}
		\right.
	\end{equation}
	The existence of a unique random $\mathfrak{D}$-pullback attractor for \eqref{SCBF_Multi} and \eqref{WZ_SCBF_Multi} is established in \cite{KM6} and Section \ref{sec3}, respectively. In this section, we prove the upper semicontinuity of attractors as $\delta\to0$. Let us denote by $\widehat{\Phi}_0$ and $\widehat{\Phi}_{\delta}$, the continuous cocycles for the systems \eqref{SCBF_Multi} and \eqref{WZ_SCBF_Multi}, respectively, and $\widehat{\mathscr{A}}_0$ and $\widehat{\mathscr{A}}_{\delta}$, random $\mathfrak{D}$-pullback attractors for the systems \eqref{SCBF_Multi} and \eqref{WZ_SCBF_Multi}, respectively. Define 
	\begin{align}
		\v(t,\mathfrak{s},\omega,\v_{\mathfrak{s}})=e^{-\omega(t)}\u(t,\mathfrak{s},\omega,\u_{\mathfrak{s}}) \ \text{	and }\ \v_{\delta}(t,\mathfrak{s},\omega,\v_{\delta,\mathfrak{s}})=e^{-\int_{0}^{t}\mathcal{Z}_{\delta}(\vartheta_{\xi}\omega)\d\xi}\u_{\delta}(t,\mathfrak{s},\omega,\u_{\delta,\mathfrak{s}}).
	\end{align}
	Then, from \eqref{SCBF_Multi} and \eqref{WZ_SCBF_Multi} (formally), we obtain
	\begin{equation}\label{CSCBF_Multi}
		\left\{
		\begin{aligned}
			\frac{\partial \v}{\partial t}+\mu\A\v+e^{\omega(t)}\B(\v)+\alpha\v+\beta e^{(r-1)\omega(t)}\mathcal{C}(\v)&=e^{-\omega(t)}\boldsymbol{f}, \ \ \ \ \  \text{ in } \mathbb{R}^d\times(\mathfrak{s},\infty), \\ 
			\v(x,\mathfrak{s})=\v_{\mathfrak{s}}(x)&=e^{-\omega(\mathfrak{s})}\u_{\mathfrak{s}}(x), \ x\in \mathbb{R}^d \text{ and }\mathfrak{s}\in\R.
		\end{aligned}
		\right.
	\end{equation}
	\begin{equation}\label{WZ_CSCBF_Multi}
		\left\{
		\begin{aligned}
			\frac{\partial \v_\delta}{\partial t}+\mu\A\v_\delta+e^{\int_{0}^{t}\mathcal{Z}_{\delta}(\vartheta_{\xi}\omega)\d\xi}\B(\v_\delta)+\alpha\v_\delta&+\beta e^{(r-1)\int_{0}^{t}\mathcal{Z}_{\delta}(\vartheta_{\xi}\omega)\d\xi}\mathcal{C}(\v_\delta)\\&=e^{-\int_{0}^{t}\mathcal{Z}_{\delta}(\vartheta_{\xi}\omega)\d\xi}\boldsymbol{f}, \ \ \ \ \ \ \text{ in }  \mathbb{R}^d\times(\mathfrak{s},\infty), \\ 
			\v_\delta(x,\mathfrak{s})=\v_{\delta,\mathfrak{s}}(x)&=e^{-\int_{0}^{\mathfrak{s}}\mathcal{Z}_{\delta}(\vartheta_{\xi}\omega)\d\xi}\u_{\delta,\mathfrak{s}}(x), \ x\in \mathbb{R}^d \text{ and }\mathfrak{s}\in\R.
		\end{aligned}
		\right.
	\end{equation}
	For all $\mathfrak{s}\in\R,$ $t>\mathfrak{s},$ and for every initial data in $\H$ and $\omega\in\Omega$, \eqref{CSCBF_Multi} and \eqref{WZ_CSCBF_Multi} have unique solution in $\mathrm{C}([\mathfrak{s},\mathfrak{s}+T];\H)\cap\mathrm{L}^2(\mathfrak{s}, \mathfrak{s}+T;\V)\cap\mathrm{L}^{r+1}(\mathfrak{s},\mathfrak{s}+T;\widetilde{\L}^{r+1}).$ Furthermore, the solution is continuous with respect to initial data and $(\mathscr{F},\mathscr{B}(\H))$-measurable in $\omega\in\Omega.$ 
	\begin{lemma}\label{LemmaUe_Multi}
		For $d=2$ with $r\geq1$, $d=3$ with $r>3$ and $d=r=3$ with $2\beta\mu\geq1$, assume that $\f\in \mathrm{L}^2_{\emph{loc}}(\mathbb{R};\V')$ satisfies \eqref{forcing2}. Then $\widehat{\Phi}_0$ possesses a closed measurable $\mathfrak{D}$-pullback absorbing set $\mathcal{K}_0=\{\mathcal{K}_0(\mathfrak{s},\omega):\mathfrak{s}\in\R, \omega\in\Omega\}\in\mathfrak{D}$ in $\H$ given by
		\begin{align}\label{ue_Multi}
			\mathcal{K}_0(\mathfrak{s},\omega)=\{\u\in\H:\|\u\|^2_{\H}\leq \mathcal{R}_0(\mathfrak{s},\omega)\}, 
		\end{align}
		where $\mathcal{R}_0(\mathfrak{s},\omega)$ is defined by
		\begin{align}\label{ue_Multi1}
			\mathcal{R}_0(\mathfrak{s},\omega)&=\frac{4}{\min\{\mu,\alpha\}} \int_{-\infty}^{0} e^{\alpha\xi-2\omega(\xi)} \|\f(\cdot,\xi+\mathfrak{s})\|^2_{\V'}\d \xi.
		\end{align}
	\end{lemma}
	\begin{proof}
		Since the existence of $\mathfrak{D}$-pullback absorbing set for $\widehat{\Phi}_0$ is proved in \cite{KM6} (see Lemma 5.6, \cite{KM6}), the rest of the proof  follows in a  similar way and hence we omit it here.
	\end{proof}
	\begin{lemma}\label{LemmaUe_Multi1}
		For $d=2$ with $r\geq1$, $d=3$ with $r>3$ and $d=r=3$ with $2\beta\mu\geq1$, assume that $\f\in \mathrm{L}^2_{\emph{loc}}(\mathbb{R};\V')$ satisfies \eqref{forcing2}. Then $\widehat{\Phi}_\delta$ possesses a closed measurable $\mathfrak{D}$-pullback absorbing set $\mathcal{K}_\delta=\{\mathcal{K}_\delta(\mathfrak{s},\omega):\mathfrak{s}\in\R, \omega\in\Omega\}\in\mathfrak{D}$ in $\H$ given by
		\begin{align}\label{ue_Multi4}
			\mathcal{K}_\delta(\mathfrak{s},\omega)=\{\u\in\H:\|\u\|^2_{\H}\leq \mathcal{R}_\delta(\mathfrak{s},\omega)\}, 
		\end{align}
		where $\mathcal{R}_\delta(\mathfrak{s},\omega)$ is defined by
		\begin{align}\label{ue_Multi5}
			\mathcal{R}_\delta(\mathfrak{s},\omega)=\frac{4}{\min\{\mu,\alpha\}} \int_{-\infty}^{0} e^{\int_{0}^{\xi}\left(\alpha-2\mathcal{Z}_{\delta}(\vartheta_{\zeta}\omega)\right)\d\zeta} \|\f(\cdot,\xi+\mathfrak{s})\|^2_{\V'}\d \xi.
		\end{align}
		Furthermore, for every $\mathfrak{s}\in\R$ and $\omega\in\Omega$,
		\begin{align}\label{ue_Multi8}
			\lim_{\delta\to0}\mathcal{R}_\delta(\mathfrak{s},\omega)=\frac{4}{\min\{\mu,\alpha\}} \int_{-\infty}^{0} e^{\alpha\xi-2\omega(\xi)} \|\f(\cdot,\xi+\mathfrak{s})\|^2_{\V'}\d \xi.
		\end{align}
	\end{lemma}
	\begin{proof}
		Since the existence of $\mathfrak{D}$-pullback absorbing set $\mathcal{K}_\delta(\mathfrak{s},\omega)$ for $\widehat{\Phi}_\delta$ is proved in Lemma \ref{PAS} and the convergence in \eqref{ue_Multi8} is proved in \cite{GLW} (see Lemma 3.7 in \cite{GLW}), hence we omit the proof here.
	\end{proof}
	\begin{lemma}\label{Solu_Conver2}
		For $d=2$ with $r\geq1$, $d=3$ with $r>3$ and $d=r=3$ with $2\beta\mu\geq1$, assume that $\f\in \mathrm{L}^2_{\emph{loc}}(\mathbb{R};\V')$. Suppose that $\{\delta_n\}_{n\in\N}$ is a  sequence such that $\delta_n\to0$. Let $\u_{\delta_n}$ and $\u$ be the solutions of \eqref{WZ_SCBF_Multi} and \eqref{SCBF_Multi} with initial data $\u_{\delta_n,\mathfrak{s}}$ and $\u_{\mathfrak{s}}$, respectively. If $\|\u_{\delta_n,\mathfrak{s}}-\u_{\mathfrak{s}}\|_{\H}\to0$ as $n\to\infty$, then for every $\mathfrak{s}\in\R$, $\omega\in\Omega$ and $t>\mathfrak{s}$,
		\begin{align*}
			\|\u_{\delta_n}(t,\mathfrak{s},\omega,\u_{\delta_n,\mathfrak{s}})-\u(t,\mathfrak{s},\omega,\u_{\mathfrak{s}})\|_{\H} \to0  \ \text{ as }\   n\to\infty.
		\end{align*}
	\end{lemma}
	\begin{proof}
		Let $\Upsilon=\v_{\delta_n}-\v$ and $\mathfrak{Z}(t)=\int_{0}^{t}\mathcal{Z}_{\delta}(\vartheta_{\xi}\omega)\d\xi-\omega(t)$. Then from \eqref{WZ_CSCBF_Multi} and \eqref{CSCBF_Multi}, we find 
		\begin{align}\label{P1}
			\frac{\d\Upsilon}{\d t}&=-\mu \A\Upsilon-\alpha\Upsilon-e^{\int_{0}^{t}\mathcal{Z}_{\delta}(\vartheta_{\xi}\omega)\d\xi}\B\big(\v_{\delta_n}\big)+e^{\omega(t)}\B\big(\v\big)-\beta e^{(r-1)\int_{0}^{t}\mathcal{Z}_{\delta}(\vartheta_{\xi}\omega)\d\xi} \mathcal{C}\big(\v_{\delta_n}\big)\nonumber\\&\quad+e^{(r-1)\omega(t)}\beta\mathcal{C}\big(\v\big) +(e^{-\int_{0}^{t}\mathcal{Z}_{\delta}(\vartheta_{\xi}\omega)\d\xi}-e^{-\omega(t)})\f,
		\end{align}
		which gives
		\begin{align}\label{P3}
			&\frac{1}{2}\frac{\d}{\d t}\|\Upsilon\|^2_{\H}\nonumber\\&=-\mu\|\nabla\Upsilon\|^2_{\H}-\alpha\|\Upsilon\|^2_{\H}-e^{\int_{0}^{t}\mathcal{Z}_{\delta}(\vartheta_{\xi}\omega)\d\xi}b(\v_{\delta_n},\v_{\delta_n},\Upsilon) +e^{\omega(t)}b(\v,\v,\Upsilon) \nonumber\\&\quad-e^{(r-1)\int_{0}^{t}\mathcal{Z}_{\delta}(\vartheta_{\xi}\omega)\d\xi}\left\langle\mathcal{C}(\v_{\delta_n}),\Upsilon\right\rangle+e^{(r-1)\omega(t)}\left\langle\mathcal{C}(\v),\Upsilon\right\rangle+(e^{-\int_{0}^{t}\mathcal{Z}_{\delta}(\vartheta_{\xi}\omega)\d\xi}-e^{-\omega(t)})\left\langle\f,\Upsilon\right\rangle\nonumber\\&=-\mu\|\nabla\Upsilon\|^2_{\H}-\alpha\|\Upsilon\|^2_{\H}-e^{\int_{0}^{t}\mathcal{Z}_{\delta}(\vartheta_{\xi}\omega)\d\xi}b(\Upsilon,\v,\Upsilon) -e^{\omega(t)}(e^{\mathfrak{Z}(t)}-1)b(\v,\v,\Upsilon) \nonumber\\&\quad-e^{(r-1)\int_{0}^{t}\mathcal{Z}_{\delta}(\vartheta_{\xi}\omega)\d\xi}\left\langle\mathcal{C}(\v_{\delta_n})-\mathcal{C}(\v),\v_{\delta_n}-\v\right\rangle-e^{(r-1)\omega(t)}(e^{(r-1)\mathfrak{Z}(t)}-1)\left\langle\mathcal{C}(\v),\Upsilon\right\rangle\nonumber\\&\quad+e^{-\omega(t)}(e^{-\mathfrak{Z}(t)}-1)\left\langle\f,\Upsilon\right\rangle.
		\end{align}
		From \eqref{MO_c}, we obtain 
		\begin{align}\label{P8}
			&-e^{(r-1)\int_{0}^{t}\mathcal{Z}_{\delta}(\vartheta_{\xi}\omega)\d\xi}\left\langle\mathcal{C}(\v_{\delta_n})-\mathcal{C}(\v),\v_{\delta_n}-\v\right\rangle\nonumber\\ &\leq-\frac{\beta}{2}e^{(r-1)\int_{0}^{t}\mathcal{Z}_{\delta}(\vartheta_{\xi}\omega)\d\xi}\||\Upsilon||\v_{\delta_n}|^{\frac{r-1}{2}}\|^2_{\H}-\frac{\beta}{2}e^{(r-1)\int_{0}^{t}\mathcal{Z}_{\delta}(\vartheta_{\xi}\omega)\d\xi}\||\Upsilon||\v|^{\frac{r-1}{2}}\|^2_{\H}\leq0.
		\end{align}
		\vskip 2mm
		\noindent
		\textbf{Case I:} \textit{$d=2$ and $r\geq1$.} Applying \eqref{b1}, Lemmas \ref{Holder} and \ref{Young}, we get
		\begin{align}
			\left|e^{\int_{0}^{t}\mathcal{Z}_{\delta}(\vartheta_{\xi}\omega)\d\xi}b(\Upsilon,\v,\Upsilon)\right|&\leq Ce^{\int_{0}^{t}\mathcal{Z}_{\delta}(\vartheta_{\xi}\omega)\d\xi}\|\Upsilon\|_{\H}\|\nabla\Upsilon\|_{\H}\|\nabla\v\|_{\H}\nonumber\\&\leq Ce^{2\int_{0}^{t}\mathcal{Z}_{\delta}(\vartheta_{\xi}\omega)\d\xi}\|\nabla\v\|^2_{\H}\|\Upsilon\|^2_{\H}+\frac{\mu}{4}\|\nabla\Upsilon\|^2_{\H},\label{P4}\\
			\left|e^{\omega(t)}(e^{\mathfrak{Z}(t)}-1)b(\v,\v,\Upsilon)\right|&\leq Ce^{\omega(t)}\left|e^{\mathfrak{Z}(t)}-1\right|\|\v\|^{1/2}_{\H}\|\nabla\v\|^{3/2}_{\H}\|\Upsilon\|^{1/2}_{\H}\|\nabla\Upsilon\|^{1/2}_{\H}\nonumber\\&\leq C\|\Upsilon\|^2_{\H}\|\nabla\Upsilon\|^2_{\H} + Ce^{\frac{4}{3}\omega(t)}\left|e^{\mathfrak{Z}(t)}-1\right|^{4/3} \|\v\|^{2/3}_{\H}\|\nabla\v\|^2_{\H}\nonumber\\&\leq C\|\nabla\v_{\delta_n}\|^2_{\H}\|\Upsilon\|^2_{\H}+C\|\nabla\v\|^2_{\H}\|\Upsilon\|^2_{\H} \nonumber\\&\quad+ Ce^{\frac{4}{3}\omega(t)}|e^{\mathfrak{Z}(t)}-1|^{4/3} \|\v\|^{2/3}_{\H}\|\nabla\v\|^2_{\H},\label{P5}\\
			\left|e^{(r-1)\omega(t)}(e^{(r-1)\mathfrak{Z}(t)}-1)\left\langle\mathcal{C}(\v),\Upsilon\right\rangle\right|&\leq e^{(r-1)\omega(t)}\left|e^{(r-1)\mathfrak{Z}(t)}-1\right|\|\v\|^r_{\widetilde{\L}^{r+1}}\|\Upsilon\|_{\widetilde{\L}^{r+1}}\nonumber\\&\leq C e^{(r-1)\omega(t)}\left|e^{(r-1)\mathfrak{Z}(t)}-1\right|\left[\|\v\|^{r+1}_{\widetilde{\L}^{r+1}}+\|\Upsilon\|^{r+1}_{\widetilde{\L}^{r+1}}\right],\label{P6}\\
			e^{-\omega(t)}(e^{-\mathfrak{Z}(t)}-1)\left\langle\f,\Upsilon\right\rangle&\leq Ce^{-2\omega(t)}\left|e^{-\mathfrak{Z}(t)}-1\right|^2\|\f\|^2_{\V'}+\frac{\min\{\mu,\alpha\}}{4}\|\Upsilon\|^2_{\V}.\label{P7}
		\end{align}
		Combining \eqref{P3}-\eqref{P7}, we get
		\begin{align}\label{P9}
			\frac{\d}{\d t}\|\Upsilon(t)\|^2_{\H}&\leq	P_1(t)\|\Upsilon(t)\|^2_{\H} + P_2(t),
		\end{align}
		for a.e. $t\in[\mathfrak{s},\mathfrak{s}+T]$, where
		\begin{align*}
			P_1(t)&=	C\left[(e^{2\int_{0}^{t}\mathcal{Z}_{\delta}(\vartheta_{\xi}\omega)\d\xi}+1)\|\nabla\v(t)\|^2_{\H}+\|\nabla\v_{\delta_n}(t)\|^2_{\H}\right],\\
			P_2(t)&=Ce^{(r-1)\omega(t)}|e^{(r-1)\mathfrak{Z}(t)}-1|\left[\|\v(t)\|^{r+1}_{\widetilde{\L}^{r+1}}+\|\v_{\delta_n}(t)\|^{r+1}_{\widetilde{\L}^{r+1}}\right]\nonumber\\&\quad+Ce^{\frac{4}{3}\omega(t)}|e^{\mathfrak{Z}(t)}-1|^{4/3} \|\v(t)\|^{2/3}_{\H}\|\nabla\v(t)\|^2_{\H}+Ce^{-2\omega(t)}|e^{-\mathfrak{Z}(t)}-1|^2\|\f(t)\|^2_{\V'}.
		\end{align*}
		\vskip 2mm
		\noindent
		\textbf{Case II:} \textit{$d= 3$ and $r\geq3$ ($r>3$ with any $\beta,\mu>0$ and $r=3$ with $2\beta\mu\geq1$).} An application of Lemmas \ref{Holder} and \ref{Young} yields
		\begin{align}\label{P11}
			\left|e^{\int_{0}^{t}\mathcal{Z}_{\delta}(\vartheta_{\xi}\omega)\d\xi}b(\Upsilon,\v,\Upsilon)\right|\leq\begin{cases}
				\frac{1}{2\beta}\|\nabla\Upsilon\|_{\H}^2+\frac{\beta}{2}e^{2\int_{0}^{t}\mathcal{Z}_{\delta}(\vartheta_{\xi}\omega)\d\xi}\||\Upsilon||\v|\|^2_{\H}, \text{ for } r=3,\\
				\frac{\mu}{4}\|\nabla\Upsilon\|_{\H}^2+\frac{\beta}{4}e^{(r-1)\int_{0}^{t}\mathcal{Z}_{\delta}(\vartheta_{\xi}\omega)\d\xi}\||\Upsilon||\v|^{\frac{r-1}{2}}\|^2_{\H}+C\|\Upsilon\|^2_{\H}, \text{ for } r>3,
			\end{cases}
		\end{align}
		and 
		\begin{align}\label{P12}
			&\left|e^{\omega(t)}(e^{\mathfrak{Z}(t)}-1)b(\v,\v,\Upsilon)\right|\nonumber\\&\leq\left|1-e^{-\mathfrak{Z}(t)}\right|e^{\int_{0}^{t}\mathcal{Z}_{\delta}(\vartheta_{\xi}\omega)\d\xi}\|\nabla\v\|_{\H}\||\Upsilon||\v|\|_{\H}\nonumber\\&\leq\begin{cases}
				\frac{1}{2\beta}\left|1-e^{-\mathfrak{Z}(t)}\right|^2\|\nabla\v\|^2_{\H}+\frac{\beta}{2}e^{2\int_{0}^{t}\mathcal{Z}_{\delta}(\vartheta_{\xi}\omega)\d\xi}\||\Upsilon||\v|\|^2_{\H}, \text{ for } r=3, \\ \frac{1}{2\beta}\left|1-e^{-\mathfrak{Z}(t)}\right|^2\|\nabla\v\|^2_{\H}+\frac{\beta}{4}e^{(r-1)\int_{0}^{t}\mathcal{Z}_{\delta}(\vartheta_{\xi}\omega)\d\xi}\||\Upsilon||\v|^{\frac{r-1}{2}}\|^2_{\H}+C\|\Upsilon\|^2_{\H},  \text{ for } r>3.
			\end{cases}
		\end{align}
		Combining \eqref{P3}-\eqref{P8}, \eqref{P6}-\eqref{P7} and \eqref{P11}-\eqref{P12}, we obtain
		\begin{align}\label{P14}
			\frac{\d}{\d t}\|\Upsilon(t)\|^2_{\H}&\leq	C\|\Upsilon(t)\|^2_{\H} + P(t),
		\end{align}
		for a.e. $t\in[\mathfrak{s},\mathfrak{s}+T]$, where
		\begin{align*}
			P(t)&=Ce^{(r-1)\omega(t)}|e^{(r-1)\mathfrak{Z}(t)}-1|\left[\|\v(t)\|^{r+1}_{\widetilde{\L}^{r+1}}+\|\v_{\delta_n}(t)\|^{r+1}_{\widetilde{\L}^{r+1}}\right]\nonumber\\&\quad+\frac{1}{\beta}\left|1-e^{-\mathfrak{Z}(t)}\right|^2\|\nabla\v\|^2_{\H}+Ce^{-2\omega(t)}|e^{-\mathfrak{Z}(t)}-1|^2\|\f(t)\|^2_{\V'}.
		\end{align*}
		Now, applying Gronwall's inequality to \eqref{P9} and \eqref{P14}, and calculating similarly as in Lemma \ref{Solu_Conver}, one can conclude the proof with the help of the convergence \eqref{N5*}.
	\end{proof}
	The following result shows the uniform compactness of family of random attractors $\widehat{\mathscr{A}}_{\delta}$, which can be obtained by similar calculations of Lemma \ref{precompact}. In fact, the proof is much easier here and hence we omit here.
	\begin{lemma}\label{precompact1}
		For $d=2$ with $r\geq1, d=3$ with $r>3$ and $d=r=3$ with $2\beta\mu\geq1$, assume that $\f\in\mathrm{L}^2_{\mathrm{loc}}(\R;\V')$ and satisfies \eqref{forcing2}. Let $\mathfrak{s}\in\R$ and $\omega\in\Omega$ be fixed. If $\delta_n\to0$ as $n\to\infty$ and $\u_n\in\widehat{\mathscr{A}}_{\delta_n}(\mathfrak{s},\omega)$, then the sequence $\{\u_n\}_{n\in\N}$ has a convergent subsequence in $\H$.
	\end{lemma}
	Finally, we are in the position of giving main result of this section.
	\begin{theorem}\label{Main_T_Multi}
		For $0<\delta\leq 1, d=2$ with $r\geq1, d=3$ with $r>3$ and $d=r=3$ with $2\beta\mu\geq1$, assume that $\f\in\mathrm{L}^2_{\mathrm{loc}}(\R;\V')$ and \eqref{forcing2} is satisfied. Then for every $\omega\in \Omega$ and $\mathfrak{s}\in\R$,
		\begin{align}\label{U-SC_m}
			\lim_{\delta\to0}\emph{dist}_{\H}\left(\widehat{\mathscr{A}}_{\delta}(\mathfrak{s},\omega),\widehat{\mathscr{A}}_0(\mathfrak{s},\omega)\right)=0.
		\end{align}
	\end{theorem}
	\begin{proof}
		By Lemma \ref{LemmaUe_Multi1}, we have for every  $\mathfrak{s}\in\R$ and $\omega\in\Omega$,
		\begin{align}\label{U-SC1_m}
			\limsup_{\delta\to0}\|\mathcal{K}_{\delta}(\mathfrak{s},\omega)\|^2_{\H}\leq\limsup_{\delta\to0}\mathcal{R}_{\delta}(\mathfrak{s},\omega)=\mathcal{R}_0(\mathfrak{s},\omega).
		\end{align}
		Consider a sequence $\delta\to0$ and $\u_{n,\mathfrak{s}}\to\u_{\mathfrak{s}}$ in $\H$. By Lemma \ref{Solu_Conver2}, we get that for every $t\geq0, \mathfrak{s}\in\R$ and $\omega\in\Omega$,
		\begin{align}\label{U-SC2_m}
			\widehat{\Phi}_{\delta}(t,\mathfrak{s},\omega,\u_{n,\mathfrak{s}}) \to \widehat{\Phi}_0(t,\mathfrak{s},\omega,\u_{\mathfrak{s}}) \ \text{ in } \ \H.
		\end{align}
		Hence, by \eqref{U-SC1_m}, \eqref{U-SC2_m} and Lemma \ref{precompact1} together with Theorem 3.2 in \cite{non-autoUpperWang}, we conclude the proof.
	\end{proof}
	\medskip\noindent
	{\bf Acknowledgments:}    The first author would like to thank the Council of Scientific $\&$ Industrial Research (CSIR), India for financial assistance (File No. 09/143(0938)/2019-EMR-I).  M. T. Mohan would  like to thank the Department of Science and Technology (DST), Govt of India for Innovation in Science Pursuit for Inspired Research (INSPIRE) Faculty Award (IFA17-MA110).

	\begin{appendix}
		\renewcommand{\thesection}{\Alph{section}}
		\numberwithin{equation}{section}
		
		\section{Random pullback attractors for Wong-Zakai approximations of 2D stochastic NSE on unbounded Poincar\'e domains} \label{sec5}\setcounter{equation}{0}
		In this appendix, we discuss the existence and uniqueness of random $\mathfrak{D}$-pullback attractors for Wong-Zakai approximations of 2D stochastic NSE on Poincar\'e domains with nonlinear diffusion term. Let $\widetilde{\mathcal{O}}\subset \mathbb{R}^2$, which is open and connected. We also assume that, there exists a positive constant $\lambda_1 $ such that the following Poincar\'e inequality  is satisfied:
		\begin{align}\label{2.1}
			\lambda_1\int_{\widetilde{\mathcal{O}}} |\phi(x)|^2 \d x \leq \int_{\widetilde{\mathcal{O}}} |\nabla \phi(x)|^2 \d x,  \ \text{ for all } \  \phi \in \H^{1}_0 (\widetilde{\mathcal{O}}),
		\end{align}
		that is, $\widetilde{\mathcal{O}}$ is a Poincar\'e domain (may be bounded or unbounded). Consider the Wong-Zakai approximations of 2D NSE on $\widetilde{\mathcal{O}}$ as
		\begin{equation}\label{WZ_NSE}
			\left\{
			\begin{aligned}
				\frac{\partial \u}{\partial t}-\nu \Delta\u+(\u\cdot\nabla)\u+\nabla \p&=\boldsymbol{f}(t) + S(t,x,\u)\mathcal{Z}_{\delta}(\vartheta_t\omega), \ \text{ in } \ \widetilde{\mathcal{O}}\times(\mathfrak{s},\infty), \\ \nabla\cdot\u&=0, \ \text{ in } \ \widetilde{\mathcal{O}}\times(\mathfrak{s},\infty), \\
				\u&=\mathbf{0}\ \ \text{ on } \ \partial\widetilde{\mathcal{O}}\times(\mathfrak{s},\infty), \\
				\u(x,\mathfrak{s})&=\u_{\mathfrak{s}}(x), \ \ \  x\in \widetilde{\mathcal{O}} \text{ and }\mathfrak{s}\in\R,
			\end{aligned}
			\right.
		\end{equation}
		where $\nu$ is the coefficient of kinematic viscosity of the fluid. In the work \cite{GGW}, authors proved the existence of unique random $\mathfrak{D}$-pullback random attractor for the system \eqref{WZ_NSE} under the Assumption \ref{NDT2}. Here, we assume that the following conditions  are satisfied:
		\begin{assumption}\label{NDT4}
			Let $S:\R\times\widetilde{\mathcal{O}}\times\R^2\to\R^2$ be a continuous function such that for all $t\in\R$, $x\in\widetilde{\mathcal{O}}$ and $\textbf{u}\in\widetilde{\mathcal{O}}$
			\begin{align*}
				|S(t,x,\textbf{u})|&\leq \mathcal{S}_1(t,x)|\textbf{u}|^{q-1}+\mathcal{S}_2(t,x),
			\end{align*}
			where $1\leq q<2$, $\mathcal{S}_1\in\mathrm{L}^{\infty}_{\emph{loc}}(\R,\L^{\frac{2}{2-q}}(\widetilde{\mathcal{O}}))$ and $\mathcal{S}_2\in\mathrm{L}^{\infty}_{\emph{loc}}(\R,\L^{2}(\widetilde{\mathcal{O}}))$. Furthermore, suppose that $S(t,x,\textbf{u})$ is locally Lipschitz continuous with respect to $\u$.
		\end{assumption}
		\begin{remark}
			Note that Assumption \ref{NDT4} is clearly different from Assumption \ref{NDT2}. Therefore, it is worth to prove the results under Assumption \ref{NDT4}.
		\end{remark}
		\begin{assumption}\label{DNFT4}
			We assume that the external forcing term $\f\in\mathrm{L}^2_{\mathrm{loc}}(\R,\L^2(\widetilde{\mathcal{O}}))$ satisfies
			\begin{itemize}
				\item [(i)] 
				\begin{align*}
					\int_{-\infty}^{\mathfrak{s}} e^{\nu\lambda_1\xi}\|\f(\cdot,\xi)\|^2_{\L^2(\widetilde{\mathcal{O}})}\d \xi<\infty, \ \ \text{ for all } \mathfrak{s}\in\R.
				\end{align*}
				\item [(ii)] for every $c>0$
				\begin{align*}
					\lim_{\tau\to-\infty}e^{c\tau}\int_{-\infty}^{0} e^{\nu\lambda_1\xi}\|\f(\cdot,\xi+\tau)\|^2_{\L^2(\widetilde{\mathcal{O}})}\d \xi=0.
				\end{align*}
			\end{itemize}
		\end{assumption}
		Since, 2D CBF equations with $r=1$ are the linear perturbation of 2D NSE,  one can prove the next theorem using the same arguments  as it is done for 2D CBF equations in subsection \ref{subsec4.3}. Moreover, there are only minor changes, hence we omit the proof here.  
		\begin{theorem}\label{WZ_RA_UB_GS_NSE}
			Assume that $\f\in\mathrm{L}^2_{\emph{loc}}(\R;\L^2(\widetilde{\mathcal{O}}))$ satisfies Assumption \ref{DNFT4} and Assumption \ref{NDT4} is fulfilled. Then there exists a unique random $\mathfrak{D}$-pullback attractor for the system \eqref{WZ_NSE}, in $\L^2(\widetilde{\mathcal{O}})$.
		\end{theorem}
	\end{appendix}

\end{document}